\newif\ifdeclaration
\declarationtrue
\documentclass[10pt]{amsart}
\usepackage[dvipdfmx]{graphicx}
\usepackage{amsmath,amscd,amsthm,amsxtra,yhmath,stackrel}
\usepackage{epsfig,graphics,color,colortbl}

\usepackage{amssymb,latexsym}
\usepackage{mathrsfs}
\usepackage{bm}
\usepackage[poly,all]{xy}
\usepackage{hyperref}
\usepackage{tikz-cd}
\usepackage[draft]{todonotes}
\usepackage{upgreek}
\usepackage{stmaryrd}
\usepackage{verbatim}
\usepackage{mathtools,xparse}
\usepackage[title]{appendix}
\usepackage{enumitem}
\usepackage[normalem]{ulem}
\usepackage{datetime}
\usepackage[pagewise,mathlines]{lineno}
\usepackage{bbm}
\usepackage[new]{old-arrows}
\usepackage{longtable}

\usepackage{subcaption}
\usepackage[centertableaux]{ytableau}

\usepackage{cleveref}

\newtheorem{thm}{Theorem}[section]
\crefname{thm}{theorem}{theorems}

\newtheorem{prop}[thm]{Proposition}
\crefname{prop}{proposition}{propositions}
\AddToHook{env/prop/begin}{\crefalias{thm}{prop}}

\newtheorem{cor}[thm]{Corollary}
\crefname{cor}{corollary}{corollaries}
\AddToHook{env/cor/begin}{\crefalias{thm}{cor}}

\newtheorem{lem}[thm]{Lemma}
\crefname{lem}{lemma}{lemmas}
\AddToHook{env/lem/begin}{\crefalias{thm}{lem}}

\crefname{conj}{conjecture}{conjectures}
\AddToHook{env/conj/begin}{\crefalias{thm}{conj}}

\theoremstyle{definition}

\newtheorem{df}[thm]{Definition}
\crefname{df}{definition}{definitions}
\AddToHook{env/df/begin}{\crefalias{thm}{df}}

\newtheorem{rem}[thm]{Remark}
\crefname{rem}{remark}{remarks}
\AddToHook{env/rem/begin}{\crefalias{thm}{rem}}

\newtheorem{ex}[thm]{Example}
\crefname{ex}{example}{examples}
\AddToHook{env/ex/begin}{\crefalias{thm}{ex}}

\numberwithin{equation}{section}

\setlist[enumerate,1]{label=\normalfont(\arabic*)}

\newcommand{\mc}{\mathcal}
\newcommand{\mf}{\mathfrak}
\newcommand{\ms}{\mathscr}

\newcommand{\pf}{\noindent{\bfseries Proof. }}
\newcommand{\ov}{\overline}

\newcommand{\cP}{\mathscr{P}}

\newcommand{\N}{\mathbb{N}}
\newcommand{\Z}{\mathbb{Z}}
\newcommand{\Q}{\mathbb{Q}}

\newcommand{\R}{\mathbb{R}}

\newcommand{\e}{\epsilon}
\newcommand{\de}{\delta}

\newcommand{\te}{\tilde{e}}
\newcommand{\tf}{\tilde{f}}
\newcommand{\td}{\widetilde}

\newcommand{\gl}{\mf{gl}}

\newcommand{\g}{\mf{g}}

\newcommand{\La}{\Lambda}

\newcommand{\la}{\lambda}

\newcommand{\ot}{\otimes}

\newcommand{\tbgwed}{{\textstyle\bigwedge}}

\newcommand{\lie}[1]{\mathfrak{#1}}
\DeclareMathOperator{\op}{op}
\DeclareMathOperator{\wt}{wt}
\DeclareMathOperator{\up}{up}
\usepackage{xparse}

\NewDocumentCommand{\QUE}{ m O{} O{q} }{U_{#3}^{#2}(#1)}

\newcommand{\Uqg}{\QUE{\lie{g}}}
\newcommand{\UqgMinus}{\QUE{\lie{g}}[-]}
\newcommand{\UqgZero}{\QUE{\lie{g}}[0]}
\newcommand{\UqgPlus}{\QUE{\lie{g}}[+]}

\newcommand{\Uqgp}{\QUE{\lie{g}, \lie{p}}}
\newcommand{\UqgpPlus}{\QUE{\lie{g}, \lie{p}}[+]}
\newcommand{\UqgpZero}{\QUE{\lie{g}, \lie{p}}[0]}
\newcommand{\UqgpMinus}{\QUE{\lie{g}, \lie{p}}[-]}
\newcommand{\UqgpPlusEqual}{\QUE{\lie{g}, \lie{p}}[\ge 0]}
\newcommand{\UqgpMinusEqual}{\QUE{\lie{g}, \lie{p}}[\le 0]}

\newcommand{\LCrys}{\ms{L}}
\newcommand{\BCrys}{\ms{B}}
\newcommand{\HCrys}[1]{\BCrys\left(#1\right)^{\rm h.w.}}

\newcommand{\LCrysZero}{\LCrys_0}
\newcommand{\BCrysZero}{\BCrys_0}
\newcommand{\HCrysZero}[1]{\BCrysZero\left(#1\right)^{\rm h.w.}}

\newcommand{\FF}{\mathcal{F}}
\newcommand{\MM}{\mathcal{M}}

\newcommand{\FFn}{\FF^{n}}
\newcommand{\BFn}{\BCrys(\FFn)}

\newcommand{\FFnM}{\FFn\otimes \MM}
\newcommand{\BFnM}{\BCrysZero(\FFnM)}

\newcommand{\FFinfty}{\FF^{\infty}}

\newcommand{\FFinftyM}{\FFinfty \otimes \MM}

\newcommand{\FFnplus}{\FF^{n+1}}
\newcommand{\FFnplusM}{\FFnplus \otimes \MM}

\newcommand{\eps}{\epsilon} 
\newcommand{\doteps}{\dot{\eps}}

\DeclareMathOperator{\bprod}{\lozenge}

\DeclareMathOperator{\Image}{Im}
\DeclareMathOperator{\soc}{soc}
\DeclareMathOperator{\Hom}{Hom}

\newcommand{\val}{\mathbbm{v}} 

\newcommand{\Uqgl}{U_q(\gl_\infty)}
\newcommand{\UqslZero}{\QUE{\lie{sl}_{\infty,0}}} 
\newcommand{\UpglPlus}[1][>0]{\QUE{\lie{gl}_{#1}}[][p]}
\newcommand{\UqglPlus}[1][>0]{\QUE{\lie{gl}_{#1}}[]}

\begin{document}
\title
[Tensor product of extremal weight modules of type $A_{+\infty}$]{Infinite-level Fock spaces, crystal bases, and tensor product of extremal weight modules of type $A_{+\infty}$}

\author{JAE-HOON KWON}
\address{Department of Mathematical Sciences and RIM, Seoul National University, Seoul 08826, Korea}
\email{jaehoonkw@snu.ac.kr}

\author{Soo-Hong Lee}
\address{Department of Mathematical Sciences, Seoul National University, Seoul 08826, Korea}
\email{shlee@crystalline.site}

\thanks{This work is supported by the National Research Foundation of Korea(NRF) grant funded by the Korea government(MSIT) (No.2020R1A5A1016126 and RS-2024-00342349).}

\begin{abstract}
We study the category $\mc{C}$ generated by extremal weight modules over $U_q(\gl_{>0})$.
We show that $\mc{C}$ is a tensor category, and {provide} an explicit description of the socle filtration of tensor product of any two extremal weight modules.
This follows from the study of Fock space $\FFinftyM$ of infinite level, which admits commuting actions of a parabolic $q$-boson algebra and $U_p(\gl_{>0})$ with $p=-q^{-1}$.
Its socle has a duality, which can be viewed as a limit of level-rank duality on the fermionic Fock space $\mc{F}^n$ of level $n$.
To describe the socle filtration of $\FFinftyM$, we introduce the notion of a saturated crystal valuation,
whose existence was observed for example in the embedding of an extremal weight module into a tensor product of fundamental weight modules of affine type due to Kashiwara and Beck-Nakajima.
\end{abstract}

\maketitle
\setcounter{tocdepth}{1}
%\tableofcontents

\noindent

\section{Introduction}
\subsection{Extremal weight modules}
Let $U_q(\g)$ be the quantized enveloping algebra associated with a  symmetrizable Kac-Moody algebra $\g$. An extremal weight module $V(\la)$ for an integral weight $\la$ is a $U_q(\g)$-module, which can be viewed as a generalization of highest or lowest weight module. It also has a crystal base and a global crystal basis \cite{Kas94'}. Especially when $\g$ is of affine type and $\la$ is of level zero, $V(\la)$ is isomorphic to a Weyl module introduced in \cite{CP01}, and  they also play an important role in understanding the cell structure of the modified quantum group of level zero \cite{BN,Kas02'}.

\subsection{A tensor category generated by extremal weight $U_q(\gl_{>0})$-modules}
Suppose that $\g=\gl_{>0}$ is a general linear Lie algebra of infinite rank, which is of type $A_{+\infty}$. 
Let $\cP$ be the set of partitions.
Let $V_{\mu,\nu}$ denote the extremal weight module $V(\la)$, where $(\mu,\nu)\in \cP^2$ corresponds to the Weyl group orbit of an integral weight $\la$ for $\gl_{>0}$.
In case of $\gl_{>0}$, the crystal $\ms{B}_{\mu,\nu}$ of $V_{\mu,\nu}$ is connected \cite{K09}, and hence $V_{\mu,\nu}$ is irreducible. 
Let $\mc{C}$ be the category of $U_q(\gl_{>0})$-modules of finite length with irreducible factors $V_{\mu,\nu}$. It is not semisimple, while the subcategories $\mc{C}^\pm$ generated by highest weight modules $V_{\mu,\emptyset}$ and lowest weight modules $V_{\emptyset,\nu}$, respectively are semisimple whose Grothendieck rings $K(\mc{C}^\pm)$ are isomorphic to the ring of symmetric functions.

The main result of this paper is an explicit description of a socle filtration of tensor product of extremal weight modules in $\mc{C}$.
We show that $\mc{C}$ is a tensor category and that the multiplicities of simples in each semisimple subquotient of the socle filtration of $V_{\mu,\nu}\ot V_{\sigma,\tau}$ for $(\mu,\nu),(\sigma,\tau)\in\cP^2$ are given in terms of Littlewood-Richardson coefficients.
In particular, $V_{\mu,\emptyset}\ot V_{\emptyset,\nu}$ is indecomposable with simple socle $V_{\mu,\nu}$, and hence the Grothendieck ring $K(\mc{C})$ is isomorphic to $K(\mc{C}^+)\ot K(\mc{C}^-)$ with two natural $\Z$-bases $\{\,[V_{\mu,\nu}]\,|\,\mu,\nu\in\cP\,\}$ and $\{\,[V_{\sigma,\emptyset}\ot V_{\emptyset,\tau}]\,|\,\sigma,\tau\in\cP\,\}$. We give a character formula of $[V_{\mu,\nu}]\in K(\mc{C})$ in terms of $[V_{\sigma,\emptyset}\ot V_{\emptyset,\tau}]$ with Littlewood-Richardson coefficients.
In particular, $\mc{C}$ provides a categorification of the universal character ring of type $GL$ introduced in \cite{Ko}, in which the classes $[V_{\sigma,\tau}]$ correspond to the basis elements given by universal characters.

We should remark that the decomposition of $\ms{B}_{\mu,\nu}\ot \ms{B}_{\sigma,\tau}$ is given in \cite{K09}, but it does not explain in general the tensor structure on $\mc{C}$. Indeed, the filtration of $V_{\mu,\nu}\ot V_{\sigma,\tau}$ does not always coincide with the decomposition of $\ms{B}_{\mu,\nu}\ot \ms{B}_{\sigma,\tau}$.
For example, we have non-split exact sequences
\begin{equation}\label{eq:example of non-split seq}
\begin{split}
0 \longrightarrow V_{(1),(1)} \longrightarrow  V_{\emptyset,(1)}\otimes V_{(1),\emptyset} \longrightarrow V_{\emptyset,\emptyset} \longrightarrow 0,
\end{split}
\end{equation}
while we have 
$\ms{B}_{(1),(1)} \cong \ms{B}_{\emptyset,(1)}\otimes \ms{B}_{(1),\emptyset}$.
In general, we have $ \ms{B}_{\mu,\nu} \cong \ms{B}_{\emptyset,\nu}\ot \ms{B}_{\mu,\emptyset} \subsetneq \ms{B}_{\mu,\emptyset}\ot \ms{B}_{\emptyset,\nu}$, where the order of tensor product depends on the choice of comultiplcation.
This phenomenon, where a (proper) embedding of modules or crystal lattices induces an isomorphism of crystals, reminds us of the embedding of an extremal weight module into certain tensor product of extremal weight modules associated to multiples of fundamental weights, when $\g$ is of affine type (conjectured in \cite{Kas02'} and proved in \cite{BN}). It is one of our motivations to understand this non-trivial difference between tensor structures in modules and crystals.

There is a non-semisimple tensor category of representations of $\mf{sl}_\infty$ with respect to a non-standrad Borel subalgebra introduced in \cite{PSt}. This category is defined in a complete different way, but it has very similar properties as $\mc{C}$. For example, the irreducible representations are parametrized by $(\mu,\nu)\in \cP^2$, say $\mathbb{V}_{\mu,\nu}$, and $\mathbb{V}_{\mu,\emptyset}\ot \mathbb{V}_{\emptyset,\nu}$ has the same socle filtration as in $\mc{C}$, where the tensor structure is studied through a non-semisimple mixed tensor power of the natural representation and its dual.  More recently, it is shown to have a nice homological property \cite{DPS} and have an interesting application to categorifying the boson-fermion correspondence  \cite{FPS}. It would be interesting to explore more direct connection between these two categories. 

\subsection{A Fock space of infinite level and a limit of level-rank duality}
Let us explain our results in more detail.
Our approach to studying the tensor structure on $\mc{C}$ is to embed the tensor products $V_{\mu,\emptyset}\ot V_{\emptyset,\nu}$ into a limit of the Fock space $\mc{F}^n$ as $n\rightarrow \infty$. Here $\mc{F}^n$ is a $q$-deformed fermionic Fock space of level $n$, which has a $U_q(\gl_{\infty})\ot U_p(\gl_n)$-module structure with $p=-q^{-1}$ (cf.\cite{U}) admitting a $q$-analogue of $(\gl_{\infty},\gl_n)$-Howe duality \cite{Fr} (also known as a level-rank duality).

To have a well-defined limit of $\mc{F}^n$ with a $U_p(\gl_{>0})$-module structure, we introduce a parabolic analogue of $q$-boson algebra $U_q(\mf{sl}_{\infty,0})$ for $\gl_{\infty}$ with respect to its maximal Levi subalgebra (naturally generalizing the $q$-boson algebra introduced in \cite{Kas91}).
It has a family of irreducible representations $V_0(\La_{\mu,\nu})$ parametrized by $(\mu,\nu)\in \cP^2$, which can be identified with maximally parabolic Verma modules of $U_q(\gl_\infty)$ as a $\mathbb{Q}(q)$-space. They form a semisimple category of $U_q(\mf{sl}_{\infty,0})$-modules with crystal bases as in the case of the usual $q$-boson algebra. Our presentation is given with respect to a general pair of $(\g,\mf{p})$ for a symmetrizable Kac-Moody algebra $\g$ and its parabolic subalgebra $\mf{p}$.

Let $\mc{M}=V_0(\La_{\emptyset,\emptyset})$ be the irreducible $U_q(\mf{sl}_{\infty,0})$-modules corresponding to the trivial highest weight.  
By using a $U_q(\gl_\infty)$-comodule structure of $U_q(\mf{sl}_{\infty,0})$, we define a directed system $\{\,\mc{F}^n\ot\mc{M}\,\}_{n\ge 0}$ with a morphism $\phi_{n,n+1} : \mc{F}^n\ot \mc{M}\longrightarrow \mc{F}^{n+1}\ot \mc{M}$ which is $U_q(\mf{sl}_{\infty,0})\ot U_p(\gl_n)$-linear, and let
\begin{equation*}
 \mc{F}^\infty\ot\mc{M} = \varinjlim_{n} \mc{F}^n\ot\mc{M}.
\end{equation*} 
First, we prove the following decomposition:
\begin{equation}\label{eq:Howe duality at infinity non-semisimple}
 \mc{F}^\infty\ot \mc{M} = \bigoplus_{(\mu,\nu)\in\cP^2}V_0(\La_{\mu,\nu})\ot (V_{\mu,\emptyset}\ot V_{\emptyset,\nu}),
\end{equation}
which is a non-semisimple $U_q(\gl_{\infty})\ot U_p(\gl_{>0})$-module (\Cref{thm:isotypic decomposition of F_inftyM}). Using the $U_q(\mf{sl}_{\infty,0})\ot U_p(\gl_n)$-crystal structure of $\mc{F}^n\ot \mc{M}$, we then construct a filtration $\{\,(\mc{F}^\infty\ot\mc{M})_{\ge -d}\,\}_{d\ge 0}$ of $U_q(\mf{sl}_{\infty,0})\ot U_p(\gl_{>0})$-submodules, and show that it has the following semisimple subquotient:
\begin{equation}\label{eq:Howe duality at infinity filtration}
 \frac{(\mc{F}^\infty\ot\mc{M})_{\ge -d}}{(\mc{F}^\infty\ot\mc{M})_{> -d}} \cong \bigoplus_{(\mu,\nu)\in \cP^2}\bigoplus_{\substack{(\zeta,\eta)\in \cP^2 \\ |\mu|-|\zeta|=|\nu|-|\eta|=d}}
 \left(V_0(\La_{\mu,\nu})\ot V_{\zeta,\eta}\right)^{\oplus n_{\zeta,\eta}^{\mu,\nu}}, 
\end{equation}
with $n_{\zeta,\eta}^{\mu,\nu}=\sum_{\sigma}c^{\mu}_{\sigma\zeta}c^{\nu}_{\sigma\eta}$, where $c^{\alpha}_{\beta\gamma}$ is the Littleweeod-Richardson coefficient (\Cref{thm:isotypic decomposition of subquotients}). In particular, when $d=0$, we obtain the following decomposition
\begin{equation}\label{eq:Howe duality at infinity}
(\mc{F}^\infty\ot\mc{M})_{\ge -0}\cong  
\bigoplus_{(\mu,\nu)\in \cP^2} V_0(\La_{\mu,\nu})\ot V_{\mu,\nu}, 
\end{equation}
which can be viewed as a limit of the level-rank duality on the Fock space $\mc{F}^n$. A combinatorial crystal model for \eqref{eq:Howe duality at infinity} (without existence of the associated representation) is given in \cite{K09}.
Note that both $\mc{F}^\infty\ot\mc{M}$ and its proper semisimple submodule $(\mc{F}^\infty\ot\mc{M})_{\ge -0}$ have crystal bases whose crystals are isomorphic.

\subsection{Saturated crystal valuation and socle filtration}
Next we prove that $\{\,(\mc{F}^\infty\ot\mc{M})_{\ge -d}\,\}_{d\ge 0}$ is the socle filtration of $\mc{F}^\infty\ot\mc{M}$ (\Cref{thm:it is the socle filtration}), that is,
\begin{equation}\label{eq:socle filtration of Fock space}
 \mc{V}^\circ_d={\rm soc}(\mc{V}_d),
\end{equation}
where 
\begin{equation*}
    \mc{V}^\circ_d = \frac{(\FFinftyM)_{\ge -d}}{(\FFinftyM)_{>-d}}\ \subset \ \mc{V}_d = \frac{\FFinftyM}{(\FFinftyM)_{>-d}}.
\end{equation*}

For this, we introduce and systematically use the notion of {\em saturated crystal valuation}. It is motivated by an observation that an $A_0$-submodule of a non-semisimple object may induce only a $\mathbb{Q}$-basis of its proper submodule at $q=0$ as in the case of  \eqref{eq:example of non-split seq} or $(\mc{F}^\infty\ot\mc{M})_{\ge -0}\subset \mc{F}^\infty\ot\mc{M}$, where $A_0$ is the subring of $f(q)\in \mathbb{Q}(q)$ regular at $q=0$.
Let $\val$ be a {valuation} on a $\Q(q)$-space $V$, which is equivalent to an $A_0$-submodule $\ms{L}$ of $V$ (not necessarily free) with no nonzero element divisible by $q$ infinitely many times. We call $\ms{L}$ a {crystal valuation} if it is stable under crystal operators and compatible with weight space decomposition when they are available, and say that the crystal valuation $\ms{L}$ is {saturated} with respect to a submodule $V^\circ$
when $\ms{L}$ is a maximal $A_0$-submodule of $V$ that restricts to $\ms{L}^\circ := \ms{L}\cap V^\circ$.
Indeed, this is equivalent to the condition that
the natural inclusion $\ms{L}^\circ \longrightarrow\ms{L}$ induces an isomorphism of $\mathbb{Q}$-spaces at $q=0$. 

We show that there exists a crystal valuation $\val_{-d}^{\infty}$ on $\mc{V}_d$ saturated with respect to $\mc{V}^\circ_d$ (\Cref{thm:valuation limit exists} and \Cref{thm:crystal lattice of a socle become isomorphic mod q}).
The existence of $\val_{-d}^{\infty}$ is obtained by analyzing the behavior of the canonical crystal valuations on the semisimple $U_q(\mf{sl}_{\infty,0})\ot U_p(\gl_n)$-module $\mc{F}^n\ot \mc{M}$ under the morphisms in the directed system and then taking a limit of appropriately shifted valuations on $\mc{F}^n\ot \mc{M}$.
This proof is the technical heart of the paper.

Then we prove \eqref{eq:socle filtration of Fock space} by using the sequence of subquotients associated to a filtration of $\mc{F}^n\ot \mc{M}$ whose limit is  $\mc{V}^\circ_d$ \eqref{eq:Howe duality at infinity filtration}, and the saturatedness of the crystal valuation on $\mc{V}_d$.

Now, it follows by considering the multiplicity space of $V_0(\La_{\mu,\nu})$ in \eqref{eq:Howe duality at infinity non-semisimple} and \eqref{eq:Howe duality at infinity filtration} that
\begin{equation}\label{eq:socle filtration of I}
\frac{{\rm soc}^{d+1}(V_{\mu,\emptyset}\ot V_{\emptyset,\nu})}{{\rm soc}^{d}(V_{\mu,\emptyset}\ot V_{\emptyset,\nu})} 
\cong \bigoplus_{\substack{(\zeta,\eta)\in \cP^2 \\ |\mu|-|\zeta|=|\nu|-|\eta|=d}} V_{\zeta,\eta}^{\oplus n_{\zeta,\eta}^{\mu,\nu}},
\end{equation}
where $n_{\zeta,\eta}^{\mu,\nu}$ is the one given in \eqref{eq:Howe duality at infinity filtration} (\Cref{cor:main application - crystal basis} ).
In particular, this implies that $V_{\mu,\emptyset}\ot V_{\emptyset,\nu}$ is indecomposable. By taking restriction of the saturated crystal valuation on $\mc{V}_d$, we obtain a saturated crystal valuation on $V_{\mu,\emptyset}\ot V_{\emptyset,\nu}/{\rm soc}^{d}(V_{\mu,\emptyset}\ot V_{\emptyset,\nu})$ for $d\in\Z_{\ge 0}$, 
whose existence is a result of its own interest. We also have the same result for $V_{\emptyset,\nu}\ot V_{\mu,\emptyset}$.
As applications of \eqref{eq:socle filtration of I}, we obtain the character formula of $[V_{\mu,\nu}]\in K(\mc{C})$, and an explicit description of the socle filtration of $V_{\mu,\nu}\ot V_{\sigma,\tau}$ for any $\mu,\nu,\sigma,\tau\in \cP$ (\Cref{thm:socle filtration of tensor product}).

\subsection{The organization}
The paper is organized as follows.
In \Cref{sec:preliminary}, we briefly review necessary background.
In \Cref{sec:parabolic q-boson algebra}, we introduce the notion of a parabolic $q$-boson algebra and its integrable representations in a general setting. Then we prove the semisimplicity and existence of a crystal base of an integrable representation.
In \Cref{sec:extremal weight modules of type A}, we recall the crystals of extremal weight modules of $U_q(\gl_{>0})$, and give a filtration of $V_{\mu,\emptyset}\ot V_{\emptyset,\nu}$, which is weaker than the socle filtration, by standard arguments using canonical basis.
In \Cref{sec:Fock space F^n}, we define the Fock space $\mc{F}^n$ as a semi-infinite limit of a $q$-deformed exterior algebra. 
In \Cref{sec:Fock space of infinite level}, we define a Fock space $\mc{F}^\infty\ot \mc{M}$, and prove the decomposition \eqref{eq:Howe duality at infinity non-semisimple}. We also define a filtration on $\mc{F}^\infty\ot \mc{M}$.
In \Cref{sec:crystal valuation}, we introduce the notion of a (saturated) crystal valuation, and prove the decomposition of a crystal valuation with respect to an isotypic decomposition of an integrable representation.
In \Cref{sec:Crystal valuation on FFinftyM}, 
we prove the decomposition of the subquotients \eqref{eq:Howe duality at infinity filtration}. Finally, we prove the existence of a saturated crystal valuation on $\mc{V}_d$ by which we prove that the filtration on $\mc{F}^\infty\ot \mc{M}$ is the socle filtration \eqref{eq:socle filtration of Fock space}.

\section{Preliminary}
\label{sec:preliminary}

\subsection{Quantized enveloping algebras}\label{subsec:quantum group}
Let $A =(a_{ij})_{i,j\in I}$ be a symmetrizable generalized Cartan matrix indexed by a set $I$, possibly infinite.
Let $P^\vee$ be the dual weight lattice, 
and $\Pi^\vee=\{\,h_i\,|\,i\in I\,\}\subset P^\vee$ the set of simple coroots.
We assume that $P^\vee / \Z\Pi^\vee$ has a finite rank.
Let $P = \{\, f\in \Hom_\Z(P^\vee,\Z) \, | \, f(h_i) = 0\text{ for all but finitely many } i\in I \, \}$, the restricted dual of $P^\vee$ with respect to $\Pi^\vee$, be the weight lattice,
and $\Pi=\{\,\alpha_i\,|\,i\in I\,\}\subset P$ the set of simple roots, which are linearly independent,
such that $\langle h_i,\alpha_j\rangle=a_{ij}$ for $i,j\in I$.

Let $(\,,\,)$ be a symmetric bilinear form on $P$ such that $(\alpha_i,\alpha_i)\in 2\Z_{+}$ for $i\in I$, $(\alpha_i,\alpha_j)\le 0$ for $i\neq j$, and $\langle h_i,\la\rangle=2(\alpha_i,\la)/(\alpha_i,\alpha_i)$ for $i\in I$ and $\la\in P$.
Let $\g$ be the Kac-Moody algebra associated with the Cartan datum $(A, P, P^\vee, \Pi, \Pi^\vee, (,))$.
Let $W$ be the Weyl group of $\g$ generated by the simple reflection $s_i$, where $s_i(\lambda)=\lambda-\langle h_i,\la \rangle\alpha_i$ for $\la\in P$.
Let $Q=\bigoplus_{i\in I}\Z\alpha_i$, and $Q_{\pm}=\pm\bigoplus_{i\in I}\Z_{\ge 0}\alpha_i$, and let $\ge$ denote the usual partial order on $P$.
Let $P^+=\{\,\la\in P\,|\,\langle h_i,\la\rangle\in \Z_{\ge 0}\ (i\in I)\,\}$ be the set of dominant integral weights.

Let $q$ be an indeterminate.
For $i\in I$, put $q_i=q^{(\alpha_i,\alpha_i)/2}$. For $a\in\Z_{\ge 0}$ and $i\in I$, let $[a]_i=\frac{q_i^a-q_i^{-a}}{q_i-q_i^{-1}}$ and $[a]_i!=[a]_i[a-1]_i\dots [1]_i$ $(a\ge 1)$ with $[0]_i=1$. If $A$ is symmetric and $(\alpha_i,\alpha_i)=2$ for all $i\in I$, then we simply write $[a]_i=[a]$.

Let $U_q(\g)$ be the associated quantized enveloping algebra, which is an associative $\Q(q)$-algebra generated by $e_i, f_i, q^{h}$ for $i\in I$ and $h\in P^\vee$ subject to the following relations:
{\allowdisplaybreaks
\begin{gather*}
q^0=1, \quad q^{h + h'}=q^{h}q^{h'},\\
 q^h e_i q^{-h}=q^{\langle h,\alpha_i\rangle}e_i,\quad
 q^h f_i q^{-h}=q^{-\langle h,\alpha_i\rangle}f_i, \\
 e_if_j - f_je_i =
\delta_{ij}\frac{t_{i} - t^{-1}_{i}}{q_i-q_i^{-1}},\\
\sum_{k=0}^{c_{ij}} (-1)^k e_i^{(k)} e_j e_i^{(c_{ij}-k)}
= \sum_{k=0}^{c_{ij}} (-1)^k f_i^{(k)} f_j f_i^{(c_{ij}-k)} =0, \label{eq:Serre-rel}
\end{gather*}
\noindent where $t_i=q^{(\alpha_i,\alpha_i)h_i/2}$, $e_i^{(k)}=e_i^k/[k]_i!$, $f_i^{(k)}=f_i^k/[k]_i!$, and $c_{ij}=1-a_{ij}$ for $i,j\in I$.}
Unless otherwise specified, we regard $U_q(\g)$ as a Hopf algebra with respect to the comultiplication $\Delta$ and the antipode $S$ given by
\begin{equation*}
\begin{split}
& \Delta(q^h)=q^h\otimes q^h, \\
& \Delta(e_i)= 1\ot e_i + e_i\ot t_i^{-1}, \\
& \Delta(f_i)= f_i\ot 1 + t_i\ot f_i , \\
S(q^h)= &\  q^{-h}, \quad S(e_i)=-e_i t_i, \quad S(f_i)=-t_i^{-1} f_i,
\end{split}
\end{equation*}
for $i\in I$ and $h\in P^\vee$,
where $\Delta = \Delta_-$ is often called the lower comultiplication.
Let $\Delta_+$ be a comultiplication of $U_q(\g)$ by
\begin{equation*}
\begin{split}
& \Delta_+(q^h)=q^h\otimes q^h, \\
& \Delta_+(e_i)= e_i\ot 1 + t_i\ot e_i , \\
& \Delta_+(f_i)= 1\ot f_i + f_i\ot t_i^{-1},
\end{split}
\end{equation*}
which is called upper comultiplication.
Let $\Delta^{\op}_\pm := \sigma \circ \Delta_\pm$ be the comultiplications
where $\sigma$ is a map $\sigma(x\ot y) = y\ot x$.

We denote by $U_q^{\pm}(\g)$ the subalgebra generated by $e_i$ and $f_i$ ($i\in I$), respectively, which is graded by $Q_{\pm}$, and denote by $U_q^0(\g)$ the subalgebra generated by $q^h$ ($h\in P^\vee$). We put $U_q^{\le 0}(\g)=U^{-}_q(\g)U^0_q(\g)$ and $U_q^{\ge 0}(\g)=U^{+}_q(\g)U^0_q(\g)$.

Let $A=\Z[q,q^{-1}]$, and let $U_q(\mf{g})_{A}$ be the $A$-subalgebra of $U_q(\mf{g})$ generated by $e_i^{(k)}$ and $f_i^{(k)}$ for $0\le i\le n$ and $k\in \Z_{\ge 0}$.

Let $\tau_\pm: \Uqg \to \Uqg$ be $\Q(q)$-linear anti-automorphisms defined by
\begin{equation}\label{eq:tau_pm for quantum groups}
    \tau_\pm(e_i) = q_i^{\pm 1}t_i^{\pm 1} f_i, \quad \tau_\pm(f_i) = q_i^{\pm 1}t_i^{\mp 1} e_i, \quad \tau_\pm(q^h) = q^h, \quad \text{for $i\in I$ and $h\in P^\vee$.}
\end{equation}
Then we have $(\tau_\pm \ot \tau_\pm) \circ \Delta_\pm = \Delta_\pm \circ \tau_\pm$.

\subsection{Crystal bases}\label{subsec:crystal base}

Let us briefly recall the notion of crystal base and its properties \cite{Kas91,Kas94',Kas02}.
Let $V$ be an integrable $U_q(\g)$-module. In other words, $V$ has a weight space decomposition $V=\bigoplus_{\mu\in P}V_\mu$, where $V_{\mu}=\{\,v\,|\,q^hv=q^{\langle h, \mu \rangle}v \ (h\in P^\vee)\,\}$, and $e_i, f_i \ (i\in I)$ act locally nilpotently on $V$.
We write $\wt(v) = \mu$ for $v\in V_\mu$.

Let $i\in I$ be given. For a weight vector $v\in V$,
we may write $v=\sum_{k \geq 0} f_i^{(k)}v_k$, where each $v_k$ satisfies $e_iv_k=0$. The (lower) crystal operators $\te^{\rm low}_i, \tf^{\rm low}_i$  or simply $\te_i$, $\tf_i$, are defined by
\begin{equation}\label{eq:Kashiwara operator lower}
\te_iv=\sum_{k\geq1}f_i^{(k-1)}v_k,\quad \tf_iv=\sum_{k\geq0}f_i^{(k+1)}v_k.
\end{equation}
and the upper crystal operators $\te^{\rm\, up}_i, \tf^{\rm\, up}_i$ are given by
\begin{equation}\label{eq:Kashiwara operator upper}
\te^{\rm\, up}_iv=\sum_{k\geq1}q_i^{-l_k+2k-1}f_i^{(k-1)}v_k,\quad
\tf^{\rm\, up}_iv=\sum_{k\geq 0}q_i^{l_k-2k-1}f_i^{(k+1)}v_k,
\end{equation}
where $l_k=\langle h_i,{\rm wt}(v_k)\rangle$.

Let $A_0$ be the subring of $\Q(q)$ consisting of $f(q)$ regular at $q=0$. A lower crystal base of $V$ is a pair $(L,B)$, where $L$ is an $A_0$-lattice of $V$, and $B$ is a $\mathbb{Q}$-basis of $L/qL$ satisfying

\begin{itemize}
\item[(1)] $L=\bigoplus_{\mu\in P} L_\mu$ and
$B=\bigsqcup_{\mu\in P} B_\mu $, where $L_\mu = L\cap V_\mu$ and $B_\mu = B\cap (L/qL)_\mu$,

\item[(2)] $\te_i L\subset L$, $\tf_iL\subset L$ and $
\te_i B \subset B\cup\{0\}$, $\tf_i B \subset B\cup\{0\}$ for $i\in I$,

\item[(3)] $\tf_ib=b'$ if and only if $\te_ib'= b$ for $i\in I$ and $b, b'\in B$,
\end{itemize}
while an upper crystal base of $V$ is defined with respect to \eqref{eq:Kashiwara operator upper}.
We call $B$ a crystal of $V$. We call an $A_0$-lattice $L$ of $V$ satisfying (1) and (2) a crystal lattice of $V$.

For $i\in I$, let $e'_i, e''_i :\UqgMinus \longrightarrow \UqgMinus$ be the $\Q(q)$-linear maps given by
\begin{equation} \label{eq:Kashiwara operators}
    [e_i, u] = \frac{t_ie''_i(u) - t_i^{-1}e'_i(u)}{q_i - q_i^{-1}},
\end{equation}
for homogeneous $u\in \UqgMinus$.
Note that ${}_ir(u) = e_i''(u)$ and $r_i(u) = q_i^{\langle h_i, \wt u\rangle + 2}e_i'(u)$, where ${}_ir,r_i$ are given in \cite{Lu93}.
We may also write $u=\sum_{k \geq 0} f_i^{(k)}v_k$, where each $v_k$ satisfies $e'_iv_k=0$. Then the crystal operators on $\UqgMinus$ are defined as in \eqref{eq:Kashiwara operator lower} and a crystal base of $\UqgMinus$ is defined in the same way. Then $\UqgMinus$ has a unique crystal base $(\LCrys(\infty), \BCrys(\infty))$.

Let $V_i$ be integrable $U_q(\g)$-modules with lower or upper crystal bases $(L_i,B_i)$ $(i=1,2)$. Then the tensor product rule states that $(L_1\ot L_2, B_1\ot B_2)$ is a crystal base of $V_1\ot V_2$ such that
\begin{equation}\label{eq:tensor product rule}
\begin{split}
{\te}_i(b_1\otimes b_2)=
\begin{cases}
{\te}_i b_1 \otimes b_2, & \text{if  $\varphi_i(b_1)\ge\varepsilon_i(b_2)$}, \\
b_1\otimes {\te}_i b_2, & \text{if $\varphi_i(b_1)<
\varepsilon_i(b_2)$},
\end{cases}\\
{\tf}_i(b_1\otimes b_2)=
\begin{cases}
{\tf}_i b_1 \otimes b_2, & \text{if  $\varphi_i(b_1)>\varepsilon_i(b_2)$}, \\
b_1\otimes {\tf}_i b_2, & \text{if $\varphi_i(b_1)\leq
\varepsilon_i(b_2)$},
\end{cases}
\end{split}
\end{equation}
for $i\in I$ and $b_1\ot b_2\in B_1\ot B_2$, where $\varepsilon_i(b)=\max\{k\geq 0 \,|\ \te_i^k b  \neq 0 \}$ and $\varphi_i(b)=\max\{k\geq 0 \,|\ \tf_i^k b \neq 0 \}$ for $b\in B_1, B_2$.

A weight vector $v\in V_\la$ is called $i$-extremal ($i\in I$) if $e_iv = 0$ or $f_iv=0$. If $v$ is $i$-extremal, we define
\begin{equation}\label{eq:extremal operators}
    S_i v =
\begin{cases}
    f_i^{(\langle h_i, \la\rangle)} v & \text{if } e_i v = 0 \\
    e_i^{(-\langle h_i, \la\rangle)} v & \text{if } f_i v = 0
\end{cases}
\end{equation}
A weight vector $v\in V$ is called an extremal vector if there exists $\{v_w\}_{w\in W}$ such that $v_w$s are $i$-extremal for all $i\in I$, and
\begin{equation*}
    v_e = v, \quad S_i v_w = v_{s_i w},
\end{equation*}
for $ i\in I, w\in W$.
An element of a crystal is also called extremal if it satisfies the same condition, where $e_i^{(n)}$ and $f_i^{(n)}$ are replaced by $\te_i^n$ and $\tf_i^n$, respectively.

For $\la\in P$, let $V(\la)$ be the $U_q(\g)$-module generated by $u_\la$ subject to the relations that $u_\la$ is an extremal vector of weight $\la$ \cite{Kas94'}. Note that $V(\la)$ is a highest weight module if $\la\in P^+$. It is shown in \cite{Kas94'} that $V(\la)$ has a crystal base $(\LCrys(\la),\BCrys(\la))$ and a global crystal basis or canonical basis $G(\la)=\{\,G_\la(b)\,|\,b\in \ms{B}(\la)\,\}$.
We often assume that $u_\la \in \BCrys(\la)_\la \pmod{q\LCrys(\la)}$.

Suppose that $\la\in P^+$. Then the canonical projection $\pi_\la: \UqgMinus \to V(\la)$ restricts to $\pi_\la|_{\LCrys(\infty)}: \LCrys(\infty) \to \LCrys(\la)$, and induces $\ov{\pi}_\la: \BCrys(\infty) \to \BCrys(\la) \sqcup \{0\}$.
There is a unique bilinear form $\langle\cdot,\cdot\rangle$ on $V(\la)$, called the $q$-Shapovalov form, which is characterized by
\begin{equation}\label{eq:q-Shapovalov form on uqg}
    \langle u_\lambda, u_\lambda\rangle = 1, \quad \langle uv,w\rangle = \langle v, \tau_-(u)w\rangle,
\end{equation}
for $u\in \Uqg, v,w\in V(\la)$.
We have $\langle \LCrys(\lambda), \LCrys(\lambda)\rangle \in A_0$, and $\LCrys(\lambda) = \{v\in V(\la) \,|\, \langle v, \LCrys(\lambda)\rangle \in A_0\}$ (cf. \cite{Kas91}).
Also recall that, for  $b\in \BCrys(\lambda)$, and $m \ge 0$, we have
\begin{equation}\label{eq:approximating action of divided powers}
    f_i^{(m)}G_\lambda(b) \in \bigoplus_{b'\in \BCrys(\lambda)_{\wt b -m\alpha_i}} q_i^{-m(\varepsilon_i(b') + m)}A_0G_\lambda(b').
\end{equation}

Let $V$ be an integrable $\Uqg$-module with a crystal base $(\LCrys(V), \BCrys(V))$.
Let $\HCrys{V} = \{ b\in \BCrys(V) \,|\, \te_i b = 0 \text{ for all } i\in I\}$, and for $\nu \in P$, let $\HCrys{V}_{\ge \nu} = \{ b\in \HCrys{V} \,|\, \wt(b) \ge \nu\}$.
The following lemma can be easily proved, and its analogues with respect to crystal operators for other algebras will be frequently used later.

\begin{lem}\label{lem:highest weight crystal representatives generate filtration}
    Suppose that $V$ has finite dimensional weight spaces and $\wt(V)$, the set of weights of $V$, is finitely dominated.
    For each $b\in \HCrys{V}$, choose a weight vector $x_b\in \LCrys(V)$ such that $x_b \equiv b\pmod{q\LCrys(V)}$.
    Then 
    \begin{equation}\label{eq:elementary_lemma}
    \begin{split}
        \Uqg\text{-span of } \{ x_b \,|\, b\in \HCrys{V}_{\ge\nu}\} &= \UqgMinus\text{-span of } \{ x_b \,|\, b\in \HCrys{V}_{\ge\nu}\}  \\
        & \cong \bigoplus_{b\in \HCrys{V}_{\ge\nu}}V(\wt(b)),
    \end{split}
    \end{equation}
    and $\{ x_b \,|\, b\in \HCrys{V}_{\ge\nu}\}$ generates a crystal lattice of $\bigoplus_{b\in \HCrys{V}_{\ge\nu}}V(\wt(b))$ under $\tf_i$ for $i\in I$.
\end{lem}
Indeed, we may assume that $V = \bigoplus_{\nu\in P^+}V(\nu)^{\oplus m_\nu}$
and $\LCrys(V) = \bigoplus_{\nu\in P^+}\LCrys(\nu)^{\oplus m_\nu}$.
Let $V_1$ (resp. $V_2$) denotes the space on LHS (resp. RHS) in the first line of \eqref{eq:elementary_lemma}, and $V_3$ denotes the space in the second line of \eqref{eq:elementary_lemma}.
It is immediate that $V_2 \subset V_1 \subset V_3$.
To see $V_3 \subset V_2$, one proceeds with induction on $\nu$, noting that the projection of $\{\,x_b\,|\,b\in\HCrys{V}_{\nu}\,\}$ to $V(\nu)^{\oplus m_\nu}$ gives a $A_0$-basis of $\LCrys(\nu)^{\oplus m_\nu}_{\nu}$,
which generates $V(\nu)^{\oplus m_\nu}$ (resp. $\LCrys(\nu)^{\oplus m_\nu}$) under the action of $U_q^-(\g)$ (resp. $\tf_i$'s).

\subsection{{Quasi-$R$-matrix, $R$-matrix, and canonical basis}}\label{subsec:quasi R}
Let us briefly review necessary materials on $R$-matrix and canonical bases of based modules \cite{Lu93}.

Let $-$ be the involution of $\Q$-algebras on $U_q(\g)$ given by $\ov{e_i}=e_i$, $\ov{f_i}=f_i$, $\ov{q^h}=q^{-h}$, and $\ov{q}=q^{-1}$ for $i\in I$ and $h\in P^\vee$. Let $M$ and $N$ be $\Uqg$-modules with weight space decomposition. 

Let $\Theta$ be the quasi-$R$-matrix \cite[Theorem 4.1.2]{Lu93}, which is given by
\begin{equation} \label{eq:quasi R matrix}
 \Theta = \Theta_+ = \sum_{\beta\in Q_+}\Theta_\beta
\end{equation}
such that $\Theta_\beta$ is a unique element in $U^-_q(\g)_{-\beta}\ot U^+_q(\g)_{\beta}$ such that $\Theta_0=1\ot 1$ and ${\Delta_+}(u)\, \Theta = \Theta\, \ov{\Delta}_+(u)$
for $u\in U_q(\g)$, where $\ov{\Delta_+}(u)=\ov{\Delta_+(\ov{u})}$ for $u\in U_q(\g)$.
Here, $\Theta$ is regarded as an element in a suitable completion $\UqgMinus\widehat{\ot}\ \UqgPlus$. We remark that the completion used in \cite{Lu93} is not suitable for $\g$ with $I$ infinite.
We instead consider a finer completion given by subspaces $\UqgPlus \UqgZero\sum_{\nu' \le \nu}\UqgMinus_{-\nu'} \otimes \Uqg + \Uqg\otimes \UqgMinus\UqgZero\sum_{\nu' \le \nu} \UqgPlus_{\nu'}$ parametrized by $\nu\in Q_+$, where $Q_+$ is regarded as a directed set by its poset structure. Then the statements for $\Theta$ in \cite{Lu93} holds for infinite $I$ with respect to this completion.

\newcommand{\UnivRMat}{R^{\rm univ}}

The quasi-$R$-matrix $\Theta$ yields a universal $R$-matrix $\UnivRMat = R =  \sigma\Pi\ov{\Theta}$,
where $\Pi$ is a $\Q(q)$-linear operator acting on $M\ot N$ by $\Pi (m\ot n) = q^{(\mu,\nu)}m\ot n$, and $\sigma$ is given by $\sigma(m\ot n)=n\ot m$ for $m\in M_\mu$ and $n \in N_\nu$.
Note that we need to extend the base field to $\Q(q^{\frac{1}{d}})$ for $d\in \Z_{>0}$ in general, but we can take $d = 1$ in type $A$.

We also need an opposite version of the quasi-$R$-matrix, given by
\begin{equation}\label{eq:opposite quasi-$R$-matrix}
    \Theta' = R\Pi\sigma = \Pi \ov{\Theta^{\op}} \Pi,
\end{equation}
where $\Theta^{\op}$ is obtained by applying $\sigma$ to $\Theta$.
Since $R^{-1} = \sigma \Pi \ov{\Theta'}$, this can be viewed as a quasi-$R$-matrix constructed out of $R^{-1}$ instead of $R$.
The following identity follows  from $R\Theta = R^{-1}\Theta' = \sigma\Pi$.
\begin{lem} \label{lem:RThetaRTheta'}
    We have $R\Theta\ov{R\Theta} = R\Theta'\ov{R\Theta'} = 1$.
\end{lem}

\begin{rem}\label{rem:$R$ matrix}
We may also construct $\Theta,\Theta'$, and $R$ for other comultiplications.
We let $\ov{\Delta}_\pm = \Delta^{\op}_\mp$ the coproducts twisted by bar-involutions $\ov{\Delta}_\pm = {\--} \circ \Delta_\pm \circ {\--}$.
Note that we have $\ov{\Delta}_\pm = \Delta^{\op}_\mp$.
Let $M \ot_\pm N$ and $M \ov{\ot}_\pm N$ be the $\Uqg$-module whose module structure is given by $\Delta_\pm$ and $\ov{\Delta}_\pm$, repectively.
These coproducts are related by natural isomorphisms below
\[
    \begin{tikzcd}
        M\ot_+ N \arrow[d, "\sigma"'] \arrow[r, "\Pi^{-1}"] & M\ot_- N \arrow[d, "\sigma"] \\
        N\ov{\ot}_- M \arrow[r, "\Pi^{-1}"']                & N\ov{\ot}_+ M
        \end{tikzcd}
\]
(cf. \cite[(2.2.9)]{KMPY}).
By pulling back $R$, $\Theta$, and $\Theta'$ along the natural isomorphisms above, we obtain 
$R_{\ot_\ast}, \Theta_{\ot_\ast}, \Theta'_{\ot_\ast}$ corresponding to each tensor products $\ot_\ast = \ot_\pm, \ov{\ot}_\pm$.
For instance, we have
\begin{equation}\label{eq:R-matrix for lower comultiplications}
    R_- := R_{\otimes_-} = \sigma \ov\Theta \Pi, \quad  \Theta_-:= \Theta_{\otimes_-} = \Pi^{-1}\Theta\Pi^{-1}, \quad \Theta'_- := \Theta'_{\ot_-} = \ov{\Theta^{\op}}.
\end{equation}
\Cref{lem:RThetaRTheta'} also holds in these cases.
\end{rem}

\begin{rem}\label{rem:bar-involution pullback}
    For a $\Uqg$-module $M$, consider a $\Uqg$-module $\ov{M}$ whose underlying set is $M$ (we denote its elements by $\ov{m}$ for $m\in M$), and the action of $\Uqg$ is given by $u \ov{m} = \ov{\ov{u}m}$ for $u\in \Uqg$ and $m\in M$.
    By \Cref{rem:$R$ matrix}, 
    \begin{equation*}
        \ov{M\ot_+ N} \cong M\ov{\ot}_+ N \cong N\ot_- M.
    \end{equation*}
    
\end{rem}

We recall the definition of a based module, introduced in \cite{Lu93}, where the condition (4) below is modified so that it is compatible with crystal lattices at $q=0$.
A lower (resp. upper) based module is a pair $(V,B)$, where $V$ is an integrable $\Uqg$-module with a $\Q(q)$-basis $B$ satisfying
\begin{enumerate}
    \item $B\cap V_\lambda$ is a basis of $V_\lambda$ for $\lambda\in P$,
    \item The $A$-submodule $V_A$ generated by $B$ is stable under the action of $U_q(\g)_A$,
    \item The $\Q$-linear involution $-$ on $V$ given by $\ov{c(q)b} = c(q^{-1})b$ for all $c(q)\in \Q(q)$ and $b\in B$ is compatible with $\Uqg$-action, that is, $\ov{u\,v} = \ov{u}\,\ov{v}$ for all $u\in \Uqg$ and $v\in V$,
    \item The $A_0$-submodule generated by $B$ is stable under lower (resp. upper) crystal operators.
\end{enumerate}

Suppose that there exist $\Q$-linear involutions $-$ on $\Uqg$-modules $V_1$ and $V_2$ compatible with the $\Uqg$-action.
For $v_1\ot v_2\in V_1\ot_+ V_2$, we define
\begin{equation}\label{eq:bar involution on tensor product}
\ov{v_1\ot v_2} = \Theta(\ov{v_1}\ot \ov{v_2}).
\end{equation}
In general, \eqref{eq:bar involution on tensor product} is well-defined only in a certain completion of $V_1\ot V_2$.
If \eqref{eq:bar involution on tensor product} gives a well-defined element in $V_1\ot V_2$ for all $v_1\ot v_2$ (for example, when $V_1$ is a lowest weight module or $V_2$ is a highest weight module), then it also gives a $\Q$-linear involution on $V_1\ot_+ V_2$,
such that $\ov{u \cdot (v_1\ot v_2)} = \ov{u} \cdot \ov{v_1\ot v_2}$ for $u\in \Uqg$.

One can define a bar-involution on a tensor product of more than two modules by applying \eqref{eq:bar involution on tensor product} inductively, which does not depend on the order of application of \eqref{eq:bar involution on tensor product} due to
$\left((1\ot \Delta_+)\Theta\right)\Theta^{23} = \left((\Delta_+\ot 1)\Theta\right)\Theta^{12}$ (see \cite[Proposition 4.2.4]{Lu90}):
Define $\Theta^{(n)}$ as a formal sum in $\Uqg^{\ot n}$ inductively by $\Theta^{(2)} = \Theta$ and $\Theta^{(n)} = (\Delta_+ \ot 1^{\ot (n-2)})\Theta^{(n-1)}(\Theta\ot 1^{\ot (n-2)})$ for $n\ge 3$.
Then the map
\begin{equation}\label{eq:bar involution on n-tensor product}
    \ov{v_1\ot\cdots\ot v_n} = \Theta^{(n)}(\ov{v_1}\ot\cdots\ot \ov{v_n}).
\end{equation}
for $v_1\ot\cdots\ot v_n \in V_1\ot_+\cdots\ot_+ V_n$ gives a $\Q$-linear involution compatible with the action of $\Uqg$.

The opposite quasi-$R$-matrix $\Theta'$ of \eqref{eq:opposite quasi-$R$-matrix} also satisfies the following so that it can be used to define an involution on a $\Uqg$-module $V_1\ot_+ V_2$ as well.
\begin{lem} For $u\in \Uqg$, we have the following identities:
    \[
        \Delta_+(\ov{u})\Theta' = \Theta' \ov{\Delta}_+(u), \quad
        \left((1\ot \Delta_+)\Theta' \right) \Theta'^{23} = \left((\Delta_+\ot 1)\Theta'\right)\Theta'^{12}.
    \]
\end{lem}
\pf By the natural isomorphism $\Pi$ in \Cref{rem:$R$ matrix}, one has $\Delta_+(u)\Pi = \Pi \ov{\Delta}_+^{\op}(u)$. Then,
\begin{equation}\label{eq:proof of first identity of opposite Theta}
    \Delta_+(\ov{u})\Pi \ov{\Theta}^{\op} \Pi = \Pi \ov{\Delta_+^{\op}(u)} \ov{\Theta}^{\op} \Pi =
    \Pi \ov{ \Theta^{\op} \ov{\Delta_+^{\op}(\ov{u})}} \Pi = \Pi \ov{\Theta}^{\op} \Delta_+^{\op}(\ov{u}) \Pi = \Pi \ov{\Theta}^{\op} \Pi \ov{\Delta_+(u)},
\end{equation}
which proves the first identity.
Next, using \eqref{eq:proof of first identity of opposite Theta} and the fact that $\Pi^{(3)} := (1\ot \Delta)(\Pi) (1\ot \Pi)$ is symmetric on each tensor component, we have
\begin{align*}
    \left((1\ot \Delta_+)\Theta' \right) \Theta'^{23} &= \left((1\ot \Delta_+)\Pi\right) \left( (1\ot \Delta_+) \ov{\Theta}^{\op} \right) \Pi^{(3)} \left( 1\ot \ov{\Theta}^{\op}\right) (1\ot \Pi) \\
    &= \Pi^{(3)} \left( (1\ot \ov{\Delta}_+^{\op}) \ov{\Theta}^{\op}\right) \left(  1\ot \ov{\Theta}^{\op}\right) \Pi^{(3)},
\end{align*}
and then $\left( (1\ot \ov{\Delta}_+^{\op}) \ov{\Theta}^{\op}\right) \left(  1\ot \ov{\Theta}^{\op}\right) = \ov{\left((\Delta_+\ot 1)\Theta\right) (\Theta\ot 1)}^{321}$ yields the second identity, where $(v_1\ot v_2\ot v_3)^{321} = v_3\ot v_2\ot v_1$.
\qed

\begin{rem}
By the natural isomorphisms in \Cref{rem:$R$ matrix}, we can also define a bar-involution on $V_1\ot_- V_2$ using $\Theta_-$ or $\Theta'_-$ in \eqref{eq:R-matrix for lower comultiplications}.
\end{rem}

The following is proved in \cite{Lu93} for $\lie{g}$ finite type, and for arbitrary $\lie{g}$ in \cite{BW}, in case of $\ot_+$.
\begin{thm}[{cf. \cite[Theorem 2.7]{BW}}]\label{thm:canonical basis of tensor product}
    Suppose that $(V_i,B_i)$ ($i=1,2$) are upper (resp. lower) based modules such that either $V_1$ is a lowest (resp. highest) weight module or $V_2$ is a highest (resp. lowest) weight module.
    Let $L_i$ be the $A_0$-span of $B_i$.
Then there exists a unique basis $B_1\bprod B_2 = \{\,b_1\bprod b_2 \,|\, b_i\in B_i\,\}$ of $V_1\ot_+V_2$ (resp. $V_1\ot_- V_2$), such that
\begin{enumerate}
\item $b_1\bprod b_2 \equiv b_1\ot b_2 \pmod{q(L_1\ot L_2)}$,
\item $\ov{b_1\bprod b_2} = b_1\bprod b_2$,
\end{enumerate}
where the bar-involution on $V_1\ot_+V_2$ (resp. $V_1\ot_- V_2$) is given by using $\Theta_+$ (resp. $\Theta'_-$).
Furthermore, $V_1\ot_+ V_2$ (resp. $V_1\ot_- V_2$) is an upper (resp. lower) based module with respect to $B_1\bprod B_2$.
\end{thm}
\pf 
In the case of $\ot_+$, the existence of $B_1\bprod B_2$ is proved in \cite{BW}, and the stability of a crystal lattice in this case is a result of a tensor product rule \eqref{eq:tensor product rule}.
The case of $\ot_-$ can be proved following arguments in \cite{BW} by changing the role of highest and lowest weight modules, since $\Theta'_-$ now lies in a completion of $\UqgPlus \ot \UqgMinus$, and for a weight vector $u\in \UqgMinus$,
\[
    \Delta_-(u) = u \ot 1 + \sum_{\wt{u_1} > \wt{u}} c_{u_1,u_1} u_1\ot u_2 \quad \left(u_1,u_2\in \UqgMinus\right).
\]
\qed

\section{Parabolic $q$-boson algebras}
\label{sec:parabolic q-boson algebra}
\subsection{Definition and basic properties}
\label{subsec:B_q}
We keep the notations in Section 2.
Let $J$ be a subset of $I$ such that $J^c:=I\setminus J$ is finite.
Let $\left(A_J = (a_{ij})_{i,j\in J}, P_J^\vee, P_J, \Pi_J^\vee, \Pi_J, (\cdot,\cdot) \right)$ be a Cartan datum of a submatrix $A_J$ of $A$
such that $\Pi_J^\vee = \{\, h_j \, | \, j\in J \,\} \subset \Pi^\vee$. Then there exists a canonical projection $P\to P_J$.
Let $\mf{l}=\g_J$ be the Kac-Moody algbra associated with the submatrix $A_J=(a_{ij})_{i,j\in J}$.
We may regard $\mf{l}$ as a subalgebra of $\g$, and let $\mf{p}=\mf{l}+\mf{b}$ denote the parabolic subalgebra, where $\mf{b}$ is the (positive) Borel subalgebra of $\g$.

Let $\Uqgp$ be an associative $\Q(q)$-algebra generated by $e'_i, e_j, f_l, q^{h}$ for $i\in J^c$, $j\in J$, $l\in I$, and $h\in P^\vee_J$ subject to the following relations:
{\allowdisplaybreaks
\begin{gather*}
q^0=1, \quad q^{h + h'}=q^{h}q^{h'} \quad  (h, h' \in P^\vee_J), \label{eq:Weyl-rel-1-Ya}\\
 q^h e'_i q^{-h}=q^{\langle h,\alpha_i\rangle}e'_i,\quad
 q^h e_j q^{-h}=q^{\langle h,\alpha_j\rangle}e_j,\quad
 q^h f_l q^{-h}=q^{-\langle h,\alpha_l\rangle}f_l,\\
 e_jf_l - f_le_j =
\delta_{jl}\frac{t_{j} - t^{-1}_{j}}{q_j-q_j^{-1}},\quad
e'_if_l = q_i^{-\langle h_i,\alpha_l \rangle} f_l e'_i + \delta_{il},\\
 \texttt{S}_{i_1,i_2}(e'_{i_1},e'_{i_2}) =
 \texttt{S}_{l_1,l_2}(f_{l_1},f_{l_2}) =
 \texttt{S}_{j_1,j_2}(e_{j_1}, e_{j_2}) = 0\quad (i_1,i_2\in J^c,\, j_1,j_2\in J,\, l_1,l_2\in I), \\
 \texttt{S}^-_{i, j}(e'_i, e_j) = \texttt{S}^+_{j,i}(e_j,e'_i) = 0,
\end{gather*}
\noindent where $t_j=q^{(\alpha_j,\alpha_j)h_j/2}$ and
\begin{align*}
    \texttt{S}_{s,t}(x,y)  & = \sum_{a+b = c_{st}} (-1)^a x^{(a)}yx^{(b)},               \\
    \texttt{S}_{s,t}^\pm(x,y) & = \sum_{a+b= c_{st}}(-1)^a q^{\pm a(\alpha_s,\alpha_t)} x^{(a)}yx^{(b)},\textsl{}
\end{align*}
with $x_l^{(a)}=x_l^a/[a]_l!$ for $x=e, e', f$.}
We call $\Uqgp$ the parabolic $q$-boson algebra associated to $(\lie{g}, \lie{p})$.
If $J=\emptyset$, then $\Uqgp$ is equal to the algebra of $q$-bosons denoted by $B_q(\g)$ in \cite{Kas91}.

Let $\UqgpMinus$ be the subalgebra of $\Uqgp$ generated by $f_l$ ($l\in I$), and let $\UqgpPlus$ be the subalgebra generated by $e'_i, e_j$ ($i\in J^c, j\in J$), where $\Uqgp^{\pm}$ is naturally graded by $Q_{\pm}$.
Let $\UqgpZero$ be the subalgebra generated by $q^h$ ($h\in P^\vee_J$).

\begin{lem}\label{lem:triangular decomposition of parabolic boson algebra}
There is an isomorphism of $\Q(q)$-spaces $\UqgpMinus\ot \UqgpZero\ot \UqgpPlus \to \Uqgp$, which is given by multiplication.
\end{lem}
\pf 
We may assume that $P^\vee_J = \Z \Pi_J^\vee$, since the general case follows from this case.
It is done by slightly modifying the arguments in \cite[Chapter 15]{Lu93}.
Let $U_q=U_q(\g)$. 
Let $\td{\tt B}_{q}^+$ be an associative $\Q(q)$-algebra generated by ${\tt e}'_i, {\tt e}_j$  ($i\in J^c, j\in J$) subject to the same relations for $e'_i$, $e_j$ in $\Uqgp$. Then it is straightforward to check that there exists an embedding of $\Q(q)$-algebras
\begin{equation*}
\begin{split}
 \xymatrixcolsep{2pc}\xymatrixrowsep{0pc}\xymatrix{
\td{\tt B}_{q}^+ \ \ar@{->}[r] &\ U_q^{\ge 0} := U_q^0U_q^+ \\
  {\tt e}'_i \ \ar@{|->}[r] &\ e_i' := -(q_i-q_i^{-1})t_ie_i\\
  {\tt e}_j \ \ar@{|->}[r] &\ e_j }
\end{split}
\end{equation*}
for $i\in J^c$ and $j\in J$. So we may identify $\td{\tt B}_q^+$ with its image in $U_q^{\ge 0}$, which is generated by $e'_i, e_j$ ($i\in J^c, j\in J$).
Let
\begin{equation}\label{eq:q-boson tilde}
 \td{\tt B}_q
 =U^-_q\,\td{\tt B}_q^0\,\td{\tt B}_q^+
 =\td{\tt B}_q^+\,\td{\tt B}_q^0\, U^-_q \subset U_q,
\end{equation}
where $\td{\tt B}_{q}^0$ denote the $\Q(q)$-subalgebra of $U_q^0$ generated by $t_i$, $q^h$ ($i\in J^c, h\in P^\vee_J$).
From the triangular decomposition of $U_q$, the multiplication in $U_q$ yields an isomorphism of $\Q(q)$-spaces $U^-_q\ot \td{\tt B}_q^0\ot \td{\tt B}_q^+ \longrightarrow U^-_q\,\td{\tt B}_q^0\,\td{\tt B}_q^+ $. We define $\tt{B}_q$ to be the quotient of $\td{\tt B}_q$ by the two-sided ideal $\tt{I}$ generated by $t_i$ ($i\in J^c$), and denote by $\tt{B}_q^+$, $\tt{B}_q^-$ and $\tt{B}_q^0$ the images of $\td{\tt B}_q^+$, $U_q^-$, and $\td{\tt B}_q^0$, respectively. Indeed, we have $\td{\tt B}_q^+\cong {\tt B}_q^+$ and $\td{\tt B}_q^-\cong U_q^-$.
Then there exists a well-defined homomorphism of $\Q(q)$-algebras
\begin{equation}\label{eq:B_q iso}
\begin{split}
 \xymatrixcolsep{2pc}\xymatrixrowsep{0pc}\xymatrix{
\Uqgp \ \ar@{->}[r] &\ \tt{B}_q
}
\end{split}
\end{equation}
from $\Uqgp$ to $\tt{B}_q$ sending $e'_i, e_j, t_j$ to $e'_i, e_j, k_j$ (in $U_q(\g)$ with the same notations) for $i\in J^c$ and $j\in J$. Since we have $\tt{B}_q^-\ot \tt{B}_q^0\ot \tt{B}_q^+\cong \tt{B}_q$ as a $\Q(q)$-space, the map in \eqref{eq:B_q iso} is an isomorphism. Hence we obtain the required isomorphism.
\qed

\begin{prop}\label{prop:comodule}
There is a homomorphism of $\Q(q)$-algebras $\Delta = \Delta_- : \Uqgp \to U_q(\g) \otimes \Uqgp$ such that
\begin{align*}
\Delta(q^{h}) & = q^h \ot q^h, \\
\Delta(e'_i)& = -(q_i - q_i^{-1})t_ie_i\ot 1 + t_i\ot e_i',\\
\Delta(e_j) & = e_j\ot t_j^{-1} + 1 \ot e_j, \\
\Delta(f_l) & = f_l\ot  1 + t_l \otimes f_l,
\end{align*}
for $i\in J^c$, $j\in J$, and $l\in I$, and $h\in P^\vee_J$.
Hence $\Uqgp$ becomes a left $U_q(\g)$-comodule algebra.
\end{prop}
\pf We may assume that $P^\vee_J = \Z \Pi_J^\vee$, since the general case follows from this case.
We see from the definition of $\td{\tt B}_q$ in the proof of the previous lemma that $\Delta(\td{\tt B}_{q})\subset U_q\ot \td{\tt B}_q$ and $\Delta({\tt I})\subset {\tt I}\ot {\tt I}$. This induces a well-defined map $\Delta : {\tt B}_q \longrightarrow U_q\ot {\tt B}_q$ given by the above formula.
\qed
\medskip

There exists an anti-involution $\tau = \tau_-$ of $\Uqgp$ defined by
\begin{equation}\label{tau for parabolic boson algebra}
        \tau(q^h) = q^h, \ 
        \tau(e'_i) = (1-q_i^2)f_i, \ 
        \tau(e_j) = q_jf_jt_j^{-1}, \ 
        \tau(f_i) = \frac{1}{1-q_i^2}e'_i, \ 
        \tau(f_j) = q_je_jt_j,
\end{equation}
for $i\in J^c$ and $j\in J$, and $h\in P^\vee_J$.
Then we can check the following.

\begin{prop}\label{prop:distributivity of tau on parabolic boson algebra}
    We have $(\tau \ot \tau) \circ \Delta_- = \Delta_- \circ \tau$ as homomorphisms of $\Q(q)$-algebras $\Uqgp \to U_q(\g)\ot \Uqgp$,
    where $\tau$ on $\Uqg$ on the left hand side is equal to $\tau_-$ in \eqref{eq:tau_pm for quantum groups}.
\end{prop}

\begin{rem}\label{rem:upper version of parabolic q-boson algebras}
    We may also define a parabolic $q$-boson algerbra associated with $\Delta_+$ of $\Uqg$.
    Let $\Uqgp^{\up}$ be a $\Q(q)$-algebra generated by $e_i'', e_j,f_l,q^h$ for $i\in J^c, j\in J, l\in I$, and $h\in P^\vee_J$,
    subject to the same relations for $e_j,f_l,q^h$ in $\Uqgp$, and
    \begin{gather*}
    e''_if_l = q^{\langle h_i,\alpha_l\rangle} f_le''_i + \delta_{il}, \\
    \texttt{S}_{i_1,i_2}(e''_{i_1},e''_{i_2}) =
    \texttt{S}^+_{i,j}(e''_i, e_j) = 
    \texttt{S}^-_{i,j}(e_i,e''_j) = 0 \quad (i_1,i_2, i\in J^c, j\in J).
    \end{gather*}
    There exists an isomorphism of $\Q$-algebras $- :\Uqgp\rightarrow \Uqgp^{\up}$ 
    given by $e'_i \mapsto e''_i, e_j \mapsto e_j, f_l\mapsto f_l, q^h\mapsto q^{-h}$, and $q \mapsto q^{-1}$.
    \Cref{lem:triangular decomposition of parabolic boson algebra} holds for $\Uqgp^{\up}$, and
    $\Delta_+$ induces a homomorphism of $\Q(q)$-algebras $\Delta_+ : \Uqgp^{\up} \to \Uqgp^{\up} \ot \Uqg$ given by
    \begin{align*}
    \Delta_+(q^{h}) & = q^h \ot q^h, \\
    \Delta_+(e''_i)& 
    = e''_i \ot t_i^{-1} + 1 \ot (q_i - q_i^{-1})t_i^{-1}e_i, \\
    \Delta_+(e_j) & = e_j\ot 1 + t_j \ot e_j, \\
    \Delta_+(f_l) & = f_l\ot t_l^{-1} + 1 \otimes f_l,
    \end{align*}
    for $i\in J^c, j\in J$, and $l\in I$, and $h\in P^\vee_J$.
    Similarly, there exists an anti-involution $\tau_+$ of $\Uqgp^{\up}$ defined by
    \begin{gather*}
        \tau_+(q^h) = q^h, \quad
        \tau_+(e''_i) = (1-q_i^{-2})f_i, \quad
        \tau_+(e_j) = q_jt_jf_j, \quad \\
        \tau_+(f_i) = \frac{1}{1-q_i^{-2}}e''_i, \quad
        \tau_+(f_j) = q_jt_j^{-1}e_j,
    \end{gather*}
    for $i\in J^c$ and $j\in J$, and $h\in P^\vee_J$.
    Then we have $(\tau_+\ot \tau_+) \circ \Delta_+ = \Delta_+ \circ \tau_+$ as homomorphisms of $\Q(q)$-algebras $\Uqgp^{\up} \to \Uqgp^{\up} \ot \Uqg$.
\end{rem}

\subsection{Integrable representations}
\label{subsec:integrable representations of parabolic q-boson algebras}
Let $U_q(\mf{l})$ be the quantized enveloping algebra associated to $(A_J, P^\vee_J, P_J, \Pi^\vee_J, \Pi_J,(\cdot,\cdot))$. We regard it as a subalgebra of $U_q(\g)$ generated by $e_i,f_i, q^h$ for $i\in J$ and $h\in P^\vee_J$.
By the construction of ${\tt B}_q$, it is not difficult to see that $U_q(\lie{l})$ is a subalgebra of ${\tt B}_q$. Hence, we may regard it as a subalgebra of $\Uqgp$ by the isomorphism \eqref{eq:B_q iso}.

Let $\mc{O}_{\Uqgp}$ be the category of $\Uqgp$-modules $V$ such that
\begin{itemize}
 \item[(1)] $V$ has a weight space decomposition $V=\bigoplus_{\la\in P_J}V_\la$ with respect to $\UqgpZero$,
 \item[(2)] given $v\in V$, $\UqgpPlus_\beta v = 0$ for all but finitely many $\beta\in Q_+$,
\end{itemize}
and let $\mc{O}^{\rm int}_{\Uqgp}$ be the subcategory of $\mc{O}_{\Uqgp}$ consisting of $V$ such that
\begin{itemize}
\item[(3)] $V$ is integrable as a $\QUE{\lie{l}}$-module, that is, $e_i,f_i$ ($i\in J$) act locally nilpotently.
\end{itemize}
We call $V \in \mc{O}_{\Uqgp}^{\rm int}$ an integrable $\Uqgp$-module.
We call a non-zero weight vector $v\in V$ singular if $\UqgpPlus v=0$.

For $\la\in P^+_J$, the set of dominant integral weights with respect to $U_q(\lie{l})$, let $V_{\lie{l}}(\la)$ be the irreducible highest weight $U_q(\mf{l})$-module with highest weight $\la$.
We regard $V_{\lie{l}}(\la)$ as a module over the subalgebra $U_q(\mf{l})\UqgpPlus$ of $\Uqgp$ by letting $e'_iu_\la=0$ for $i\in J^c$, where $u_\la$ is a highest weight vector of $V_{\lie{l}}(\la)$. Let
\begin{equation}\label{eq:parabolic Verma-1}
 V_J(\la)=\Uqgp\ot_{U_q(\mf{l})\UqgpPlusEqual}V_{\lie{l}}(\la),
\end{equation}
where $\UqgpPlusEqual=\UqgpZero \UqgpPlus$.
Since $U_q(\mf{l})\UqgpPlusEqual\cong U_q^-(\mf{l}) \ot \UqgpPlusEqual$ and $\Uqgp\cong U_q^-(\g) \ot  \UqgpPlusEqual$ as $\Q(q)$-spaces, we have as $\Q(q)$-spaces
\begin{equation}\label{eq:parabolic Verma-2}
  V_J(\la)
  \cong U_q^-(\g)\ot_{U^-_q(\mf{l})}V_{\lie{l}}(\la)
  \cong U^-_q(\g)\bigg/ \sum_{j\in J} U^-_q(\g) f_j^{\langle h_j,\la \rangle + 1}.
\end{equation}

\begin{lem}
For $\la\in P_J^+$, $V_J(\la)$ is a $\Uqgp$-module in $\mc{O}^{\rm int}_{\Uqgp}$.
\end{lem}
\pf It is clear from \eqref{eq:parabolic Verma-1} that $V_J(\la)$ belongs to $\mc{O}_{\Uqgp}$, so it suffices to show the condition (3). This follows from the fact that $V_J(\la)$ is generated by a singular vector $1\ot u_\la$ as a $\UqgMinus$-module, and the identity (3.4.1) of \cite{Kac}.
\qed
\medskip

The following lemma justifies our notation $e'_i$ for the generator of $\Uqgp$.

\begin{lem}\label{lem:identifying Kashiwara operator and an action of a parabolic boson algebra}
    Let $u\in \UqgMinus$ be given.
    Under the identification of \eqref{eq:parabolic Verma-2}, the action of $e'_i$ on the image of $u$ in $V_J(\la)$ equals the image of $e'_i(u)$ for $i\in J^c$, where $e'_i$ is understood as a linear endomorphism of $\UqgMinus$ defined in \eqref{eq:Kashiwara operators}.
\end{lem}
\pf
For $u\in \UqgMinus$, we can prove $e'_iu - q_i^{\langle h_i \wt u\rangle}ue'_i = e'_i(u)$ by induction on the height of $\wt(u)$, where both sides are understood as an element of $\Uqgp$.
Then $e'_i(u\ot u_\la) = e'_iu \ot u_\la = (e'_iu - q_i^{\langle h_i, \wt u\rangle}ue'_i) \ot u_\la = e'_i(u) \ot u_\la$, and the statement follows.
\qed

\begin{rem}\label{rem:analogue of identification of Kashiwara operators for upper parabolic boson algebras}
    We may define the category $\mc{O}^{\rm int}_{\Uqgp^{\up}}$ of integrable $\Uqgp^{\up}$-modules and $V_J(\la)$ for $\la\in P_J^+$ in the same way.
    We have an analogue of \Cref{lem:identifying Kashiwara operator and an action of a parabolic boson algebra} with $e'_i$ replaced by $e''_i$ in \eqref{eq:Kashiwara operators}.
\end{rem}

\subsection{Complete reducibility}\label{subsec:complete reducibility}
Let us introduce an analogue of the quantum Casimir operator on $\Uqgp$-modules, which is defined in a similar way as in \cite{Lu93}.

For $\alpha\in Q_+$, let $U_q^+(\g)_{\ge \alpha}=\bigoplus_{\beta\ge \alpha} {U^+_q(\g)}_{\beta}$ and
${\UqgpPlus}_{\ge \alpha}= \bigoplus_{\beta\ge \alpha} {\UqgpPlus}_{\beta}$.
Let $\widehat{U}_q(\lie{g})$ be the inverse limit of an inverse system $(\Uqg/U_q^{\le 0}(\g)U^+_q(\g)_{\ge \alpha})_{\alpha\in Q_+}$ with respect to the maps $f_{\alpha\beta}: \Uqg/U_q^{\le 0}(\g)U^+_q(\g)_{\ge \beta} \to \Uqg/U_q^{\le 0}(\g)U^+_q(\g)_{\ge \alpha}$ indexed by a directed set $Q_+$.
Let $\widehat{U}_q(\lie{g},\lie{p})$ be defined in a similar way with respect to $\{\,\Uqgp/\UqgpMinusEqual\UqgpPlus_{\ge \beta}\,|\,\beta\in Q_+\}$.
Let
\begin{equation}\label{eq:Casimir}
 \Omega=m((S_+\ot 1)\Theta_+),
\end{equation}
where $m$ denotes the multiplication in $U_q(\g)$, $\Theta_+$ is given in \eqref{eq:quasi R matrix}, and $S_+$ is the antipode for $\Delta_+$ given by $S_+ = (\ov{S})^{-1}$ with $\ov{S}(u)=\ov{S(\ov{u})}$ for $u\in U_q(\g)$. It is a well-defined element in $\widehat{U}_q(\g)$.

\begin{lem}\label{lem:Omega} 
$\Omega$ induces a well-defined element in $\widehat{U}_q(\lie{g},\lie{p})$
satisfying
\begin{enumerate}
 \item $e_j\Omega = t_j^{2}\Omega e_j$, $f_j\Omega = t_j^{-2}\Omega f_j$, $q^h \Omega =\Omega q^h$  for $j\in J$ and $h\in Q_J$,

 \item $e'_i\Omega = \Omega f_i=0$  for $i\in J^c$.
\end{enumerate}
\end{lem}
\pf Let $\td{\tt B}_q$ be the subalgebra of $U_q(\g)$ given in \eqref{eq:q-boson tilde}. We define the completion $\td{\tt B}_q^{{}^{\wedge}}$ in a similar way. We may regard $\td{\tt B}_q^{{}^{\wedge}}$ as a subalgebra $\widehat{U}_q(\g)$, and $\widehat{U}_q(\lie{g},\lie{p})$ is a quotient of $\td{\tt B}^{{}^{\wedge}}_q$ by the two-sided ideal $\tt{I}^{{}^{\wedge}}$ generated by $t_i$ ($i\in J^c$).

Recall from \cite[4.1]{Lu93} that we have
\begin{equation*}
(\Theta_+)_\beta = \sum_{b\in {\bf B}_\beta}c_b b^*\ot b \quad (\beta\in Q_+),
\end{equation*}
for some $c_b\in \Q(q)$, where ${\bf B}_\beta$ is a $\Q(q)$-basis of $U^+_q(\g)_\beta$ and ${\bf B}^*_\beta:=\{\,b^*\,|\,b\in {\bf B}_\beta\,\}$ is a dual basis of $U^-_q(\g)_{-\beta}$, and hence
\begin{equation*}
\begin{split}
 \Omega &= \sum_{b\in {\bf B}_\beta}c_b S_+(b^*) b.
\end{split}
\end{equation*}
Since $S_+(f_i)=-f_it_i$ for $i\in I$, we see that $S_+(b^*)b\in \td{\tt B}_q$, and $\Omega$ induces a well-defined element in $\td{\tt B}_q^{{}^{\wedge}}$, which we still denote by $\Omega$.

Now, the equations in (1) follow from \cite[6.1.2]{Lu93}. This also implies $e'_i\Omega = \Omega f_i=0$ in $\widehat{U}_q(\lie{g},\lie{p})$ since $e'_i\Omega = t_i^2 \Omega e'_i$ and $t_i^2  f_i\Omega = \Omega f_i$ in $\td{\tt B}_q$ ($i\in J^c$), which belong to $\tt{I}^{{}^{\wedge}}$. This proves (2).
\qed
\medskip

\begin{rem}
Although $\Uqgp$ is compatible with $\Delta_-$ (for example, \Cref{prop:comodule} and \Cref{prop:distributivity of tau on parabolic boson algebra}), we have to use $\Omega$ associated to $\Delta_+$ in order to have $\Omega \in \td{\tt B}_q^{{}^{\wedge}}$.
Indeed, the quantum Casimir element associated with $\Delta_-$ does not induce a well-defined element in $\td{\tt B}_q^{{}^{\wedge}}$.
\end{rem}

Let $V\in \mc{O}_{\Uqgp}$ be given.
Let $\Xi$ be a linear operator on $V$ given by $\Xi v= q^{(\lambda+2\rho, \lambda)}v$ for $v\in V_\la$, where $\rho\in P$ is given by $(\rho,\alpha_i)=(\alpha_i,\alpha_i)/2$ for $i\in I$.

\begin{prop}\label{prop:Casimir}
We have the following.
\begin{enumerate}
 \item If $v\in V_\la$ is singular, then $\Omega\Xi v = q^{(\la+2\rho,\la)} v$.

 \item $\Omega\Xi$ commutes with the $U_q(\mf{l})$-action on $V$.
\end{enumerate}
\end{prop}
\pf (1) is clear from the definition of $\Omega$. (2) follows from \Cref{lem:Omega} (1).
\qed
\medskip

\begin{prop} \label{prop:single_simple}
If $V\in\mc{O}^{\rm int}_{\Uqgp}$ is generated by a singular vector, then $V$ is irreducible. In particular, $V_J(\la)$ is irreducible for $\la\in P^+_J$.
\end{prop}
\pf Let $v$ be a singular vector of weight $\la$, which generates  $V$.
    Suppose that there is a proper submodule $W$. Then $W$ has a singular vector $w\in \UqgpMinus v$ of weight $\mu$.
    Since $W$ is integrable as a $U_q(\lie{l})$-module, we have $\mu\in P^+_J$ (cf. \cite[Proposition 3.5.8]{Lu93}).

    Let $u\in \UqgpMinus$ be such that $uv = w$. By \Cref{prop:Casimir}(1), we have $(\lambda + 2\rho, \lambda) = (\mu + 2\rho, \mu)$. Since $\la\ge \mu$ and $\mu\in P^+_J$, we have $\la=\mu$ by standard arguments (cf.~\cite{Kac}), and hence $u$ is a linear combination of monomials in $f_i$'s ($i\in J^c$). Since $W$ is a proper submodule, $u$ is not a constant and $\Omega u =0$ by \Cref{lem:Omega}(2). So $\Omega\Xi w   =  0$, which is again a contradiction.
\qed

\begin{lem}\label{lem:Casimir 2}
The following identities hold on $V$:
\[
\Omega \Xi \Omega = \Omega_{\lie{l}} \Xi \Omega = \Omega \Xi \Omega_{\lie{l}},
\]
where $\Omega_{\mf l}$ denotes the Casimir element \eqref{eq:Casimir} for $U_q(\mf{l})$.
\end{lem}
\pf It suffices to show that the identities hold for $v\in V_\la$. Since $\Xi$ is a scalar multiplication, the statement follows from \Cref{lem:Omega}(2).
\qed

\begin{thm}\label{thm:semisimple}
    The category $\mc{O}^{\rm int}_{\Uqgp}$ is semisimple with irreducibles $V_J(\la)$ for $\la\in P^+_J$.
\end{thm}
\pf Let $V\in \mc{O}^{\rm int}_{\Uqgp}$, and let $V'$ be the submodule generated by the singular vectors. Then $V'$ is semisimple by \Cref{prop:single_simple}. We claim that $V = V'$.

	Suppose that $V\ne V'$. There exists a nonzero $v\in V$ of weight $\mu$, which gives a singular vector of $V/V'$. Consider a $U_q(\mf{l})$-submodule $W$ generated by $\Omega v$. Since $W$ is a semisimple $U_q(\lie{l})$-module, we may write $\Omega v = v_1 + \cdots + v_n$, where $v_1,\dots, v_n$ belong to mutually non-isomorphic isotypic components.
	Hence $\Omega_{\lie{l}}\Xi v_i = q^{a_i} v_i$ for $a_i\in \Q$ ($i=1,\dots,n$), where $a_1,\dots,a_n$ are pairwise distinct.
    By \Cref{prop:Casimir} (1) and \Cref{lem:Casimir 2}, the following identities hold for all $N\ge 1$:
    \begin{gather*}
        \sum_{k=1}^n q^{Na_k}e_jv_k = e_j (\Omega_{\lie{l}} \Xi)^N\Omega v = (\Omega_{\lie{l}}\Xi)^N\Omega t_j^2e_jv \in V' \quad (j\in J), \\
        \sum_{k=1}^n q^{Na_k}e'_iv_k = e'_i (\Omega_{\lie{l}} \Xi)^N \Omega v = e'_i (\Omega\Xi)^N \Omega v = 0 \quad (i\in J^c).
    \end{gather*}
    Thus, $e_jv_k \in V'$ and $e'_iv_k = 0$ for $1\le k \le n$.
    Since $v\not\in V'$ and $\Omega\Xi v - q^{(\la+2\rho,\la)}v\in V'$, at least one among $v_k$ is not in $V'$. So by replacing $v$ with $v_k$, we may assume that $e'_iv=0$ for all $i\in J^c$, and that $\Omega_{\mf l}\Xi v = q^a v$ for some $a$.
    By applying $\Omega_{\mf l}\Xi$ on $V/V'$, we have $(\mu + 2\rho,\mu) = a$.

    By the condition (2) in the definition of $\mathcal{O}_{\Uqgp}$, there exists a $U_q(\lie{l})$-highest weight vector $w \in U_q(\mf{l}) v$ of weight $\lambda$.
    Then $\Omega_{\mf l}\Xi w = q^{(\lambda + 2\rho, \lambda)}w$, and  $(\mu + 2\rho, \mu) = (\lambda + 2\rho, \lambda)$,
    and $\lambda \ge \mu$, which implies $\lambda = \mu$. So $v$ is a singular vector with respect to $U_q(\mf{l})$-action, and hence a singular vector with respect to $\Uqgp$-action, but this contradicts the assumption that $v\not\in V'$.

    Finally, let $V$ be an irreducible $\Uqgp$-module in $\mc{O}^{\rm int}_{\Uqgp}$, which is generated by a singular vector $v$ of weight $\la$. Since $U_q(\mf{l})v$ is isomorphic to $V_{J}(\la)$ with $\la\in P_J^+$, there exists a surjective homomorphism $V_J(\la)$ onto $V$, and hence it is an isomorphism.
\qed

\begin{rem}
    The complete reducibility also holds for $\mc{O}^{\rm int}_{\Uqgp^{\up}}$,
    since twisting the action by $-:\Uqgp \to \Uqgp^{\up}$ in \Cref{rem:upper version of parabolic q-boson algebras} induces an equivalence from $\mc{O}^{\rm int}_{\Uqgp^{\up}}$ to $\mc{O}^{\rm int}_{\Uqgp}$.
\end{rem}

\subsection{A parabolic analogue of $q$-derivations}
\label{subsec:parabolic q-derivations}
We introduce a parabolic analogue of $q$-derivations  ${}_ir, r_i$ for $\UqgMinus$.

Let $\MM_J = V_J(0)$, and denote its singular vector by $1\in \MM_J$.
For $i\in J^c$, $f_i \in \MM_J$ is a highest weight vector of weight $-\alpha_i$ with respect to the $U_q(\lie{l})$-action, which is unique up to scalar multiplication.
Hence, there exists a $U_q(\lie{l})$-linear projection $\pi_i: \MM_J \longrightarrow V_{\lie{l}}(-\alpha_i)$. 
Regarding $\MM_J$ as a $\Uqg$-module, that is, as a parabolic Verma module induced from $V_{\lie{l}}(0)$, we have an injective $\Uqg$-linear map $\MM_J \longrightarrow \MM_J\ot_\pm \MM_J$, sending $1$ to $1\ot 1$.
Then we define $U_q(\lie{l})$-linear maps
\begin{equation}\label{eq:def of derivations}
    \begin{tikzcd}[row sep=tiny]
        {r_i^\pm: \MM_J} & {\MM_J\ot_{\pm} \MM_J} & {\MM_J\ot_\pm V_{\lie{l}}(-\alpha_i)}, \\
        {{}_ir^\pm : \MM_J} & {\MM_J \ot_\pm \MM_J} & {V_{\lie{l}}(-\alpha_i) \ot_\pm \MM_J}.
        \arrow[from=1-1, to=1-2]
        \arrow["{1\ot \pi_i}", from=1-2, to=1-3]
        \arrow[from=2-1, to=2-2]
        \arrow["{\pi_i \ot 1}", from=2-2, to=2-3]
    \end{tikzcd}
\end{equation}
If $J = \emptyset$, then $V(-\alpha_i)$ is one-dimensional $\Q(q)$-subspace of weight $-\alpha_i$, $r^+_i$ and ${}_ir^+$ coincides with the maps in \cite{Lu93}.
Note that
\begin{equation}\label{eq:intertwining upper and lower versions of derivations}
\sigma \circ (\-- \ot \--) \circ r_i^\pm = {}_i r^\mp,
\end{equation}
where $\sigma$ is the $\Q(q)$-linear map given by $\sigma(x\ot y) = y\ot x$.
From the definition, we clearly have
$({}_{i_1}r^\pm\ot 1) \circ r^\pm_{i_2} = (1\ot r^\pm_{i_2}) \circ {}_{i_1}r^\pm$ ($i_1,i_2\in J^c$),
which can be viewed as an analogue of $e_i'e_j'' = q_i^{\langle h_i,\alpha_j\rangle}e_j''e_i'$ in \cite[Proposition 3.4.5]{Kas91}.

The following is an analogue of \cite[Lemma 3.4.7]{Kas91}, which plays an important role in \Cref{sec:Crystal valuation on FFinftyM}.

\begin{lem}\label{lem:ri-annihilator is 1}
    If $u\in \MM_J$ satisfies $r^+_i(u) = 0$ for all $i\in J^c$, then $u$ is a scalar multiple of $1$. The same holds for $r^-_i$, ${}_ir^+$, and ${}_ir^-$.
\end{lem}
\pf
    Note that $\MM_J$ has a canonical $Q_-$ grading induced from that of $\UqgMinus$ by \eqref{eq:parabolic Verma-2}.
    Suppose that $u\in (\MM_J)_{-\xi}$ with respect to this grading.
    Suppose that $\xi\ne 0$. Since $r^+_i$'s are $\QUE{\lie{l}}$-linear, any element of $\QUE{\lie{l}}u$ satisfies the same
    condition as $u$.
    By the condition (2) in the definition of $\mathcal{O}_{\Uqgp}$, there exists $\alpha \in Q^+_J$ such that
    $(\QUE{\lie{l}}[+]u)_{-\alpha} \ne 0$ and $(\QUE{\lie{l}}[+]u)_{-\beta} = 0$ for $0\le \beta < \alpha$.
    Note that $\alpha\ne 0$ since $(\MM_J)_{-\alpha^\vee_j} = 0$ for $j\in J$. Therefore, $\QUE{\lie{l}}[+]u$ contains a $\QUE{\lie{l}}$-highest weight vector $u'$ with $u'\in (\MM_J)_{-\xi'}$ with $\xi'\ne 0$ such that $r_i(u') = 0$ for all $i\in J^c$.
    Hence, by replacing $u$ with $u'$, we may assume that $u$ satisfies $e_ju = 0$ for all $j\in J$.

    Recall that
    \[
        \Delta_+(u) \in u \ot 1 + c e'_i(u) \ot f_i + \sum_{\xi \in Q_+, \xi \ne 0, \alpha_i} \Uqg[\le 0] \ot \UqgMinus_{-\xi}
    \]
    for some $c\in \Q(q)$ (cf. \cite[1.2.13]{Lu93}). Therefore, $r_i^+(u) = 0$ implies $e'_i(u) = 0$.
    By \Cref{lem:identifying Kashiwara operator and an action of a parabolic boson algebra}, 
    $u$ is a singular vector of $\MM_J$ with respect to the $\Uqgp$-action.
    Then by \Cref{prop:single_simple}, $u$ is a scalar multiple of 1, which is a contradiction.

    A similar argument using $e''_i$ instead of $e'_i$ and an application of \Cref{rem:analogue of identification of Kashiwara operators for upper parabolic boson algebras} proves the case for ${}_ir^+$. The cases for $r_i^-$ and ${}_ir^-$ follows from \eqref{eq:intertwining upper and lower versions of derivations}.
\qed

\subsection{Crystal base of $V_J(\la)$ and tensor product rule}
\label{subsec:crystal base of parabolic Verma modules}
Let $\mc{O}^{\rm int}_{B_{q}(\g)} = \mc{O}^{\rm int}_{\Uqgp}$ when $J = \emptyset$.
Recall that
$U^-_q(\g)$ is a unique irreducible $B_q(\g)$-module in $\mc{O}^{\rm int}_{B_{q}(\g)}$ with a crystal base $(\ms{L}(\infty),\ms{B}(\infty))$ \cite{Kas91}, where the crystal operators $\te_i,\tf_i$ for $i\in I$ are given by
\begin{equation}\label{eq:Kashiwara operator for parabolic}
\te_i u=\sum_{k\geq1}f_i^{(k-1)}u_k,\quad
\tf_i u=\sum_{k\geq0}f_i^{(k+1)}u_k,
\end{equation}
for $u\in U_q^-(\g)_{\beta}$ ($\beta\in Q_-$) with $e'_iu_k=0$.
Let $G(\infty)=\{\,G(b)\,|\,b\in \ms{B}(\infty)\,\}$ be the global crystal basis or canonical basis of $U^-_q(\g)$.

Let $V\in \mc{O}^{\rm int}_{\Uqgp}$ be given.
We define a crystal base of $V$ in the same way as in the case of integrable $U_q(\g)$-modules in \Cref{subsec:crystal base} with respect to $\te_i, \tf_i$ in \eqref{eq:Kashiwara operator for parabolic} for $i\in J^c$ and \eqref{eq:Kashiwara operator lower} for $i\in J$, where $P$ is replaced by $P_J$.

Let $\la\in P_J^+$ be given. By \eqref{eq:parabolic Verma-2}, $V_J(\la)$ can be identified with
\begin{equation*} 
U^-_q(\g)\bigg/ \sum_{j\in J} U^-_q(\g) f_j^{\langle h_j,\la \rangle + 1}
\end{equation*}
as a $\Q(q)$-space, where the action of $e'_i, f_i$ for $i\in J^c$ on $V_J(\la)$ coincide with those induced from $U^-_q(\g)$ (cf. \Cref{lem:identifying Kashiwara operator and an action of a parabolic boson algebra}).
Let $u_\la$ denote the highest weight vector of $V_J(\la)$.

Let 
\[
    \pi^J_\la: U^-_q(\g) \longrightarrow V_J(\la)
\]
be the canonical projection, which is a homomorphism of $B_q(\lie{g}_{J^c})$-modules.
Let $\ms{B}_J(\la)=\{\,b\in \ms{B}(\infty)\,|\,\varepsilon_j^*(b)\le \langle h_j,\la \rangle\ (j\in J)\,\}$, where $\varepsilon_j^*(b)=\max\{\,k\,|\,\te_j^k b^* \neq 0 \,\}$ and $\ast$ denotes the involution on $\ms{B}(\infty)$ induced from the $\ast$-involution on $U_q^-(\g)$ \cite{Kas93}.
Then $\pi^J_\la(G(b))\neq 0$ if and only if  $b\in \ms{B}_J(\la)$, and $\{\,G_{J,\la}(b):=\pi^J_\la(G(b))\,|\, b\in \ms{B}_{J}(\la)\,\}$ forms a $\Q(q)$-basis of $V_J(\la)$. Let
\begin{equation*}
\begin{split}
    \LCrys_J(\la)&=\bigoplus_{b\in \ms{B}_J(\la)}A_0 G_{J,\la}(b), \\
    \BCrys_J(\la)&=\{\, G_{J,\la}(b) \!\!\!\! \pmod{q \LCrys_J(\la)}\,|\,b\in \ms{B}_{J}(\la)\,\}\setminus\{0\}.
\end{split}
\end{equation*}

\begin{thm}\label{thm:crystal base of V_J}
For $\la\in P_J^+$, $(\LCrys_J(\la),\BCrys_J(\la))$ is a crystal base of $V_J(\la)$.
\end{thm}
\pf It suffices to show that
\begin{equation}\label{eq:crystal lattice of V_J}
\begin{split}
& \td{x}_i\LCrys_J(\la) \subset \LCrys_J(\la),\quad
 \td{x}_i\BCrys_J(\la) \subset \BCrys_J(\la)\cup\{0\}\quad (i\in I, x=e,f)
\end{split}
\end{equation}
since the other conditions for crystal base follow immediately.

Since $\pi^J_{\la}$ commutes with $e'_i, f_i$ for $i\in J^c$, it also commutes with $\te_i, \tf_i$ for $i\in J^c$. This implies \eqref{eq:crystal lattice of V_J} for $i\in J^c$.

Let $\La\in P^+$ such that $\langle h_j,\La \rangle = \langle h_j,\la \rangle$ for $j\in J$. 
Let 
\begin{equation}\label{eq:pi lower J lambda}
    \pi_{J,\La}: V_J(\lambda) \longrightarrow V(\La)
\end{equation}
be the canonical projection of $U_q^-(\lie{g})$-modules sending $u_\lambda$ to $u_{\Lambda}$, so that $\pi_{J,\La} \circ \pi^J_\la$ is the canonical projection $\pi_\La: U^-_q(\g) \longrightarrow V(\La)$.
Recall that $\ms{B}(\La)=\{\,b\in \ms{B}(\infty)\,|\,\varepsilon_i^*(b)\le \langle h_i,\La \rangle\ (i\in I)\,\}$ and  $G(\La)=\{\,\pi_\La(G(b))\,|\,b\in \ms{B}(\La)\,\}$.

Let $\beta\in Q_+$ be given.
Choose $\La$ such that $\langle h_j,\La \rangle\gg 0$ for $j\in J$ so that $\LCrys_J(\la)_{\la-\beta+k\alpha_j}$ ($k=0,\pm1$) is isomorphic to $\ms{L}(\La)_{\La-\beta+k\alpha_j}$ under $\pi_{J,\La}$. Since $\pi_{J,\La}$ is $U_q(\mf{l})$-linear by \eqref{eq:parabolic Verma-2}, $\pi_{J,\La}$ commutes with $\te_j, \tf_j$ for $j\in J$. This implies that $\td{x}_j\LCrys_J(\la)_{\la-\beta} \subset \LCrys_J(\la)_{\la-\beta\pm\alpha_j}$ and $\td{x}_jG_{J,\la}(b)\in \ms{B}_{J}(\la)\cup\{0\} \!\!\pmod{q\LCrys_J(\la)}$ for $x=e,f$ and $j\in J$. Hence \eqref{eq:crystal lattice of V_J} holds for $j\in J$.
\qed

\begin{cor}
We have
\begin{equation*}
\begin{split}
\LCrys_J(\la)&=\sum_{r\geq 0,\, i_1,\ldots,i_r\in I}A_0 \tf_{i_1}\cdots\tf_{i_r}u_\la, \\
\BCrys_J(\la)&=\{\, \,\tf_{i_1}\cdots\tf_{i_r}u_\la\!\!\! \pmod{q \LCrys_J(\la)}\,|\,r\geq 0, i_1,\ldots,i_r\in I\,\}\setminus\{0\}.
\end{split}
\end{equation*}
In particular, $V_J(\la)$ has a unique crystal base up to scalar multiplication.
\end{cor}

We also have an analogue of $q$-Shapovalov form on $V_J(\la)$.

\begin{prop}\label{prop:q-Shapovalov form on parabolic Verma modules}
There exists a unique nondegenerate symmetric bilinear form $\langle \cdot, \cdot \rangle$ on $ V_J(\la)$ such that
\[
    \langle u_\la ,u_\la \rangle = 1,\quad \langle uv, w \rangle = \langle v, \tau(u) w\rangle \quad (u\in \Uqgp, \, v,w\in V_J(\la)),
\]
where $\tau$ is given in \eqref{tau for parabolic boson algebra}
\end{prop}
\pf
It can be proved by similar arguments as in the case of $J=I$ using \Cref{prop:single_simple}.
\qed

\medskip

\begin{prop}\label{thm:characterization of crystal lattice of parabolic Verma modules in terms of q-Shapovalov form}
    We have the following.
    \begin{enumerate}
    \item $\langle \tf_i u, v\rangle \equiv \langle u, \te_i v\rangle \pmod{qA_0}$ for $u,v\in \LCrys_J(\la)$ and $i\in I$.
    \item $\LCrys_J(\la) = \{\, x \in V_J(\la) \,|\, \langle x, \LCrys_J(\la) \rangle \subset A_0 \,\}$.
    \end{enumerate}
\end{prop}
\pf
(1) As in the proof of \cite[Proposition 5.1.1, 5.1.2]{Kas91}, we may show by induction on the height of $\xi \in -Q_+$ that
\begin{equation}\label{crystal operator adjoint}
\langle \tf_i u, v\rangle \equiv \langle u, \te_i v\rangle \pmod{qA_0}
\end{equation}
for $u\in \LCrys_J(\la)_{\la + \xi +\alpha_i}, v\in \LCrys_J(\la)_{\la + \xi}$.
For $i\in J$, the proof is identical to the one in \cite{Kas91}.
For $i\in J^c$, we have
\[
    \left\langle f_i^{(n+1)}u_0, f_i^{(m)} v_0 \right\rangle = \delta_{n+1,m}(1-q_i^2)^{-m}\left(q_i^{\frac{m(m-1)}{2}}[m]_i!\right)^{-1} \left\langle u_0, v_0\right\rangle,
\]
where $u = f_i^{(n)}u_0, v = f_i^{(m+1)}v_0$ with $e'_iu_0 = e'_iv_0 = 0$. This provides an analogue of \cite[(5.1.2)]{Kas91}, which is needed for the proof of \eqref{crystal operator adjoint}.
(2) can be proved by the same arguments as in \cite[Proposition 5.1.1]{Kas91}.
\qed
\medskip

Let $V_1$ be an integrable $U_q(\g)$-module, and let $V_2$ be an integrable $\Uqgp$-module. By \Cref{prop:comodule}, $V_1\ot V_2$ is a $\Uqgp$-module. Then we have the following.

\begin{thm}\label{thm:tensor product rule}
Let $(L_i,B_i)$ be a crystal base of $V_i$ $(i=1,2)$. Then $(L_1\ot L_2, B_1\ot B_2)$ is a crystal base of $V_1\ot V_2$ such that $\te_i$ and $\tf_i$ $(i\in I)$ act on $B_1\ot B_2$ by the same formula as in \eqref{eq:tensor product rule}.
\end{thm}
\pf Let $(L,B)=(L_1\ot L_2, B_1\ot B_2)$.
If $j\in J$, then $(L,B)$ is a crystal base of $V_1\ot V_2$ as a module over the subalgebra $\langle e_j,f_j,t_j^{\pm 1} \rangle$.
If $i\in J^c$, then $(L,B)$ is also a crystal base of $V_1\ot V_2$ as a module over the subalgebra $\langle e'_i,f_i \rangle$ satisfying \eqref{eq:tensor product rule} \cite[Section 3.5]{Kas91}.
\qed
\begin{rem}    \label{rem:characterizing highest weight elements}
In particular, $b_1\ot b_2\in B_1\ot B_2$ satisfies $\td{e}_i(b_1\ot b_2) = 0$ for all $i\in I$ if and only if 
$\td{e}_i b_1 = 0$ and $\varepsilon_i(b_2) \le \langle h_i, \mathrm{wt}(b_1) \rangle$ for all $i\in I$.
\end{rem}

\section{A category generated by extremal weight modules}
\label{sec:extremal weight modules of type A}

\newcommand{\iset}[1]{\mathsf{#1}} 
\newcommand{\dset}[1]{\mathsf{#1} \cap (\mathsf{#1} - 1)} 

\subsection{Notations}\label{sec:gln}
Let $\gl_\infty$ be a Lie algebra of $(\Z \times \Z)$-matrices spanned by the elementary matrices $E_{ij}$ ($i,j\in \Z$).
Let $U_q(\gl_{\infty})$ be the quantized enveloping algebra associated with
\begin{enumerate}
 \item $P=\bigoplus_{i\in\Z}\Z\epsilon_i\oplus\Z\Lambda_0 \oplus \Z\delta$, $P^\vee=\bigoplus_{i\in\Z}\Z E_{ii}$ with
\begin{gather*}
\langle E_{ii},\e_j \rangle=\de_{ij}\quad (i,j\in \Z),\quad
\langle E_{kk},\La_0\rangle = \begin{cases} 1 & k \le 0 \\ 0 & k > 0 \end{cases},\quad
\langle E_{kk}, \delta\rangle = 1 \quad (k\in \Z), \\
(\epsilon_i,\epsilon_j) = \delta_{ij}\quad (i,j\in \Z),\quad
(\epsilon_k, \La_0) =  \begin{cases} 1 & k \le 0 \\ 0 & k > 0 \end{cases},\quad
(\La_0,\La_0) = (\La_0,\delta) = (\delta,\delta) = 0,
\end{gather*}
 \item $\Pi=\{\,\alpha_i=\epsilon_i-\epsilon_{i+1}\,|\,i\in \Z\,\}$, $\Pi^\vee=\{\,h_i=E_{ii}-E_{i+1i+1}\,|\,i\in \Z\,\}$.
\end{enumerate}
Then $(A,P^\vee,P,\Pi^\vee,\Pi,(,))$ satisfies conditions in \Cref{subsec:quantum group}.
For $i\in \Z\setminus\{0\}$, let $\Lambda_i\in
P^+$ be the $i$-th fundamental weight given by
\begin{equation}\label{eq:def of fundamental weight}
\Lambda_i=
\begin{cases}
\Lambda_0+\sum_{k=1}^{i}\epsilon_k & \text{if $i>0$}, \\
\Lambda_0-\sum_{k=i+1}^{0}\epsilon_k & \text{if $i<0$}. \\
\end{cases}
\end{equation}
For $n\ge 1$, let
$\Z_+^n=\{\,(\lambda_1,\cdots,\lambda_n)\,|\,\lambda_i\in\Z,\
\lambda_1\geq \cdots\geq \lambda_n\,\}$. For $\lambda\in\Z_+^n$, we put
\begin{equation*}
\begin{split}
\Lambda_{\lambda}& =\Lambda_{\lambda_1}+\cdots+\Lambda_{\lambda_n}\in P^+.
\end{split}
\end{equation*}

For an interval $S\subset \mathbb{R}$, denote $\iset{S} = S\cap \Z$.
Let $A_S = (a_{ij})_{i,j\in \dset{S}}$, where $\iset{S}-1$ denotes the translate of $\iset{S}$ by $-1$.
Let $P^\vee_S = \bigoplus_{i\in \iset{S}} \Z E_{ii} \subset P^\vee$, and $P_S$ the dual of $P^\vee_S$ inside $P$.
Note that $P_S$ is a quotient of $P$, where $\epsilon_i=0$ for $i\in \Z\setminus \iset{S}$ together with some relations expressing $\La_0$, $\delta - \La_0$, or $\delta$ as linear combination of $\epsilon_i$'s (depending on $S$). Therefore, we may identify $P_S$ with a subgroup of $P$ containing $\{\epsilon_i\}_{i\in \iset{S}}$, and by restricting $(,)$ to $P_S$, we obtain a bilinear form $(,)_S$ on $P_S$.
We define $U_q(\gl_S)$ to be the quantized enveloping algebra associated with this realization, which is canonically identified as a subalgebra of $U_q(\gl_\infty)$.
We denote $U_q(\gl_S)$ by $U_q(\gl_{>0}), U_q(\gl_{\le 0})$, and $U_q(\gl_n)$ when $S = \R_{>0}, \R_{\le 0}$, and $[1,n]$, respectively.

Let $\cP$ denote the set of partitions $\la=(\la_i)_{i\ge 1}$.
For $\la\in \cP$, let $\ell(\la)$ be the length of $\la$ and $|\la|=\sum_{i\ge 1}\la_i$.

For $\la\in \cP$, let $\eps_\la = \sum_{i\ge 1}\la_i\epsilon_i$.
For $\lambda\in \Z^\ell_+$, we also write $\eps_{\lambda}=\lambda_1\epsilon_{1}+\lambda_2\epsilon_{2}+\cdots+\lambda_{\ell}\epsilon_{\ell}$.
For $\mu,\nu\in \cP^2$ with $n\ge \ell(\mu) + \ell(\nu)$, let $\eps_{\mu,\nu}^n = \mu_1\e_1 +\dots + \mu_{\ell(\mu)}\e_{\ell(\mu)}  - \nu_{\ell(\nu)} \e_{n-\ell(\nu)+1} -\dots -\nu_1\e_n$.

\subsection{Extremal weight modules $V_{\mu,\nu}$ over $\UqglPlus$}\label{subsec:ext wt module}
Let $P_{>0} = P_{\R_{ >0}}$ be the weight lattice of $\gl_{>0}$.
Let $P_{>0}/W_{>0}$ be the set of $W_{>0}$-orbits in $P_{>0}$, where $W_{>0}$ is the Weyl group of $\gl_{> 0}$.  Given
$\lambda=\sum_{i\geq 1}\lambda_i\epsilon_i\in P_{>0}$, we have
partitions $\mu$ and $\nu$, which are uniquely determined by sorting the coefficients $\la_i$ with $\la_i > 0$ and $\la_i <0$, respectively.
Then we have a well-defined bijection
\begin{equation}\label{eq:extremal wt bijection}
\xymatrixcolsep{2pc}\xymatrixrowsep{0pc}\xymatrix{
P_{>0}/W_{>0} \ \ar@{->}[r] &\ \cP^2=\cP\times\cP \\
W_{>0}\la \ \ar@{->}[r] &\ (\mu,\nu)
}.
\end{equation}
For $(\mu,\nu)\in \cP^2$, let $V_{\mu,\nu}$ denote the extremal weight module $V(\la)$, where $W_{>0}\la$ corresponds to $(\mu,\nu)$.
It has a crystal base $(\LCrys(V_{\mu,\nu}),\BCrys(V_{\mu,\nu}))$.

\begin{prop}\label{prop:irreducibility of extremal weight modules}
For $(\mu,\nu)\in \cP^2$, $V_{\mu,\nu}$ is irreducible.
\end{prop}
\pf Let $V$ be a non-zero submodule of $V_{\mu,\nu}$ and let $v\in V$ be a non-zero weight vector.
Let $s=\ell(\mu)$ and $t=\ell(\nu)$.
For $n\geq s + t$, let $u_{\mu,\nu}^n$ be an extremal weight vector of weight $\eps^n_{\mu,\nu}$,
which is a highest weight vector with respect to $U_q(\gl_n)$.
Since $V_{\mu,\nu}$ is generated by $u_{\mu,\nu}^n$, we may assume that $v$ belongs to 
\begin{equation}\label{eq:Vnmn}V^n_{\mu,\nu} := U_q(\gl_n)u_{\mu,\nu}^n\end{equation}
for a sufficiently large $n$. Since $V^n_{\mu,\nu}$ is an irreducible $U_q(\gl_n)$-module, this implies that $V$ contains $u_{\mu,\nu}^n$ and hence $V=V_{\mu,\nu}$.
\qed

\smallskip

Let us briefly recall a combinatorial realization of $\BCrys(V_{\mu,\nu}) = \ms{B}(\la)$ for $\la\in P_{>0}$ in \cite{K09}.
Let us regard $\Z_{>0}$ as the crystal of the natural representation of $U_q(\gl_{>0})$: $1 \stackrel{1}{\longrightarrow} 2
\stackrel{2}{\longrightarrow}3\stackrel{3}{\longrightarrow}
\cdots$ where ${\rm wt}(k)=\e_k$ for $k\in \Z_{>0}$. Let $\Z_{>0}^\vee=\{\,k^\vee\,|\,k\in \Z_{>0}\,\}$, and regard it as the dual crystal of $\Z_{>0}$, that is, $\cdots \stackrel{3}{\longrightarrow}
3^\vee\stackrel{2}{\longrightarrow}2^\vee\stackrel{1}{\longrightarrow}
1^\vee$ where ${\rm wt}(k^\vee)=-\e_k$ for $k\in \Z_{>0}$.

Let $(\mu,\nu)\in \cP^2$ correspond to $W_{>0}\la$ under \eqref{eq:extremal wt bijection}. Note that $V_{\mu,\emptyset}$ is a highest weight $\UqglPlus$-module with highest weight $\eps_\mu$ and $V_{\emptyset,\nu}$ is a lowest weight $\UqglPlus$-module with lowest weight $-\eps_\nu$.

For $\mu\in \cP$, we identify $\mu$ with its Young diagram and let $SST_{\Z_{>0}}(\mu)$ be the set of all
semistandard tableaux of shape $\mu$ with entries in $\Z_{>0}$ \cite{Fu}. Then the crystal $\ms{B}(V_{\mu,\emptyset})$ can be identified with $SST_{\Z_{>0}}(\mu)$. Indeed, we identify $S\in SST_{\Z_{>0}}(\mu)$ with $w_1\ot w_2\ot \dots \ot w_{|\mu|}$, where $w_1w_2\dots w_{|\mu|}$ is the column word of $S$, that is, a word given by reading the entries column by column
from right to left and from top to bottom in each column.
Then $SST_{\Z_{>0}}(\mu)$ is the connected component of $H_\mu$ in $SST_{\Z_{>0}}((1))^{\ot |\mu|}$, where $H_\mu$ is the unique element in $SST_{\Z_{>0}}(\mu)$ such that each $i$-th row (from the top) is filled with $i$ for $1\leq i \leq \ell(\mu)$, that is, ${\rm wt}(H_\mu)=\eps_\mu$.

Similarly, for $\nu\in \cP$, let $SST_{\Z_{>0}^\vee}(\nu^\pi)$ be the set of semistandard tableaux of shape $\nu$ with entries in $\Z_{>0}$, where $\nu^\pi$ is the skew diagram obtained by $180^\circ$-rotation of $\nu$ and $\Z_{>0}$ has a linear order $a^\vee<b^\vee$ for $a>b$. Then the crystal $\ms{B}(V_{\emptyset,\nu})$ can be identified with $SST_{\Z_{>0}^\vee}(\nu^\pi)$, which is the connected component of $H_\nu^\vee$ in $SST_{\Z_{>0}^\vee}((1))^{\ot |\nu|}$ (as its column word), where $H_\nu^\vee$ is the unique element in $SST_{\Z_{>0}^\vee}(\nu^\pi)$ such that each $i$-th row (from the bottom) is filled with $i^\vee$ for $1\leq i \leq \ell(\nu)$,
that is, ${\rm wt}(H_\nu^\vee)=-\eps_\nu$.
\begin{figure}[htbp]
    \centering
    \ytableausetup{boxsize=1.3em, centertableaux}
    \begin{subfigure}[b]{0.48\textwidth}
        \centering
        $\scalebox{0.8}{
        \begin{ytableau}
            1 & 1 & 1 & 1 \\
            2 & 2
        \end{ytableau}}
        \ \;\in \ SST_{\mathbb{Z}_{>0}}(\mu)$
        \caption*{$H_\mu$ for $\mu=(4,2)$}
    \end{subfigure}
    \hfill 
    \begin{subfigure}[b]{0.48\textwidth}
        \centering
        $\scalebox{0.8}{
        \begin{ytableau}
            \none & \raisebox{-0.2ex}{$2^\vee$} & \raisebox{-0.2ex}{$2^\vee$} \\
            \raisebox{-0.2ex}{$1^\vee$} & \raisebox{-0.2ex}{$1^\vee$} & \raisebox{-0.2ex}{$1^\vee$}
        \end{ytableau}}
        \  \;\in \ SST_{\mathbb{Z}_{>0}^\vee}(\nu^\pi)
        $
        \caption*{$H_\nu^\vee$ for $\nu=(3,2)$}
    \end{subfigure}
    \label{fig:semistandard_examples}
\end{figure}

Let $\ms{B}_{\mu,\nu}$ to be the set of bitableaux $(S,T)$ such that
\begin{itemize}
\item[(1)] $S\in SST_{\Z_{>0}}(\mu)$ and $T\in SST_{\Z_{>0}^\vee}(\nu^\pi)$,

\item[(2)] $\big|\,\{\,i\,|\,S^{(i,1)}\leq k\,\}\,\big|+\big|\,\{\,i\,|\,T_{(i,1)}\geq k^\vee\,\}\,\big|\leq k$ for all $k\geq 1$.
\end{itemize}
Here $S^{(i,j)}$ is the entry of $S$ in the $i$-th row from
the top and the $j$-th column from the left, and $T_{(i,j)}$ is the
entry of $T$ in the $i$-th row from the bottom and the $j$-th
column from the right.
We regard $\ms{B}_{\mu,\nu}\subset SST_{\Z_{>0}}(\mu)\ot SST_{\Z_{>0}^\vee}(\nu^\pi)$ and apply $\te_i$ and $\tf_i$ for $i\in\Z_{>0}$.

\begin{thm}[\cite{K09}]\label{thm: ext wt crystal of >0}
For $(\mu,\nu)\in \cP^2$, we have
\begin{enumerate}
 \item $\ms{B}_{\mu,\nu}\cup\{0\}$ is stable under $\te_i$ and $\tf_i$ for $i\in\Z_{>0}$ and it is connected,

 \item $\ms{B}_{\mu,\nu}$ is isomorphic to $\ms{B}(V_{\mu,\nu})$,

 \item $\ms{B}_{\mu,\nu}$ is isomorphic to $\ms{B}_{\emptyset,\nu}\ot \ms{B}_{\mu,\emptyset}$.

\end{enumerate}
\end{thm}

\begin{rem}\label{rem:tensor product}
It is shown in \cite{K09} that $\ms{B}_{\mu,\nu}\ot \ms{B}_{\sigma,\tau}$ is isomorphic to a disjoint union of $\ms{B}_{\zeta,\eta}$, where the multiplicity for each $\ms{B}_{\zeta,\eta}$ is given in terms of Littlewood-Richardson coefficients $c^{\gamma}_{\alpha\beta}$ for $\alpha,\beta,\gamma\in \cP$ with $|\alpha|+|\beta|=|\gamma|$ (cf. \cite{Fu}). For example, multiplicity for $\ms{B}_{\zeta,\eta}$ in $\ms{B}_{\mu,\emptyset}\ot \ms{B}_{\emptyset,\nu}$ is given by
\begin{equation}\label{eq:tensor multiplicity}
m^{\mu,\nu}_{\zeta,\eta}:=\sum_{\sigma}c^{\mu}_{\sigma\zeta}c^{\nu}_{\sigma\eta}.
\end{equation}
But we should remark  that $V_{\mu,\nu}\ot V_{\sigma,\tau}$ is not semisimple in general.
\end{rem}

\subsection{A monoidal category generated by $V_{\mu,\nu}$}\label{subsec:monoidal cat C}
\label{subsec:monoidal category generated by Vmunu}
In this subsection, we show that there exists a filtration on $V_{\mu,\emptyset}\ot V_{\emptyset,\nu}$, which is compatible with the decomposition of $\ms{B}_{\mu,\emptyset}\ot \ms{B}_{\emptyset,\nu}$ (\Cref{rem:tensor product}).

\begin{lem}\label{lem:embedding of ext wt module}
Let $V$ be an integrable $U_q(\gl_{>0})$-module with a crystal base $(L,B)$.
Let $v$ be a weight vector such that
\begin{enumerate}
 \item $v$ is an extremal weight vector of weight $\la$,
 \item $v\in B \pmod{qL}$ and its connected component in $B$ is isomorphic to $\ms{B}(V(\la))$.
\end{enumerate}
Then there exists an injective $U_q(\gl_{>0})$-linear map $\phi: V(\la)\longrightarrow V$ with $\phi(u_\la)=v$.
\end{lem}
\pf Let $V'$ be the $U_q(\gl_{>0})$-submodule of $V$ generated by $v$.
Then there exists a surjective $U_q(\gl_{>0})$-linear map $\phi: V(\la)\longrightarrow V'$ such that $\phi(u_\la)=v$.
Let
\begin{equation*}
\begin{split}
L'&=\sum_{r\geq 0,\, i_1,\ldots,i_r\in \Z_{>0}}A_0 \td{x}_{i_1}\cdots\td{x}_{i_r}v, \\
B'&=\{\, \,\td{x}_{i_1}\cdots\td{x}_{i_r} v\!\!\! \pmod{qL}\,|\,r\geq 0, i_1,\ldots,i_r\in\Z_{>0}\,\}\setminus\{0\}.
\end{split}
\end{equation*}
We have $\phi(\ms{L}(V(\la)))=L'$ and $\ov{\phi}(\ms{B}(V(\la)))=B'$ since $\phi$ is $U_q(\gl_{>0})$-linear, where $\ov{\phi}:\LCrys(\la)/q\LCrys(\la) \to L/qL$ is the map induced from $\phi$. Also $\ov{\phi}$ is a bijection since $B'$ is isomorphic to $\ms{B}(V(\la))$ and $\ov{\phi}$ commutes with $\te_i$ and $\tf_i$ for $i\in \Z_{>0}$.

Let $u\in V(\la)$ be given such that $\phi(u)=0$. Let $c\in \Q(q)^\times$ such that $cu \in \ms{L}(V(\la))$ and $cu =\sum_{b\in \ms{B}(V(\la))}c_b b \pmod{q\ms{L}(V(\la))}$ for some $c_b\in \mathbb{Q}$, which are not all $0$. Then $\ov{\phi}(cu) = \sum_{b\in \ms{B}(V(\la))}c_b \ov{\phi}(b)=0$, but this contradicts the fact that $\ov{\phi}$ is a bijection and $\ov{\phi}(\ms{B}(V(\la)))$ is linearly independent. Hence, $\phi$ is injective.
\qed
\medskip

Let $(\mu,\nu)\in \cP^2$ be given.
Let $(\ms{L}_{\mu,\nu},\ms{B}_{\mu,\nu})$ and $G_{\mu,\nu}$ denote the crystal base and global crystal basis of $V_{\mu,\nu}$, respectively.
Let
\begin{equation*}
\begin{split}
 V&=V_{\mu,\emptyset}\ot V_{\emptyset,\nu},\ \
 (\ms{L},\ms{B})=(\ms{L}_{\mu,{\emptyset}}\ot\ms{L}_{{\emptyset},\nu},\ms{B}_{\mu,{\emptyset}}\ot\ms{B}_{{\emptyset},\nu}).
\end{split}
\end{equation*}
Then $(\ms{L},\ms{B})$ is a crystal base of $V$. Let ${G}$ be the canonical basis of $V$ so that $(V,G)$ is a based module with respect to $(\ms{L},\ms{B})$ (cf. \Cref{thm:canonical basis of tensor product}).

\begin{lem}\label{lem:embedding of ext wt module-2}
There exists a unique injective $U_q(\gl_{>0})$-linear map $\phi: V_{\mu,\nu}\longrightarrow V_{\mu,\emptyset}\ot V_{\emptyset,\nu}$ such that $\phi(G_{\mu,\nu})\subset G$.
\end{lem}
\pf Let $u_{\mu,\nu}$ be an extremal weight vector of $V_{\mu,\nu}$ with weight $\eps^n_{\mu,\nu}$, say $\la$, where $n=\ell(\mu)+\ell(\nu)$.
Let $\la^+=\mu_1\e_1 +\dots + \mu_s\e_s$ and $\la^-= - \nu_t \e_{n-t+1} -\dots -\nu_1\e_n$.
Let $u_1\in V_{\mu,\emptyset}$ and $u_2\in V_{\emptyset,\nu}$ be unique weight vectors (up to scalar multiplication) such that ${\rm wt}(u_1)=\la^+ $ and ${\rm wt}(u_2)= \la^-$, respectively. We may assume that $u_1\ot u_2\in \ms{B} \pmod{q\ms{L}}$.

Then $u_1\ot u_2$ is an extremal weight vector of weight $\la$ in $V$, and the connected component of $u_1\ot u_2$ in $\ms{B}$ is isomorphic to $\ms{B}_{\mu,\nu}$ by \Cref{thm: ext wt crystal of >0}.
Hence there exists a unique injective $U_q(\gl_{>0})$-linear map $\phi: V_{\mu,\nu}\longrightarrow V$ such that $\phi(u_{\mu,\nu})=u_1\ot u_2$ by \Cref{lem:embedding of ext wt module}.

Since $u_{\mu,\nu}\in G_{\mu,\nu}$ and $u_1\ot u_2\in G$, we have
\begin{itemize}
\item[(1)] $\phi(G_{\mu,\nu}(b))\in \ms{L}\cap  U_q(\gl_{>0})_A(u_1\ot u_2)$,

\item[(2)] $\ov{\phi(G_{\mu,\nu}(b))}=\phi(G_{\mu,\nu}(b))$,

\item[(3)] $\phi(G_{\mu,\nu}(b))\equiv \ov{\phi}(b)\in \ms{B} \pmod{q\ms{L}}$,
\end{itemize}
for all $b\in \ms{B}_{\mu,\nu}$.
This implies that $\phi(G_{\mu,\nu}(b))\subset G$ (cf. \cite[27.1.5]{Lu93}).
\qed

\medskip

Let us introduce some partial orders on $\cP^2$. For $(\mu,\nu), (\zeta,\eta)\in \cP^2$, we first define a partial order $\succeq$ by
\begin{equation*}
    (\mu,\nu) \succeq (\zeta,\eta) \iff \eps_{\mu,\nu}^n - \eps_{\zeta,\eta}^n \in Q_+ \text{ for } n\ge \max(\ell(\mu) + \ell(\nu) , \ell(\zeta) + \ell(\eta)).
\end{equation*}
Note that this condition is independent of $n$.
Next, we define a partial order $\ge$ by
\begin{equation} \label{eq:weaker partial order on cP^2}
    (\mu,\nu) \ge (\zeta,\eta) \iff |\mu| - |\nu| = |\zeta| - |\eta| \text{ and } \mu \supset \zeta, \nu \supset \eta,
\end{equation}
where $\mu\supset\zeta$ means $\mu_i \ge \zeta_i$ for all $i$. Note that $(\mu,\nu) \ge (\zeta,\eta)$ implies $(\mu,\nu) \succeq (\zeta,\eta)$.

Let $\cP(\mu,\nu) = \{(\zeta,\eta)\in \cP^2 \,|\, (\zeta,\eta) \le (\mu,\nu)\}$, which is a finite set. Define $\cP_k(\mu,\nu)$ for $k\ge 0$ inductively by letting $\cP_0(\mu,\nu) = \{(\mu,\nu)\}$ and $\cP_k(\mu,\nu)$ the set of maximal elements in $\cP(\mu,\nu) \setminus \bigsqcup_{l < k} \cP_l(\mu,\nu)$ with respect to $\succeq$.

\begin{prop}\label{prop:filtration on highest tensor lowest}
There exists a sequence of $U_q(\gl_{>0})$-submodules $0=F_{-1} \subset F_0\subset F_1\subset \cdots \subset F_k\subset \cdots \subset V=V_{\mu,\emptyset}\ot V_{\emptyset,\nu}$ such that
\begin{equation*}
\bigcup_{k\ge 1}F_k = V,\quad
F_k/F_{k-1} \cong
\bigoplus_{(\zeta,\eta)\in \cP_k(\mu,\nu)} V_{\zeta,\eta}^{\oplus m^{\mu,\nu}_{\zeta,\eta}},
\end{equation*}
where $m^{\mu,\nu}_{\zeta,\eta}$ is given in \eqref{eq:tensor multiplicity}.
\end{prop}
\pf For $n\ge \max\{ \ell(\mu),\ell(\nu)\}$, consider
\begin{equation*}
V^n\subset V^{n+1}\subset \cdots \subset V,
\end{equation*}
where $V^n=V^n_{\mu,0}\ot V^n_{0,\nu}$ with $V^n_{\mu,0}=U_q(\gl_n)u_{\mu,0}^n$ and $V^n_{0,\nu}=U_q(\gl_n)u_{0,\nu}^n$ (cf. \eqref{eq:Vnmn}).
Let $\ms{B}^n$ be the crystal of $V^n$. We may regard $\ms{B}^n \subset \ms{B}$.

By \eqref{eq:tensor multiplicity}, we have
\[
    \BCrys \cong \bigsqcup_{(\zeta,\eta) \in \cP(\mu,\nu)} \BCrys_{\zeta,\eta}^{\oplus m^{\mu,\nu}_{\zeta,\eta}}.
\]
For $k\ge 0$, let
\[
    \BCrys_k = \bigsqcup_{l\le k} \bigsqcup_{(\zeta,\eta) \in \cP_l(\mu,\nu)} \BCrys_{\zeta,\eta}^{\oplus m^{\mu,\nu}_{\zeta,\eta}} \subset \BCrys.
\]
We inductively construct $F_k$ so that $F_k$ contains $G(b)$ for $b\in \ms{B}_k$.

For $k=0$, there exists a submodule isomorphic to $V_{\mu,\nu}$ by \Cref{lem:embedding of ext wt module-2}, which we denote by $F_0$.
Suppose that we have constructed $F_{k-1}$ for $k\ge 1$.
Let $(\zeta,\eta)\in \cP_k(\mu,\nu)$ be given and let $\ms{B}'\subset \ms{B}$ be a connected component of $\ms{B}$ isomorphic to $\ms{B}_{\zeta,\eta}$.
Choose a sufficiently large $n$ so that $\ms{B}'\cap \ms{B}^n\neq \emptyset$.
Let $b\in \ms{B}'\cap \ms{B}^n$ be a highest weight element.
Then $\wt(b) = \eps_{\zeta,\eta}^n$, which is maximal among the weights of $\BCrys^n \setminus \BCrys_{k-1}$. Hence, $G(b)$ is a highest weight vector in $V^n / (V^n\cap F_{k-1})$.

Suppose that $n' > n$ and $b' \in \BCrys'\cap \BCrys^{n'}$ is a highest weight element of weight $\eps_{\zeta,\eta}^{n'}$.
By the same argument as above, $G(b')$ is a highest weight vector in $V^{n'} / (V^{n'}\cap F_{k-1})$, which is an extremal weight vector with respect to the action of $\UqglPlus[n']$. Since
\[\frac{V^n}{V^n\cap F_{k-1}} \subset \frac{V^{n'}}{V^{n'}\cap F_{k-1}},\]
$G(b)$ is an extremal weight vector in $\UqglPlus[n']G(b')$.
This implies that $G(b)$ is an extremal weight vector in $V/(V\cap F_{k-1})$ with respect to the action of $\UqglPlus$, and $\UqglPlus G(b)$ is isomorphic to $V_{\zeta,\eta}$ by \Cref{thm: ext wt crystal of >0} (2) and \Cref{lem:embedding of ext wt module}.

Let $F'_k$ be the sum of submodules of $V/V\cap F_{k-1}$ corresponding to each connected component $\BCrys'\cong \BCrys_{\zeta,\eta}$ for all $(\zeta,\eta)\in \cP_k(\mu,\nu)$. Then $F'_k \simeq \bigoplus_{(\zeta,\eta) \in \cP_k(\mu,\nu)} V_{\zeta,\eta}^{\oplus m^{\mu,\nu}_{\zeta,\eta}}$. Now we take $F_k = \pi^{-1}(F'_k)$ where $\pi:V\to V/(V\cap F_{k-1})$ is the canonical projection. This completes the induction.
\qed

\begin{cor} \label{cor:identical filtration on lowest tensor highest}
    There exists a filtration on $V_{\emptyset, \nu} \ot V_{\mu, \emptyset}$ with the same property as in \Cref{prop:filtration on highest tensor lowest}. 
\end{cor}
\pf
The functor $-$ in \Cref{rem:bar-involution pullback} preserves weight spaces, and is exact.
Also, we have $\ov{V_{\zeta,\eta}} \cong V_{\zeta,\eta}$ for $(\zeta,\eta)\in \cP^2$, since $\ov{u_{\zeta,\eta}} \in \ov{V_{\zeta,\eta}}$ is also extremal for an extremal weight vector $u_{\zeta,\eta}\in V_{\zeta,\eta}$.
Hence the filtration $\{F_k\}_{k\ge 0}$ on $V_{\mu,\emptyset}\ot V_{\emptyset,\nu}$ in \Cref{prop:filtration on highest tensor lowest} is mapped to a filtration $\{\ov{F}_k\}_{k\ge 0}$ on $\ov{V_{\mu,\emptyset}\ot V_{\emptyset,\nu}}$ such that $\ov{F}_k/\ov{F}_{k-1} \cong F_k/F_{k-1}$.

Now, we apply the following natural isomorphisms
\[
    \ov{V_{\mu,\emptyset}\ot_- V_{\emptyset,\nu}} \cong V_{\mu,\emptyset} \ov{\ot}_{-} V_{\emptyset,\nu} \cong V_{\emptyset,\nu} \ot_+ V_{\mu,\emptyset} \cong V_{\emptyset,\nu} \ot_- V_{\mu,\emptyset}
\]
(see \Cref{rem:$R$ matrix,rem:bar-involution pullback}) to have a filtration on $V_{\emptyset,\nu}\ot V_{\mu,\emptyset}$.
\qed

\begin{rem}\label{rem:R-matrix on lohimod}
    The R-matrix $R = R_-$ induces a well-defined $\UqglPlus$-linear isomorphism $V_{\emptyset,\nu}\ot V_{\mu,\emptyset} \cong V_{\mu,\emptyset}\ot V_{\emptyset,\nu}$, whose inverse $R^{-1}$ is also well-defined. We may also deduce \Cref{cor:identical filtration on lowest tensor highest} from this isomorphism.
\end{rem}

\newcommand{\CatC}{\mathcal{C}}

Let $\td{\CatC}$ be a category of $\UqglPlus$-modules with weight space decomposition,
and let $\CatC$ be the full subcategory of $\UqglPlus$-modules of finite length whose simple subquotients are $V_{\lambda,\mu}$ for $(\lambda,\mu)\in \cP^2$,
that is, the Serre subcategory generated by $V_{\lambda,\mu}$'s.
\begin{prop}
The category $\CatC$ is a monoidal subcategory of $\td{C}$, whose Grothendieck ring $K(\CatC)$ is commutative.
\end{prop}
\pf 
It is enough to show that $\CatC$ is closed under tensor product.
By \Cref{prop:filtration on highest tensor lowest} and \Cref{cor:identical filtration on lowest tensor highest},
we have $V_{\mu,\emptyset}\ot V_{\emptyset,\nu}, V_{\emptyset,\nu}\ot V_{\mu,\emptyset} \in \CatC$ for $\la,\mu\in \cP$.
Then we have $V_{\mu,\nu}\ot V_{\zeta,\emptyset}\in \CatC$ for $\zeta\in \cP$, since $V_{\mu,\nu} \subset V_{\emptyset,\nu} \ot V_{\mu,\emptyset}$ and $V_{\mu,\emptyset}\ot V_{\zeta,\emptyset}$ is a direct sum of $V_{\nu,\emptyset}$'s.
Similarly, $V_{\mu,\nu}\ot V_{\emptyset,\zeta} \subset (V_{\mu,\emptyset}\ot V_{\emptyset,\nu}) \ot V_{\emptyset,\eta}\in \CatC$ for $\eta\in \cP$.

Now for $V\in \CatC$, we have $V\ot V_{\mu,\emptyset}, V\ot V_{\emptyset,\nu} \in \CatC$, which implies that $V\ot V_{\mu,\nu}\subset V\ot V_{\mu,\emptyset}\ot V_{\emptyset,\nu}\in \CatC$. Similarly, we have $V_{\mu,\nu}\ot V \in \CatC$. Therefore $\CatC$ is closed under tensor product.
The commutativity of $K(\CatC)$ follows from the same argument as in the proof of \Cref{cor:identical filtration on lowest tensor highest}.
\qed

\section{Fock space $\FFn$}
\label{sec:Fock space F^n}
\subsection{$R$-matrix and $q$-deformed exterior and symmetric algebras}\label{subsec:q-deformed exterior and symmetric algebras}
We introduce a uniform construction of $q$-deformed exterior or symmetric algebras using $R$-matrices, which carries a commuting action of two quantum groups.

Let $U_q(\g)$ be the quantized enveloping algebra as in \Cref{subsec:quantum group}.
Let $p$ be another formal variable.
Let $\dot{A}=(\dot{a}_{ij})_{i,j\in \dot{I}}$ be another symmetrizable generalized Cartan matrix indexed by $\dot{I}$, and let $U_p(\dot{\g})$ be the quantized enveloping algebra associated with a Cartan datum $(\dot{A},\dot{P},\dot{P}^\vee,\dot{\Pi},\dot{\Pi}^\vee, (,))$ over $\Q(p)$.

Let $V$ and ${W}$ be $U_q(\g)$ and $U_p(\dot{\g})$-modules such that the universal $R$ matrices $R:=R^{\rm univ}$ and $\dot{R}:=\dot{R}^{\rm univ}$ for $U_q(\g)$ and $U_p(\dot{\g})$ yield well-defined maps on $V\ot V$ and ${W}\ot {W}$, respectively.
Here the comultiplications for $U_q(\g)$ and $U_p(\dot{\g})$ are assumed to be any $\ot_*$ in \Cref{rem:$R$ matrix}.

For $k\in \Z_{>0}$, we define
\begin{equation}\label{eq:symmetric algebra}
 \mc{A}^k(V\ot {W}) = (V\ot {W})^{\ot k} \Bigg/ \sum_{i=1}^{k-1} {\rm Im}\left( R_{i,i+1} - \dot{R}_{i,i+1} \right),
\end{equation}
where $R_{i,i+1}$ denotes $R$ acting on the $i$th and $(i+1)$st component $V\ot V$, and $\dot{R}_{i,i+1}$ is defined in the same way. It is a $U_q(\g)\ot U_p(\dot{\g})$-module.
Let
\begin{equation*}
\mc{A}(V\ot {W})=\bigoplus_{k\ge 0}\mc{A}^k(V\ot {W}),\quad  \mc{A}^0(V\ot {W})=\Q(q)\ot \Q(p),
\end{equation*}
which is a $\Z_{\ge 0}$-graded $\Q(q)\otimes \Q(p)$-algebra.

When $p=q$ (resp. $p=-q^{-1}$), as we shall see, $\mc{A}(V\ot {W})$ behaves like a $q$-deformed symmetric (resp. exterior) algebra generated by $V\ot {W}$, so we denote it by $S(V\ot {W})=\bigoplus_{k\ge 0}S^k (V\ot {W})$ and $\tbgwed(V\ot {W})=\bigoplus_{k\ge 0}\tbgwed^k (V\ot {W})$, respectively.
In these cases, the algebras $\mc{A}(V\ot W)$ will be understood as algebras over $\Q(q)$, by tensoring with $\Q(q)\ot \Q(p) / (p-q)$ (resp. $\Q(q)\ot \Q(p) / (p+q^{-1})$).

\begin{rem}

When $\lie{g} = \dot{\lie{g}} = \lie{gl}_n$ and $V$ and $W$ are standard representations, we recover the quantized coordinate ring of $GL_n$, see \cite{B} and references therein.
See \cite{BZ} for an analogous construction using only $\Uqg$ and $V$ without $\QUE{\dot{\lie{g}}}$ and $W$.
Although it is not needed in this paper, we expect that it can be extended to the case of quantum affine algebras analogously using normalized $R$-matrices.
When $W = \Q(q)$ with $\dot{R} = -1$, this will give an exterior algebra introduced in \cite{KMPY}.
We also refer the reader to \cite{U} for a $q$-deformed exterior algebra with a $\QUE{\widehat{\lie{gl}}_m} \ot \QUE{\widehat{\lie{gl}}_n}[][-q^{-1}]$-action constructed by using affine Hecke algebras.
\end{rem}

\begin{lem}\label{lem:application of hexagon identity}
    For $x\in \mc{A}^k(V\ot W)$ and $y\in \mc{A}^l(V\ot W)$, we have $R(x \ot y) = \dot{R}(x \ot y)$ in $\mc{A}^{k+l}(V\ot W)$.
\end{lem}
\pf
    By the hexagon property of $R$ (cf. \cite[Proposition 32.2.4]{Lu93}), $R(x\ot y)$ is given by a series of application of $R_{i, i+1}$'s on $x\ot y$. The same holds for $\dot{R}(x\ot y)$, and the action of $R_{i,i+1}$ agrees with that of $\dot{R}_{i,i+1}$, so the equality follows.
\qed
\medskip

Suppose that $V$ and $W$ have $\Q$-linear involutions $-$ compatible with the actions of $\Uqg$ and $\QUE{\dot{\lie{g}}}[][p]$, respectively.
We define an involution $-$ on $(V\ot W)^{\ot k}$ using $\Theta\dot{\Theta}'$ following \eqref{eq:bar involution on tensor product} and \eqref{eq:bar involution on n-tensor product}, where $\Theta$ is the quasi-$R$-matrix and $\dot{\Theta}'$ is the one associated with $\dot{R}^{-1}$ (see \eqref{eq:opposite quasi-$R$-matrix}), with respect to given comultiplications for $\Uqg$ and $\QUE{\dot{\lie{g}}}[][p]$.

\begin{lem} \label{lem:bar-involution descends to exterior algebra}
    For $k\ge 0$, the involution $-$ on $(V\ot W)^{\ot k}$ induces a well-defined involution on $\mc{A}^k(V\ot W)$, which we still denote by $-$.
\end{lem}
\pf We assume that the coproduct of $U_q(\lie{g})$ and $U_p(\dot{\lie{g}})$ are $\Delta_+$, since the proofs for the other cases are similar. Since $\Theta^{(k)} = (1^{\ot (i-1)} \ot \Delta \ot 1^{\ot (k-i-1)})( \Theta^{(k-1)} ) \Theta^{(i,i+1)}$, where $\Theta^{(i,i+1)}$ denotes $\Theta$ acting on the $(i,i+1)$ compononents, it suffices to show that $\Theta\dot{\Theta}' \ov{\left( \Image (R-\dot{R}) \right)} \subset \Image(R-\dot{R})$. It follows from the identity $\Theta\dot{\Theta}'(\ov{R} - \dot{\ov{R}}) = -(R - \dot{R})R^{-1}\dot{R}^{-1}\Theta\dot{\Theta}'$ due to \Cref{lem:RThetaRTheta'}.
\qed
\medskip

We denote by $\wt$ and $\dot{\wt}$ the weights for $\Uqg$ and $\QUE{\dot{\lie{g}}}[][p]$, respectively.
    We call $x \in \mc{A}^k(V\ot W)$ a weight vector if
    $q^h x = q^{\langle h, \wt(x) \rangle} x$ and $p^{\dot{h}} x = p^{\langle \dot{h}, \dot{\wt}(x) \rangle} x$ for all $h \in P$ and $\dot{h} \in \dot{P}$.
We have the following formula for the involution.

\begin{lem}\label{lem:bar-involution reverses order}
    For weight vectors $x\in \mc{A}^k(V\ot W)$ and $y\in \mc{A}^l(V\ot W)$,
    \begin{equation}\label{eq:bar-involution reverses order}
        \ov{x \cdot y} = q^{\varepsilon (\wt(x), \wt(y))} p^{\dot{\varepsilon} (\dot{\wt}(x), \dot{\wt}(y))} \ov{y} \cdot \ov{x},
    \end{equation}
    where $x\cdot y$ denotes the product in $\mc{A}(V\ot W)$, and the sign $\varepsilon$ is 
    $+$ (resp. $-$) if the coproduct on $\Uqg$ is $\ot_+$ (resp. $\ot_-$) or $\ov{\ot}_{-}$ (resp. $\ov{\ot}_+$). The sign $\dot{\varepsilon}$ is determined by the coproduct on $\QUE{\dot{\lie{g}}}[][p]$ in the same way.
\end{lem}
\pf We assume that the coproduct of $U_q(\lie{g})$ and $U_p(\dot{\lie{g}})$ are $\Delta_+$, since the proofs for the other cases are similar.
\begin{equation}\label{eq:proof of bar-involution reverses order}
\begin{aligned}
    \Theta\dot{\Theta}' (\ov{x}\cdot \ov{y}) &= R^{-1}\Pi\sigma \dot{R}\dot{\Pi}\dot{\sigma} (\ov{x}\cdot \ov{y}) \\
    &= q^{(\wt(x), \wt(y))} p^{(\dot{\wt}(x), \dot{\wt}(y))} R^{-1}\dot{R}(\ov{y}\cdot \ov{x}) \\
    &= q^{(\wt(x), \wt(y))} p^{(\dot{\wt}(x), \dot{\wt}(y))} \ov{y}\cdot \ov{x},
\end{aligned}
\end{equation}
where the last equality follows from \Cref{lem:application of hexagon identity}.
\qed

\begin{lem}
    The involution $-$ on $\mc{A}(V\ot W)$ coincides with that induced from $\Theta'\dot{\Theta}$.
\end{lem}
\pf
    By similar arguments as in the proof of \Cref{lem:bar-involution descends to exterior algebra}, we see that $\Theta'\dot{\Theta}$ induces a well-defined bar-involution on $\mc{A}(V\ot W)$.
    We have an analogue of \eqref{eq:proof of bar-involution reverses order}, where $\Theta$, $\dot{\Theta}'$, $R^{-1}$, and $\dot{R}$ are replaced by $\Theta'$, $\dot{\Theta}$, $R$, and $\dot{R}^{-1}$. This shows that \Cref{lem:bar-involution reverses order} also holds for the bar-involution induced from $\Theta'\dot{\Theta}$.
\qed

\subsection{$q$-deformed exterior algebras of type $A$}\label{subsec:q-ext alg of type A}

Suppose that $S$ and $T$ are intervals of $\mathbb{R}$, not necessarily bounded. Let $\QUE{\lie{gl}_S}$ and $\QUE{\lie{gl}_T}[][p]$ be as in \Cref{sec:gln} with $p = -q^{-1}$. Denote $\iset{S} = S\cap \Z$ and $\iset{T} = T\cap \Z$.
\begin{equation}\label{eqref:our choice}
    \parbox{0.9\textwidth}{
        From now on, we assume that the comultiplications for $\QUE{\lie{gl}_S}$ and $\QUE{\lie{gl}_T}[][p]$ are $\Delta_-$ and $\ov{\Delta}_+$, respectively.
    }
\end{equation}
Let $V_S$ be the natural representation of $U_q(\gl_S)$ with basis $\{\,v_a\,|\,a\in \iset{S} \,\}$, where $e_{i}v_a=\de_{i+1\, a}v_{i}$, $f_{i}v_a=\de_{ia}v_{i+1}$ and $q^{E_{ii}}v_a = q^{\langle E_{ii},\e_a \rangle}v_a$ for $i\in \dset{S}$. Let $\dot{V}_T$ be the natural representation of $U_p(\gl_T)$ with basis $\{\,\dot{v}_b\,|\,b\in \iset{T}\,\}$.
Let
\begin{equation*} 
\tbgwed_{S,T}=\bigoplus_{k\ge 0}\tbgwed^k_{S,T},\quad\text{where} \quad
\tbgwed^k_{S,T}=\tbgwed^k(V_S\ot\dot{V}_T).
\end{equation*}
It is well-known \cite{J} that $R$ acts on $V_S\ot V_S$ by
\begin{equation}\label{eq:R matrix for type A natural}
R(v_a\ot v_b)=
\begin{cases}
q v_a\ot v_b & \text{if $a=b$},\\
v_b\ot v_a + (q-q^{-1}) v_a\ot v_b& \text{if $a<b$},\\
v_b\ot v_a  & \text{if $a>b$},
\end{cases}
\end{equation}
for $a,b\in \iset{S}$, and the formula for $\dot{R}$ on $\dot{V}_T\ot\dot{V}_T$ is given by applying the flipping $\sigma$ to \eqref{eq:R matrix for type A natural}.

Let $w_{(a,b)}$ denote the image of $v_a\ot \dot{v}_b$ in $\bigwedge^1_{S,T}$ for $(a,b)\in \iset{S}\times \iset{T}$, and
$w_{(a_1,b_1)}\wedge\dots\wedge w_{(a_k,b_k)}$ denote the image of $(v_{a_1}\ot \dot{v}_{b_1})\ot \cdots\ot(v_{a_k}\ot \dot{v}_{b_k})$ in $\bigwedge^k_{S,T}$.
By \eqref{eq:R matrix for type A natural}, we have the following relations in $\bigwedge^2_{S,T}$:
\begin{equation}\label{eq:relations for q-wedge}
\begin{split}
w_{(a,b)}\wedge w_{(c,d)}
=
\begin{cases}
0 & \text{if $(a,b)=(c,d)$},\\
-q w_{(c,b)}\wedge w_{(a,d)} & \text{if $a>c$ and $b=d$},\\
q^{-1} w_{(a,d)}\wedge w_{(c,b)} & \text{if $a=c$ and $b<d$},\\
w_{(c,d)}\wedge w_{(a,b)} & \text{if $a<c$ and $b<d$},\\
w_{(c,d)}\wedge w_{(a,b)}-(q-q^{-1})w_{(c,b)}\wedge w_{(a,d)} & \text{if $a>c$ and $b<d$},
\end{cases}
\end{split}
\end{equation}
for $a,a'\in \iset{S}$ and $b,b'\in \iset{T}$.
It follows from \eqref{eq:symmetric algebra} that $\bigwedge_{S,T}$ is the $\Q(q)$-algebra generated by $w_{(a,b)}$ for $(a,b)\in \iset{S}\times\iset{T}$ subject to \eqref{eq:relations for q-wedge}.

By \Cref{lem:bar-involution reverses order}, the involution $-$ on $\bigwedge^k_{S,T}$ is given as follows:
\begin{equation}\label{eq:bar-involution on type A exterior algebras}
    \ov{w_{(c_1,d_1)} \wedge \cdots \wedge w_{(c_k, d_k)}} = q^{-l(w_{\bm{c}})}p^{-l(w_{\bm{d}})} w_{(c_k,d_k)} \wedge \cdots \wedge w_{(c_1,d_1)},
\end{equation}
where $w_{\bm{c}}$ is the element in the symmetric group $S_k$ of maximal length among the ones fixing $k$-tuple of integers $\bm{c} = (c_1,\ldots, c_k)$, and $w_{\bm{d}}$ is defined for $\bm{d}= (d_1,\ldots, d_k)$ similarly.

We define a total order $<$ and $<'$ on $\iset{S}\times \iset{T}$ by
\begin{equation*}
\begin{split}
\text{$(a,b)< (c,d)$ if and only if ($b<d$) or ($a> c$ and $b=d$),} \\
\text{$(a,b)<' (c,d)$ if and only if ($a>c$) or ($a=c$ and $b<d$).}
\end{split}
\end{equation*}
Let 
\begin{equation*}
    B_{S,T} =\left\{\,M=(m_{ab})\,|\, m_{ab}\in \{0, 1\} \  (a\in \iset{S}, b\in \iset{T}), \, \textstyle{\sum_{a,b}} m_{ab} < \infty  \, \right\}.
\end{equation*}
For $M=(m_{ab})\in B_{S,T}$, we put
\[ w_{M}=\vec{\bigwedge}\, w^{m_{ab}}_{(a,b)},\] 
where $\vec{\bigwedge}$ denotes the wedge product over $\iset{S}\times \iset{T}$ with respect to the total order $<$ (that is, we read the entries in $M$ column by column from left to right, and then from bottom to top in each column). By \eqref{eq:relations for q-wedge}, we have $w_M=w'_M$, where $w'_M$ is the product defined with respect to $<'$  (that is, we read the entries in $M$ row by row from bottom to top, and then from left to right in each row).
Then the set of standard monomials $\{\,w_M\,|\,M\in B_{S,T}\,\}$ is a $\Q(q)$-linear basis of $\bigwedge_{S,T}$ by standard arguments \cite{Be}.

\begin{ex}
Consider the case when $S = [1,m]$ and $T = [1,n]$ for $m,n\in \Z_{> 0}$. We denote $\bigwedge_{S,T}$ by $\tbgwed_{m,n}$, and $B_{S,T}$ by $B_{m,n}$.
Let $M\in B_{m,n}$ be given. For $1\le a\le m$ and $1\le b\le n$, let $M_a$ and $M^b$ be the $a$-th row and $b$-th column of $M$, respectively. Then
\begin{equation*}
w_M=w_{M^{1}}w_{M^{2}}\dots w_{M^n}=w_{M_{m}}w_{M_{m-1}}\dots w_{M_1},
\end{equation*}
where we regard $M_a$ (resp. $M^b$) as an element in $B_{m,n}$ which is equal to $M$ in the $a$-th row (resp. $b$-th column) and zero elsewhere.
Then we have an isomorphism of $U_q(\gl_m)$-modules
\begin{equation}\label{eq:iso U_q}
\begin{split}
\xymatrixcolsep{3pc}\xymatrixrowsep{0pc}\xymatrix{
 \bigwedge_{m,n} \ \ar@{->}[r] &\ \bigwedge_{m,1}^{\ot n}  \\
 w_M \ \ar@{|->}[r] &\  w_{M^{1}}\ot \dots\ot w_{M^n}
},
\end{split}
\end{equation}
and also an isomorphism of $U_p(\gl_n)$-modules
\begin{equation}\label{eq:iso U_p}
\begin{split}
\xymatrixcolsep{3pc}\xymatrixrowsep{0pc}\xymatrix{
 \bigwedge_{m,n} \ \ar@{->}[r] &\ \bigwedge_{1,n}^{\ot m}  \\
 w_M \ \ar@{|->}[r] &\  w_{M_{m}}\ot \dots\ot w_{M_1}
}.
\end{split}
\end{equation}
\end{ex}
\medskip

Let us construct a crystal base of $\tbgwed_{S,T}$.
Let $\te_i$ and $\tf_i$ be the lower crystal operators \eqref{eq:Kashiwara operator lower} on $\tbgwed_{S,T}$ with respect to the action of $U_q(\gl_S)$ \eqref{eq:Kashiwara operator lower} for $i\in\dset{S}$.

To define crystal operators for the $\gl_T$, let us twist the action of $U_p(\lie{gl}_T)$ by the isomorphism of $\Q$-algebras
\begin{equation}\label{eq:pull back Uqgln}
    \psi : U_q(\gl_T)\longrightarrow U_{p}(\gl_T)
\end{equation}
such that $\psi(e_i)=e_i$, $\psi(f_i)=f_i$, $\psi(q^{h})=p^{-h}$ and $\psi(q)=p^{-1}$ for $i$ and $h$.
Then we define $\dot{\te}_j$ and $\dot{\tf}_j$  for $j \in \dset{T}$ to be the upper crystal operators \eqref{eq:Kashiwara operator upper} with respect to this action of $\QUE{\gl_T}$.
Since $(\ov{\Delta}_+)^\psi := (\psi^{-1}\ot \psi^{-1}) \circ (\ov{\Delta}_+)\circ \psi = \Delta_+$, the upper crystal lattice for $U_q(\gl_T)$ is compatible with \eqref{eqref:our choice}.
Note that $\dot{\te}_j$, $\dot{\tf}_j$ commute with $\te_i$, $\tf_i$,
although \eqref{eq:pull back Uqgln} does not extend to a homomorphism of $\Q$-algebras $U_q(\lie{gl}_m)\ot_{\Q(q)} U_q(\lie{gl}_n) \to U_q(\lie{gl}_m)\ot_{\Q(q)} U_p(\lie{gl}_n)$.

Let
\begin{equation*}
\ms{L}(\tbgwed_{S,T})= \bigoplus_{M\in B_{S,T}}A_0 w_M,\quad
\ms{B}(\tbgwed_{S,T})= \{\,w_M \!\!\!\pmod{q \ms{L}(\tbgwed_{S,T})}\,|\,M\in B_{S,T}\,\}.
\end{equation*}

\begin{prop}\label{prop:crystal base of a wedge}
The pair $(\ms{L}(\tbgwed_{S,T}),\ms{B}(\tbgwed_{S,T}))$ is a lower crystal base of $\tbgwed_{S,T}$ as a $U_q(\gl_S)$-module and also upper crystal base of $\tbgwed_{S,T}$ as a $U_q(\gl_T)$-module with respect to \eqref{eq:pull back Uqgln}.
\end{prop}
\pf The proof is essentially the same as in \cite[Propositions 4.2, 4.3]{K14}.
It suffices to consider the case when $S$ and $T$ are bounded, since any element of $\tbgwed_{S,T}$ is contained in a subspace $\tbgwed_{S',T'}$ for bounded intervals $S'$ and $T'$. So we may assume $S = [1,m]$ and $T = [1,n]$.

First, it is clear that $(\ms{L}(\tbgwed_{m,1}),\ms{B}(\tbgwed_{m,1}))$ is a lower crystal base of $\tbgwed_{m,1}$ as a $U_q(\gl_m)$-module, and hence so is $(\ms{L}(\tbgwed_{m,n}),\ms{B}(\tbgwed_{m,n}))$ is a lower crystal base of $\tbgwed_{m,n}$ as a $U_q(\gl_m)$-module by \eqref{eq:iso U_q}.

Next, it is also clear that $(\ms{L}(\tbgwed_{1,n}),\ms{B}(\tbgwed_{1,n}))$ is a upper crystal base of $U_q(\gl_n)$ with respect to $\dot{\te}_j$ and $\dot{\tf}_j$ for $j=\dset{T}$ since $\tbgwed_{1,n}$ is a direct sum of fundamental representations. Hence $(\ms{L}(\tbgwed_{m,n}),\ms{B}(\tbgwed_{m,n}))$  is also a upper crystal base of $\tbgwed_{m,n}$ as a $U_q(\gl_n)$-module by \eqref{eq:iso U_p}.
\qed\medskip

\begin{ex}
We may identify $\ms{B}(\tbgwed_{m,n})$ with $B_{m,n}$. Note that for $M=(m_{a1})_{1\le a \le m}\in B_{m,1}$ and $i$, we have $\te_i M = M + {\bf e}_i - {\bf e}_{i+1}$ and $\tf_i M = M - {\bf e}_i + {\bf e}_{i+1}$, where ${\bf e}_i$ denotes the standard basis of $\Z^m$ and $\td{x}_i M$ ($x = e,f$) is assumed to be zero if it does not belong to $B_{m,1}$. Then we apply \eqref{eq:tensor product rule} to describe $\te_i$ and $\tf_i$ on $B_{m,n}=B_{m,1}^{\ot n}$. We have a similar description of $\dot{\te}_j$ and $\dot{\tf}_j$ on $B_{m,n}$.
\end{ex}

\begin{rem}
    There is an affine analogue of $\bigwedge_{m,n}$ in \cite[Section 3.3]{U} 
    which is a $U'_q(\widehat{\lie{sl}}_n)\otimes U'_p(\widehat{\lie{sl}}_l)$-module algebra of level 0 defined in terms of affine Hecke algebras. Also, its bar-involution is similar to the one here.
    We also expect to recover the higher level Fock space in \cite{U} following the construction \eqref{eq:symmetric algebra}.
    The algebra $\bigwedge_{m,n}$ also appears as the subalgebra of the negative half of the quantum superalgebra associated with $\gl_{m|n}$ generated by odd root vectors \cite[Section 2.4]{K14}, where the linear order on $[n]$ is reversed.
\end{rem}

Define a $\Q(q)$-bilinear form on $\bigwedge_{S,T}$ by 
\begin{equation}\label{eq:q-Shapovalov form on Fn}
    \langle w_M, w_{M'} \rangle = \de_{M,M'} \quad (M,N'\in B_{S,T}).
\end{equation}
Then for $x,y\in \bigwedge_{S,T}$ and $u\in \QUE{\lie{gl}_S}, \dot{u}\in \QUE{\lie{gl}_T}[][p]$, we have
    \begin{equation} \label{eq:antiauto relation of q-Shapovalov form}
      \langle ux,y\rangle = \langle x, \tau(u)y\rangle, \quad \langle x, \dot{u}y\rangle = \langle \tau(\dot{u})x, y\rangle.
    \end{equation}
(cf. \cite[Proposition 5.7]{U}),
This also can be seen from the observation that \eqref{eq:antiauto relation of q-Shapovalov form} holds for the case of a single-column or single-row, and that $(\tau\otimes \tau )\circ \Delta = \Delta \circ \tau$.

\subsection{Semi-infinite limit and the Fock space $\mc{F}^n$}\label{sec:semi-infinite limit}

Suppose that $S$ and $T$ are intervals with $T$ bounded. Consider a directed system of $\Q(q)$-spaces \[ \{\, \tbgwed_{S', T} \,|\, S' \subset S, S'\text{ is an interval which is bounded below} \, \}, \] with respect to the maps $\psi_{S',S''}$ for $S' \subset S''$ with $(S'' \setminus S') \cap \Z = [a,b]\cap \Z$ for some $a,b$, which are defined by

\begin{equation}\label{eq:directed system}
    \begin{split}
    \xymatrixcolsep{3pc}\xymatrixrowsep{0pc}\xymatrix{
     \bigwedge_{S',T} \ \ar@{->}[r]^{\psi_{S',S''}} &\ \bigwedge_{S'',T}  \\
     w \ \ar@{|->}[r] &\ w \wedge w_{[a,\min(b,0)]\times T}
    }
\end{split}
\end{equation}
Here $[a,b]\times T$ denotes the matrix in $B_{S'',T}$ with $1$'s in the rows from $a$ to $b$ and $0$'s elsewhere. Since $\iset{T}$ is finite, $w_{[a,b]\times T}$ is well-defined.
In particular, this sends $w_M$ ($M\in B_{S',T}$) to another standard monomial $w_{M'}$, where $M'\in B_{S'',T}$ is obtained by extending $M$ with $0$'s on the rows with positive indices, and $1$'s on the rows with non-positive indices.

Let $\FF_{S,T}$ be the limit of \eqref{eq:directed system}, and $\psi_{S'}: \bigwedge_{S',T} \longrightarrow \FF_{S,T}$ the canonical embedding.
The map in \eqref{eq:directed system} is $U_q(\gl_{S'})\ot U_p(\lie{sl}_T)$-linear, and hence $\mc{F}_{S,T}$ becomes a well-defined $U_q(\gl_S)\ot U_p(\lie{sl}_T)$-module. 
Let $w \in \FF_{S,T}$ be given such that  $w$ is an image of a weight vector $w' \in \bigwedge_{S',T}$ for some interval $S' \subset \R$ that is bounded below.
We define $\dot{\wt}(w) = \dot{\wt}(w') - (\min (S' \cap \Z) - 1)\sum_{i=1}^n \eps_i$,
which is independent of $S'$.
Hence, $\FF_{S,T}$ becomes a $U_q(\gl_S)\ot U_p(\gl_T)$-module with respect to $\wt$ and $\dot{\wt}$.

\begin{lem}
    We have
    $\ov{\psi_{S',S''}(w)} = \psi_{S',S''}(\ov{w})$ for $w\in \bigwedge_{S',T}$.
\end{lem}
\pf It can be shown either by combinatorial arguments using \eqref{eq:bar-involution on type A exterior algebras}, or by using quasi-$R$-matrices as follows. Let $w_\circ = w_{[a,\min(b,0)]\times T}$. It is easy to see that $w_\circ$ is bar-invariant.
By \Cref{lem:bar-involution descends to exterior algebra}, we have
\[
    \ov{w \wedge w_\circ} = \Theta_{\ot_-} \dot{\Theta}'_{\ov{\ot}_+}(\ov{w} \wedge \ov{w_\circ}).
\]
Note that $\Theta_{\ot_-} = \Pi^{-1}\Theta\Pi^{-1}$.
First, $\Pi(\ov{w} \wedge w_\circ) = \ov{w}\wedge w_\circ$, since $U_q(\gl_{S''})$-weights of $\ov{w}$ and $\ov{w_\circ}$ are orthogonal.
Since $\Theta$ is a formal sum of elements of $\QUE{\gl_{S''}}[-]^-\ot \QUE{\gl_{S''}}[+]$, and $\QUE{\gl_{S''}}[+]$ acts trivially on $\ov{w_\circ}$,
$\Theta(\ov{w}\wedge w_{\circ}) = \ov{w}\wedge w_{\circ}$.
Similarly, $\dot{\Theta}'_{\ov{\ot}_+} = \dot{\ov{\Theta}}$ acts trivially $\ov{w}\wedge w_{\circ}$. Therefore, $\ov{w\wedge w_\circ}=\ov{w}\wedge w_\circ$.
\qed
\medskip

Then, we have the following.

\begin{prop}\label{prop:bar-involution on Fock space}
    $\FF_{S,T}$ has a $\Q$-linear involution $-$ such that $\ov{\psi_{S'}(w)} = \psi_{S'}(\ov{w})$ for $w\in \bigwedge_{S',T}$.
\end{prop}

Let
\begin{equation*}
F_{S,T} =\left\{\,M=(m_{ab})\,|\,a\in \iset{S}, b\in \iset{T},\,\text{$m_{ab}=0$ for $a\gg 0$, and $m_{ab}=1$ for $a\ll 0$}\, \right\}.
\end{equation*}
For $M=(m_{ab})\in F_{S,T}$, we define $w_M$ to be a unique element in $\FF_{S,T}$ such that
\begin{equation*}
 w_M = \psi_{S'}(w_{M_{S'}}),
\end{equation*}
with $M_{S'}=(m_{ab})_{a\in \iset{S'}, b\in \iset{T}}$ for some interval $S'$ which is bounded below such that $m_{ab}=1$ ($a \le \min{S'}$, $b\in \iset{T}$). 
It immediately follows from the construction that
$\{\,w_M\,|\,M\in F_{S,T}\,\}$ is a $\Q(q)$-basis of $\mc{F}_{S,T}$,

Let
\begin{equation}\label{eq:crystal base of Fock space}
\ms{L}(\mc{F}_{S,T})= \bigoplus_{M\in F_{S,T}}A_0 w_M,\quad
\ms{B}(\mc{F}_{S,T})= \{\,w_M \!\!\!\pmod{q \ms{L}(\mc{F}_{S,T})}\,|\,M\in F_{S,T}\,\}.
\end{equation}

\begin{prop}
The pair $(\ms{L}(\mc{F}_{S,T}),\ms{B}(\mc{F}_{S,T}))$ is a lower crystal base of $\mc{F}_{S,T}$ as a $U_q(\gl_S)$-module and also upper crystal base as a $U_q(\lie{gl}_T)$-module with respect to \eqref{eq:pull back Uqgln}.
\end{prop}
\pf Let $S',S''$ be intervals bounded below with $S'\subset S'' \subset S$. Then $\psi_{S',S''}$ in \eqref{eq:directed system} preserves the crystal bases, and $(\ms{L}(\tbgwed_{S',T}),\ms{B}(\tbgwed_{S',T}))$ forms a directed system, whose limit is equal to  $(\ms{L}(\mc{F}_{S,T}),\ms{B}(\mc{F}_{S,T}))$. Hence it is a lower (resp. upper) base of $\mc{F}_{S,T}$ as a $U_q(\gl_S)$-module (resp. $U_q(\gl_T)$-module).
\qed\medskip

Note that $\mc{F}_{S,n}:=\FF_{S,[1,n]}$ is isomorphic to $\mc{F}_{S,\{1\}}^{\ot n}$ as a $U_q(\gl_S)$-module for $n\in \Z_{>0}$.
Let
\[
    \FFn = \FF_{\R,n}.
\]
Put $\mc{F} = \mc{F}^1$.
The $U_q(\gl_\infty)\ot U_p(\gl_n)$-crystal $\BCrys(\FFn)$ (that is, a crystal with respect to $\tilde{x}_i,\dot{\tilde{y}}_j$ for $i\in \Z, 1\le j\le n-1$) can be identified with $F_{\R,n}$, and it is isomorphic to the one introduced in \cite[Section 6]{K09}. In particular, we have a crystal isomorphism
\begin{equation}\label{eq:level n Fock space crystal decomposition}
    \BCrys(\FFn)\cong \bigsqcup_{\la\in \Z_+^n} \ms{B}(\La_\la)\times \ms{B}(\doteps_\lambda),
\end{equation}
(\cite[Theorem 6.1]{K09}), where $\ms{B}(\La_\la)$ and $\ms{B}(\doteps_\la)$ are the crystals of the highest weight modules $V(\La_\la)$ and $V(\doteps_\la)$ of $U_q(\gl_\infty)$ and $U_q(\gl_n)$ with highest weights $\La_\la$ and $\doteps_\la$, respectively.
Here $\doteps_\la$ means the $\eps_\la$ with respect to $\dot{\wt}$.
The corresponding decomposition of $\FFn$ is a quantum analogue of the level-rank duality, whose non-quantum version is originally due to \cite{Fr}:
\begin{equation}\label{thm:level n Fock space decomposition}
    \FFn \cong \bigoplus_{\la\in \Z_+^n} V(\La_\la)\otimes V(\doteps_\la).
\end{equation}
Indeed, for $\lambda\in \Z_+^n$, let $M(\lambda) \in \BFn$ be the matrix whose $(i,j)$ entry is given by
\begin{equation}\label{eq:highest weight element in FFn}
    M(\lambda)_{ij} = \begin{cases}
            1 & i \le \lambda_j, \\
            0 & \text{otherwise.}
    \end{cases}
\end{equation}
Then $M(\lambda)$ is the highest weight element in the connected component of $\BFn$ isomorphic to $\BCrys(\Lambda_\lambda)\ot \BCrys(\doteps_\lambda)$.
The corresponding monomial $w_{M(\lambda)}$ is a highest weight vector, and generates a submodule isomorphic to $V(\Lambda_\lambda)\ot V(\doteps_\lambda)$.
Since $w_{M(\lambda)}$'s generate $\BFn$, they also generate $\FFn$, and hence \eqref{thm:level n Fock space decomposition} follows.

Note that the bilinear form \eqref{eq:q-Shapovalov form on Fn} induces a bilinear form on $\FF_{S,T}$,
which restricts to the $q$-Shapovalov form \eqref{eq:q-Shapovalov form on uqg} on each irreducible component $V(\La_\la)\ot V(\doteps_\la)$ in case of $\FFn$.

Let us end this subsection with some notations which will be used in later sections.

For $M = (m_{ij})\in \BFn$, we define the charge of $M$ to be $\#\{ (i,j) \,|\, m_{ij}=1 \ (i\ge 1)\} - \#\{ (i,j) \,|\, m_{ij} = 0\  (i\le 0) \}$.
For $k\in \Z$, the $\Q(q)$-span of $w_M$ for $M\in \BFn$ of charge $k$ is invariant under the action of $\Uqgl \ot \UpglPlus[n]$.
Put $|k\rangle = w_{M((k))} \in \FF$, which is a $\Uqgl$-highest weight vector of weight $\La_k$.

For $(\mu,\nu)\in \cP^2$ with $\ell(\mu)+\ell(\nu)\le n$, put
\begin{equation}\label{eq:integer partitions to two partitions}
\Lambda^n_{\mu,\nu} = \Lambda_\la, \quad M^n(\mu,\nu) = M(\la),
\end{equation}
where $\la$ is the unique element in $\Z^n_+$ such that $\mu = (\max\{ \la_1,0\}, \max\{\la_2,0\}, \cdots)$ and $\nu = (- \min\{\la_n,0\}, -\min\{\la_{n-1},0\},\cdots)$.
We also let
\begin{equation*}
\La_{\mu,\nu} = \La^n_{\mu,\nu} - (n-\ell(\mu) - \ell(\nu))\La_0.
\end{equation*}

For $M = (m_{ij}) \in \BFn$ and $M' = (m'_{ij}) \in \BCrys(\FF^{n'})$, let $M\ot M'$ be the element in $\BCrys(\FF^{n+n'})$ whose $(i,j)$-entry is given by
\[
    (M\ot M')_{ij} = \begin{cases} m_{ij} & \text{if $1\le j\le n$},\\
        m'_{i, j-n} & \text{if $n < j\le n+n'$.}
    \end{cases}
\]
This gives a canonical isomorphism $\BCrys(\FFn)\ot \BCrys(\FF^{ n'}) \cong \BCrys(\FF^{ n+n'})$ as a $U_q(\gl_\infty)$-crystal.

Put $\FFn_+ = \FF_{\R_{>0}, n}$ and $\FFn_- = \FF_{\R_{\le0} ,n}$ which are representations of $\UqglPlus[{>0}] \ot \UpglPlus[n]$ and $\UqglPlus[{\le 0}] \ot \UpglPlus[n]$, respectively. There is a $\Q(q)$-linear isomorphism
\begin{equation}\label{eq:horizontal split of FFn}
\begin{tikzcd}[row sep=tiny]
	{\FFn_+ \ot \FFn_-} & \FFn \\
	{w_{M_+} \ot w_{M_-} } & {w_{M_+ \ot M_-}}
	\arrow[from=1-1, to=1-2]
	\arrow[maps to, from=2-1, to=2-2]
\end{tikzcd},
\end{equation}
for $M_\pm\in \BCrys(\FFn_\pm)$.
It is an isomorphism of $\left( \UqglPlus[>0] \ot \UqglPlus[\le 0] \right)\ot \UpglPlus[n]$-modules,
where we regard $\FFn$ as a module over $\UqglPlus[>0] \ot \UqglPlus[\le 0]$ $\subset U_q(\lie{gl}_\infty)$.
Recall that $\ot$ in \eqref{eq:horizontal split of FFn} is $\ov{\ot}_+$ by \eqref{eqref:our choice} when it is understood with respect to the action of $\UpglPlus[n]$.

\newcommand{\Mnlm}{M^n(\mu,\nu)}
\newcommand{\Mnpluslm}{M^{n+1}(\mu,\nu)}

For $\mu\in \cP$ and $n\ge \ell(\mu)$, let $H^n_+(\mu) \in \BCrys(\FFn_+)$ be given by
\[
    H^n_+(\mu)_{ij} = \begin{cases}
        1 & i \le \mu_j ,\\
        0 & \text{otherwise}.
    \end{cases}
\]
Similarly, for $\nu \in \cP$ and $n\ge \ell(\nu)$, let $H^n_-(\nu)\in \BCrys(\FFn_-)$ be given by
\[
    H^n_-(\nu)_{ij} = \begin{cases}
        1 & i \le -\nu_{n+1-j}, \\
        0 & \text{otherwise}.
    \end{cases}
\]
For $n\ge \ell(\mu)+\ell(\nu)$, we have  $\Mnlm = H^n_+(\mu)\ot H^n_-(\nu) \in \BCrys(\FFn)$, and
$w_{\Mnlm} = w_{H^n_+(\mu)} \otimes w_{H^n_-(\nu)}$ under the isomorphism in \eqref{eq:horizontal split of FFn}.

\section{Fock space $\FFinftyM$ of infinite level}
\label{sec:Fock space of infinite level}
\subsection{Parabolic $q$-boson algebra $\UqslZero$ and $\MM$}\label{sec:uqslzero and mm}

Let $\QUE{\lie{sl}_\infty}$ be the subalgebra of $\QUE{\lie{gl}_\infty}$ generated by $e_i, f_i, q^{h_i}$ ($i\in \Z$).
We consider the parabolic $q$-boson algebra $\Uqgp$ in \Cref{subsec:B_q},
when $\lie{g} = \lie{sl}_\infty$ and $\lie{p}$ is the parabolic subalgebra associated with $J = \Z \setminus \{0\}$.
We may identify its weight and coweight lattices as $P^\vee = \bigoplus_{i\in \Z}\Z(E_{ii} - E_{i+1i+1})$, $P = \bigoplus_{i\in \Z} \Z\epsilon_i \oplus \Z \La_0$, 
$P^\vee_J = \bigoplus_{i\in \Z, i\ne 0}\Z(E_{ii} - E_{i+1i+1}) \subset P^\vee$, and $P_J =  \bigoplus_{i\in \Z} \Z\epsilon_i \subset P$.

For simplicity, we write $\UqslZero = \Uqgp$, $P_0 = P_J$ and $P_0^+ = P_J^+$.
Note that $P_0^+ = \{\, \La_{\mu,\nu} \,|\, \mu,\nu \in \cP\,\}$.
We also write $V_0(\La) = V_J(\La)$ and $(\LCrys_0(\La), \BCrys_0(\La)) = (\LCrys_J(\La), \BCrys_J(\La))$ for $\La \in P^+_0$.

Let $\MM  = V_0(0)$.
By \eqref{eq:parabolic Verma-2}, we identify $\MM$ with $\QUE{\lie{sl}_\infty}[-]/\sum_{i\in \Z\setminus \{0\} }\QUE{\lie{sl}_\infty}[-] f_i$, which induces a natural $Q_-$-grading on $\MM$, and denote the image of $1\in \QUE{\lie{sl}_\infty}[-]$ by the same symbol.

From the proof of \Cref{thm:crystal base of V_J}, we may regard $\BCrys_0(\MM) = \{\, b\in \BCrys(\infty) \,|\, \eps^\ast_i(b) = 0 \text{ for } i\ne 0\,\}$,
where $\BCrys(\infty)$ is the crystal of $\QUE{\lie{sl}_\infty}[-] = \QUE{\gl_\infty}[-]$.
By \cite{SST}, $\BCrys_0(\MM)$ can be realized as
\begin{equation}\label{eq:crystal realization of MM}
    \BCrysZero(\MM) = \{\, M \in \mathrm{Mat}_{\Z_{<0} \times \Z_{\ge 0}}(\Z_{\ge 0}) \,|\, M_{ij} = 0 \text{ for all but finitely many } i,j\,\}.
\end{equation}
Here, we regard each row (resp. column) of $M \in \BCrysZero(\MM)$ as an element of the crystal of a $\UqglPlus$-module $V(s\eps_1)$ 
where $s$ is the sum of the entries in that row (resp. $\QUE{\gl_{\le 0}}$-module $V(-t\eps_0)$ where $t$ is the sum of entries in that column), and then apply the tensor product rule row by row (resp. column by column) from top to bottom (resp. left to right). Moreover $\tf_0M$ is given by increasing the $(0,1)$-entry of $M$ by $1$.

Note that \[
 \wt(\MM)  = \left\{ \, \sum_{i\in \Z} c_i\eps_i \, \middle| \, c_i \ge 0 \  (i>0), \ c_j \le 0 \  (j \le 0), \ \sum_{i\in \Z} c_i = 0 \, \right\}.
 \]
 This gives the following alternative description of the partial order on $\cP^2$.

\begin{lem}\label{lem:alternative description of the partial order on P2}
    For $(\mu,\nu),(\zeta,\eta)\in \cP^2$,
    $(\mu,\nu) \ge (\zeta,\eta)$ if and only if $\Lambda_{\mu,\nu} - \Lambda_{\zeta,\eta} \in \wt(\MM)$.
\end{lem}
\begin{proof}
We may write $\La_{\mu,\nu} - \La_{\zeta,\eta} = \sum_{i\in \Z} c_i \eps_i + c\La_0$ for some $c_i\in \Z$ and $c\in \Z$ (cf. \eqref{eq:def of fundamental weight}).
The condition $|\mu| - |\nu| = |\zeta| - |\eta|$ is equivalent to $c = 0$. 
The condition $c_i \ge 0$ for $i > 0$ (resp. $c_j\le 0$ for $j \le 0$) is equivalent to $\mu \supset \zeta$ (resp. $\nu \supset \eta$).
\end{proof}

\newcommand{\UqglInfMinusZero}{\UqglPlus \ot \QUE{\gl_{\le 0}}}
We may also regard $\MM$ as a module over $\UqglInfMinusZero \subset \UqslZero$. Here, $\UqglInfMinusZero = U_q(\lie{l})$ in \Cref{subsec:integrable representations of parabolic q-boson algebras} with $J = \Z\setminus \{0\}$. By \eqref{eq:crystal realization of MM}, we see that $\BCrysZero(\MM)$ is a disjoint union of the connected components of $\textrm{diag}(\la_1,\la_2,\cdots)$ for $\la\in \cP$ as a $\UqglInfMinusZero$-crystal.
This implies $\MM = \bigoplus_{\la\in \cP} \MM(\la)$, where $\MM(\la)$ is the irreducible $\UqglInfMinusZero$-submodule of $\MM$ corresponding to $\la$, that is, $\MM(\la) = V_{\lie{l}}(\sum_{i\ge 1} \lambda_i(\eps_i - \eps_{-i+1}))$.

Recall that the map $r:= {}_0r^-$ in \Cref{subsec:parabolic q-derivations} induces a homomorphism of $\UqglInfMinusZero$-modules
\[
    \MM \longrightarrow V_{\lie{l}}(-\alpha_0)\ot \MM.
\]
By \Cref{lem:ri-annihilator is 1}, it is injective on $\bigoplus_{\la \ne (0)} \MM(\la)$. Hence it gives an embedding of $\MM(\la)$ into $V_{\lie{l}}(-\alpha_0)\ot \MM$ for $\la \in \cP\setminus \{(0)\}$.

This yields an embedding of $\UqglInfMinusZero$-crystals.
More precisely, since the decomposition of $\MM$ is multiplicity-free, it follows from \cite[Theorem 3]{Kas91} that
\begin{equation}\label{eq:rzerominus induces crystal map up to scaling}
r^{-1}\left( \LCrys(V_{\lie{l}}(-\alpha_0))\ot \LCrysZero(\MM) \right) = \LCrys(\MM((0))) \oplus \bigoplus_{\la\in \cP\setminus \{\emptyset\}} q^{n_\la} \LCrys(\MM(\la)),
\end{equation}
for some $n_\la\in \Z$, where $\LCrys(\MM(\mu)) = \LCrysZero(\MM) \cap \MM(\mu)$ for $\mu\in \cP$.
So we have a morphism of (abstract) $\UqglInfMinusZero$-crystals
\begin{equation} \label{eq:crystal morphism induced by r-}
    \overline{r}: \BCrysZero(\MM) \setminus \{1\} \longrightarrow \BCrys(V_{\lie{l}}(-\alpha_0)) \ot \BCrysZero(\MM),
\end{equation}
since $\BCrysZero(\MM)$ is a disjoint union of crystals of $\MM(\la)$.

\subsection{Fock space $\FFinftyM$ of infinite level}\label{subsec:limit of Fock space}

We introduce a Fock space of infinite level, which is the main object in this paper.

First, note that $\FF\otimes \MM$ has a structure of $\UqslZero$-module by \Cref{prop:comodule}.
Recall that $\FF$ and $\MM$ are integrable $\QUE{\lie{gl}_{>0}}\ot \QUE{\gl_{\le 0}}$-module by construction, and hence so is $\FF\ot\MM$ by \Cref{prop:comodule}. This implies that $\FF\ot \MM$ belongs to $\mc{O}^{\rm int}_{\UqslZero}$. It has a singular vector $|0\rangle\otimes 1$ of weight $0$, which generates a submodule isomorphic to $\MM$ by \Cref{thm:semisimple}. Hence, there exists an injective $\UqslZero$-linear map
\[
    \phi: \MM \longrightarrow \FF\otimes \MM.
\]
By tensoring with $\FFn$ for $n\ge 1$, this naturally extends to an injective $\UqslZero$-linear map
\[
    \phi_{n,n+1}: \FFnM \longrightarrow \FFnplusM.
\]
Then we define $\FFinftyM$ as the direct limit of $\{\FFnM \,|\,n\ge 0\}$ along with $\phi_{n,n+1}$'s:
\[
    \FFinftyM = \varinjlim_{n} \FFnM.
\]
Let \[ \phi_n:\FFnM \longrightarrow \FFinftyM \] be the canonical embedding.
Note that $\FFnM$ carries a structure of a $\QUE{\lie{gl}_n}[][p]$-module, acting on $\FFn$ with $p = -q^{-1}$, which commutes with the action of $\UqslZero$.
By definition of $\phi_{n,n+1}$, it is clear that $\phi_n$ is $\QUE{\lie{gl}_n}[][p]$-linear.
Hence, we may naturally define a $\UpglPlus$-module structure on $\FFinftyM$ as follows:
For $u\in \UpglPlus$ and $x\in \FFinftyM$, take a sufficiently large $n\in \N$ such that $u\in \QUE{\gl_n}[][p]$ and $x = \phi_n(x_n)$ for some $x_n \in \FFnM$, and then define $ux = \phi_n(ux_n)$.

As a $\UqslZero$-module, $\FFinftyM$ has a weight space decomposition with weights in $P_0 = P/\Z\La_0$.
We may assume that $\FFnM$ has a weight space decomposition with weights in $P$ via the identification of $P_0$ as a subgroup of $P$ (see \Cref{sec:uqslzero and mm}), which is given by $\eps_i + \Z\La_0 \to \eps_i$.
This weight space decomposition of $\FFnM$ in terms of $P$ is compatible with $\phi_n$, and hence induces the one on $\FFinftyM$. We denote by $\wt$ and $\dot{\wt}$ the weights for the action of $\UqslZero$ and $\UpglPlus$ on $\FFinftyM$, respectively.

\subsection{Crystal base of $\FFinftyM$}
\label{subsec:crystal base of Fock spaces}

Let us describe a crystal base of $\FFinftyM$.
For $n\ge 1$, let \[ (\LCrysZero(\FFnM), \BCrysZero(\FFnM)) = (\LCrys(\FFn) \ot \LCrysZero(\MM), \BCrys(\FFn) \ot \BCrysZero(\MM)), \]
which is a crystal base of $\FFnM$ as a $\UqslZero$-module by \Cref{thm:tensor product rule}.
It is also a crystal base as a $\QUE{\lie{gl}_n}[][p]$-module for $n\ge 2$.
Since $\phi(1) = |0\rangle \ot 1$ and $\BCrysZero(\MM)$ is connected, we have $\phi(\LCrysZero(\MM))\subset \LCrysZero(\FF\ot\MM)$ and $\ov\phi(\BCrysZero(\MM))\subset \BCrysZero(\FF\ot\MM)$, where $\ov{\phi}$ is the induced map at $q = 0$. This implies that for $n\ge 1$
\begin{equation*}
    \phi_{n,n+1}(\LCrysZero(\FFnM)) \subset \LCrysZero(\FFnplusM), \quad \ov{\phi}_{n,n+1}(\BCrysZero(\FFnM)) \subset \BCrysZero(\FFnplusM),
\end{equation*}
where $\ov{\phi}_{n,n+1}$ commutes with the crystal operators $\te_i,\tf_i$ ($i\in \Z$) for $\UqslZero$ and $\dot{\te}_j,\dot{\tf}_j$ ($j = 1, \cdots , n-1$) for $\UpglPlus[n]$.
Define
\[
    \LCrysZero(\FFinftyM) = \varinjlim_n \LCrysZero(\FFnM), \quad \BCrysZero(\FFinftyM) = \varinjlim_n \BCrysZero(\FFnM).
\]
Then $(\LCrysZero(\FFinftyM), \BCrysZero(\FFinftyM))$ is a crystal base of $\FFinftyM$ as a $\UqslZero\ot \UpglPlus$-module.
Also, we have the following maps induced from $\phi_n:\FFnM \to \FFinftyM$.
\[\begin{tikzcd}[row sep=tiny]
	{\LCrysZero(\FFnM)} & {\LCrysZero(\FFinftyM)}, \\
	{\BCrysZero(\FFnM)} & {\BCrysZero(\FFinftyM)}.
	\arrow["{\phi_n}", from=1-1, to=1-2]
	\arrow["{\ov{\phi}_n}", from=2-1, to=2-2]
\end{tikzcd}\]

Let us describe $\BCrysZero(\FFinftyM)$ more explicitly.
First, let $(\FFinftyM)_0$ be the $\UqslZero \ot \UpglPlus$-submodule of $\FFinftyM$ generated by
\begin{equation}\label{eq:extremal weight vectors of zeroth filtration}
w^n_{\mu,\nu} := \phi_n (w_{M^n(\mu,\nu)} \ot 1),
\end{equation}
for $(\mu,\nu)\in \cP^2$ and $n \ge \ell(\mu) + \ell(\nu)$.

\begin{prop}\label{prop:decomposition of zeroth filration}
As a $\UqslZero\ot \UpglPlus$-module, we have
\[
    (\FFinftyM)_0 \cong \bigoplus_{(\mu,\nu)\in \cP^2} V_0(\La_{\mu,\nu})\ot V_{\mu,\nu}.
\]
\end{prop}
\pf
Suppose that $(\mu,\nu)\in \cP^2$ is given.
For $n\ge \ell(\mu)+\ell(\nu)$, let $v = \phi_n(w_{M^n(\mu,\nu)}\ot 1)$.
Since $\FFinftyM \in \mc{O}^{\rm int}_{\UqslZero}$ and $v$ is a singular vector of weight $\La_{\mu,\nu}$, $v$ generates $V_0(\La_{\mu,\nu})$ as a $\UqslZero$-module.
On the other hand, $v$ generates $V_{\mu,\nu}$ as a $\UpglPlus$-module by \Cref{lem:embedding of ext wt module}.
Hence $v$ generates $V_0(\La_{\mu,\nu})\ot V_{\mu,\nu}$.
For $n'\ge n$, let $v' = \phi_{n'}(w_{M^{n'}(\mu,\nu)} \ot 1)$.
Then $v$ is an extremal weight vector with respect to the action of $\UpglPlus$ in the $\UqslZero \ot \UpglPlus$-submodule generated by $v'$.
Hence $v$ and $v'$ generate the same irreducible component.
This implies that $(\FFinftyM)_0$ is semisimple together with the required decomposition.
\qed

\medskip

Let $\LCrysZero((\FFinftyM)_0)$ be the $A_0$-span of $\td{x}_{i_1}\cdots \td{x}_{i_k} \dot{\td{y}}_{j_1}\cdots \dot{\td{y}}_{j_l} w^n_{\mu,\nu}$'s for $(\mu,\nu)\in \cP^2$, $n\ge \ell(\mu) + \ell(\nu)$, and $i_1,\cdots, i_k\in \Z$, $j_1,\cdots, j_l\in \Z_{>0}$ ($k,l\ge 0$) with $x,y\in \{e,f\}$, and let $\BCrysZero((\FFinftyM)_0)$ be the union of connected components of $w^n_{\mu,\nu}$'s in $\BCrysZero(\FFinftyM)$.
Then $(\LCrysZero((\FFinftyM)_0), \BCrysZero((\FFinftyM)_0))$  is a crystal base of $(\FFinftyM)_0$ by \Cref{prop:decomposition of zeroth filration}.

\begin{prop}\label{thm:identifying crystal basis of FFinftyM with that of its socle}
    We have $\BCrysZero(\FFinftyM) = \BCrysZero((\FFinftyM)_0)$.
\end{prop}
\pf
Let $b\in \BCrysZero(\FFinftyM)$ be given.
We have $b \equiv \ov{\phi}_l(M_1\ot b_1) \pmod{q\LCrysZero(\FFinftyM)}$ for some $M_1\ot b_1 \in \BCrysZero(\FF^{ l}\ot \MM)$. 
Suppose that $b_1 \ne 1$. Then $b_1 = \tf_{i_1}\cdots \tf_{i_{k_1}} 1$ for some $k_1\ge 1$ and $i_1,\cdots, i_{k_1}\in \Z$.
By \Cref{thm:tensor product rule}, $\ov{\phi}(b_1) = M_2\ot b_2$ for some $M_2\in \BCrys(\FF)$ with $M_2 \ne |0\rangle$ and $b_2 \in \BCrysZero(\MM)$.
Since $M_2 \ne |0\rangle$, we have $b_2 = \tf_{j_1}\cdots \tf_{j_{k_2}} 1$ for some $k_2 < k_1$ and $j_1\cdots, j_{k_2} \in \Z$.
By induction on $k_1$, we conclude that if $n$ is large enough, then
\begin{equation}\label{eq:iterated phi on BM}
    \ov{\phi}_{n-1,n} \cdots \ov{\phi}_{0,1}(b) = M\ot 1,
\end{equation}
for some $M\in \BFn$.
Therefore, we have by \eqref{eq:iterated phi on BM}
\[
    b = \ov{\phi}_l(M_1\ot b_1) = \ov{\phi}_{l+n}(M_1\ot \ov{\phi}_{n-1,n} \circ\cdots \circ \ov{\phi}_{0,1}(b_1)) = \ov{\phi}_{l+n}((M_1\ot M)\ot 1).
\]
Since $M_1\ot M$ is connected to $w_{M^{l+n}(\mu,\nu)}$ for some $\mu,\nu$ under $\td{e}_i$ and $\dot{\td{e}}_j$'s by \eqref{eq:level n Fock space crystal decomposition}, $b$ belongs to a connected component of $w^{l+n}_{\mu,\nu}$. This implies that $b \in \BCrysZero((\FFinftyM)_0)$.
\qed
\medskip

Let 
\begin{equation}\label{eq:FRR def}
    F_{\R, \R_{>0}} = \left\{ M = (M_{ab})_{a\in \Z, b\in \Z_{>0}} \middle|\,
    \begin{aligned}
        & (M_{ab})_{a\in \Z} \in \BCrys(\FF) \text{ for all } b\\
        & (M_{ab})_{a\in \Z} = |0\rangle \text{ for } b\gg 0
    \end{aligned}
    \right\}.
\end{equation}
Since we may regard the submatrix $(M_{ab})_{a\in \Z, b\in[1,n]}$ for $M = (M_{ab}) \in F_{\R,\R_{>0}}$, as an element of $F_{\R,n}$, 
it follows from the proof of \Cref{thm:identifying crystal basis of FFinftyM with that of its socle} that
the crystal $\BCrysZero(\FFinftyM)=\BCrysZero((\FFinftyM)_0)$ can be identified with $F_{\R, \R_{>0}}$.

In \Cref{sec:Crystal valuation on FFinftyM}, we explicitly construct a socle filtration of $\FFinftyM$ and show that $\soc (\FFinftyM) = (\FFinftyM)_0$.

\subsection{A bar-invariant basis of $\FFinftyM$ with upper-triangularity}

Let us start with a $\Q$-linear involution $-$ on $\FFinftyM$ compatible with the action of $\UqslZero\ot\UpglPlus$.
We first construct an involution on each $\FFnM$ that is compatible under $\phi_{n,n+1}$. Note that $\MM=V_0(0)$ has a canonical bar-involution that fixes $1$, and $\FFn$ has a canonical bar-involution by \Cref{prop:bar-involution on Fock space}. We define
\begin{equation}\label{eq:bar-involution on FFnM}
    \ov{w\otimes m} = \Theta'_-(\ov{w} \otimes \ov{m})  \quad  (w\in \FFn, m\in \MM),
\end{equation}
where $\Theta'_-$ is the opposite quasi-$R$-matrix for $U_q(\lie{sl}_\infty)$ in \eqref{eq:R-matrix for lower comultiplications}.
It is well-defined since $\Theta'_-$ is a summation of elements of $\QUE{\lie{sl}_\infty}[+]\otimes \QUE{\lie{sl}_\infty}[-]$,
and $\QUE{\lie{sl}_\infty}[-] = \UqslZero[-]$.

\begin{lem}\label{lem:permanance of bar-involution}
    For $w\in \FFn$ and $m\in \MM$, we have $\phi_{n,n+1}(\ov{w\ot m}) = \ov{\phi_{n,n+1}(w\ot m)}$.
\end{lem}
\pf
    Write $\phi(m) = \sum_{i = 1}^N w_i\otimes m_i$ for some $w_i\in \FF$ and $m_i\in \MM$. Since $\phi$ commutes with the action of $\UqslZero[-]$ and the bar-involutions,
    \begin{align*}
        \phi_{n,n+1}(\ov{w\otimes m}) &= \phi_{n,n+1}(\Theta'_-(\ov{w}\otimes \ov{m})) = (\Theta'_-)^{(3)} (\ov{w} \otimes \sum_{i=1}^N\ov{w_i}\otimes \ov{m_i}).
    \end{align*}
    On the other hand, we have
    \begin{align*}
        \ov{\phi_{n,n+1}(w\otimes m)} &= (\Delta_- \otimes 1)(\Theta'_-)\sum_{i=1}^N \ov{w\otimes w_i} \otimes m_i \\
        &= (\Theta'_-)^{(3)} (\dot{\Theta}_{\ov{\otimes}_+} \otimes 1) \sum_{i=1}^N \ov{w}\otimes \ov{w_i} \otimes \ov{m_i}.
    \end{align*}
Under the canonical identification $\FFn\otimes \FF \to \FFnplus$, $w_i$ is sent to an element that is a linear combination of products of
$w_{(i,n+1)}$ for $i\in \Z$. So the action of $\QUE{\lie{gl}_{n+1}}[+][p]$ on $w_i$ vanishes.
Since $\dot{\Theta}_{\ov{\otimes}_+} = \dot{\Pi}^{-1} \dot{\Theta}^{\op} \dot{\Pi}^{-1}$ and $\dot{\Theta}^{\op}$ is a summation of elements of $\QUE{\lie{gl}_{n+1}}[-][p] \otimes \QUE{\lie{gl}_{n+1}}[+][p]$, we see that $\dot{\Theta}_{\ov{\otimes}_+} (\ov{w}\otimes \ov{w_i}) = \ov{w}\otimes \ov{w_i}$.
\qed

\begin{prop}
There exists a $\Q$-linear involution $-$ on $\FFinftyM$ such that $\ov{u\cdot x} = \ov{u}\cdot \ov{x}$,  $\ov{\dot{u}\cdot x} = \ov{\dot{u}}\cdot \ov{x}$ for
$u\in \UqslZero$, $\dot{u}\in \UqglPlus$, and $x\in \FFinftyM$.
\end{prop}
\pf
Let $x\in \FFinftyM$ be given. We have $x = \phi_n(w\ot m)$ for some $n$ and $w\ot m\in \FFnM$, and then define $\ov{x} = \phi_n(\ov{w\ot m})$.
By \Cref{lem:permanance of bar-involution}, this gives a well-defined $\Q$-linear involution satisfying the properties.
\qed

\begin{rem}
    Our convention for the comultiplication on $\UpglPlus$ in \eqref{eqref:our choice} is crucial for \Cref{lem:permanance of bar-involution} to hold.
\end{rem}

Now we construct a canonical basis of $\FFinftyM$ by the same arguments as \cite[Section 24.2]{Lu93}.
For $M\in \BCrys(\FFn)$, define a sequence of non-negative integers $c(M): \Z \to \Z_{\ge 0}$ by $c(M)(i) = \sum_{k=1}^n M_{i+k, k}$. For $M,N\in\BCrys(\FFn)$, we define $M\ge N$ if $\wt M = \wt N$, $\dot{\wt} M = \dot{\wt} N$, and $c(M) \ge c(N)$ with respect to the lexicographic order starting from $-\infty$. It is a partial order on $\BCrys(\FFn)$ that is locally finite, that is, for $M_1$ and $M_2\in \BCrys(\FFn)$, there exists finitely many $N$ such that $M_1< N < M_2$.

By \eqref{eq:relations for q-wedge} and \eqref{eq:bar-involution on type A exterior algebras}, $\ov{w_M}$ is a linear combination of $w_M$ and $w_N$'s for $N\ne M$, where $N$ is obtained from $M$ by replacing two 1's at diagonal position to anti-diagonal position. Moreover, the coefficient of $w_M$ is 1. Thus,
\[
   \ov{w_M} \in w_M + \sum_{N < M} A w_N.
\]
By \cite[Section 2.4.2]{Lu93}, there exists a unique $G(M)\in \LCrys(\FFn)$ such that $\ov{G(M)} = G(M)$ and $G(M) \equiv M \pmod{q\LCrys(\FFn)}$.

For $b\in \BCrysZero((\FFinftyM)_0)$, let $n(b)$ be the minimal $n$ such that $b\in \ov{\phi}_n(\BCrysZero(\FFnM))$.
Let $b_i \in \BCrysZero((\FF\ot \MM)_0)$ be given such that $\wt(b_1) = \wt(b_2)$ and $\dot{\wt}(b_1) = \dot{\wt}(b_2)$.
Define $b_1 \ge b_2$ if $n(b_1) > n(b_2)$ or
$n(b_1) = n(b_2)$ with $\wt(b'_1) > \wt(b'_2)$ or
$n(b_1) = n(b_2)$ with $\wt(b'_1) = \wt(b'_2)$ and $M'_1 \ge M'_2$,
where $b_i = M_i \ot b'_i\in \BCrysZero(\FFnM)$ with $n = n(b_i)$ for $i = 1,2$.
It is a partial order on $\BCrysZero((\FFinftyM)_0)$, which can be easily seen to be locally finite.

Let $\{G_0(b) \,|\, b\in \BCrysZero(\MM)\}$ be the canonical basis of $\MM$.
By \eqref{eq:bar-involution on FFnM}, we have
\[
    \begin{aligned}
        \ov{w_M \ot G_0(b)} &= \Theta'_-(\ov{w_M} \ot G_0(b)) =  \ov{\Theta^{\op}} (\ov{w_M} \ot G_0(b)) \\
        & \in \ov{w_M} \ot G_0(b) + \left( \QUE{\lie{sl}_\infty}[+]_A \ot \QUE{\lie{sl}_\infty}[-]_A \right) (\ov{w_M} \ot G_0(b)).
    \end{aligned}
\]
Then, it is easy to see that $\ov{w_M \ot G_0(b)}$ is an $A$-linear combination of $w_N\ot G_0(b')$'s with $M\ot b \ge N\ot b'$.

Therefore, by the standard argument together with \Cref{lem:permanance of bar-involution}, we see that there exists a unique $\Q(q)$-basis  
$\{ \, G(b) \, | \, b\in \BCrysZero(\FFinftyM) \, \}$ of $\FFinftyM$ satisfying
\begin{enumerate}
    \item $\ov{G(b)}=G(b)$,
    \item $G(b) \in w_{M} \ot G_0(b_0)  + \sum_{(M'\ot b'_0 )\le b} q\Z[q] w_{M'} \ot G_0(b'_0)$ where $b = \ov{\phi}_n (M\ot b_0) \in \BCrysZero(\FFnM)$.
\end{enumerate}

\subsection{$\UqslZero$-isotypic decomposition of $\FFinftyM$}

Since $\FFinftyM$ is an integrable $\UqslZero$-module, it admits a $\UqslZero$-isotypic decomposition.

\begin{thm}\label{thm:isotypic decomposition of F_inftyM}
    As a $\UqslZero\ot \UpglPlus$-module, we have
    \begin{align*}
        \FFinftyM & \cong \bigoplus_{(\mu,\nu)\in \cP^2} V_0(\Lambda_{\mu,\nu})\otimes \left( V_{\mu,\emptyset} \otimes V_{\emptyset, \nu} \right).
    \end{align*}
\end{thm}

\pf
For $\la\in \cP$ and $n\ge \ell(\la)$, let $L^n_-(\la)\in \BCrys(\FFn_-)$ be given by
\[
    L^n_-(\lambda)_{ij} = \begin{cases}
        1 & i \le -\lambda_j, \\
        0 & -\lambda_j < i.
    \end{cases}
\]
By convention, we assume $\lambda_i = 0$ for $i > \ell(\lambda)$.

For $(\mu,\nu)\in \cP^2$, we write $E^n(\mu,\nu) = H^n_+(\mu) \ot L^n_-(\nu) \in \BCrys(\FFn)$.
We have $w_{E^n(\mu,\nu)} = w_{H^n_+(\mu)} \otimes w_{L^n_-(\nu)}$ under the isomorphism $\FFn\cong \FFn_+ \ot \FFn_-$ in \eqref{eq:horizontal split of FFn}.
Note that $\phi_n(w_{E^n(\mu,\nu)} \ot 1) \in \FFinftyM$ and $\ov{\phi}_n(E^n(\mu,\nu)\ot 1)\in \BCrysZero(\FFinftyM)$ are independent of $n$,
which we denote by $w_{E(\mu,\nu)}$ and $E(\mu,\nu)$, respectively.

We claim that $E := \{ \, w_{E(\mu,\nu)} \, | \, \mu,\nu\in \cP \, \}$ generates $\FFinftyM$ as a $\UqslZero \otimes \UpglPlus$-module.
Let $S$ be the submodule of $\FFinftyM$ generated by $E$.

First, we claim that $E^n := \{\, E^n(\mu,\nu) \,|\, \mu,\nu\in \cP, \ell(\mu) + \ell(\nu) \le n\,\}$ generates $\FFn$ as a $(\UqglInfMinusZero) \ot \UpglPlus[n]$-module.
Note that
\begin{equation*}
\begin{aligned}
\FFn_+ \cong \bigoplus_{\substack{ \mu \in \cP \\ \ell(\mu) \le n}} V(\gl_{>0}, \La_{\mu,\emptyset}) \ot V^n_{\mu,\emptyset}, \\
\FFn_- \cong \bigoplus_{\substack{ \nu \in \cP \\ \ell(\nu) \le n}} V(\gl_{\le 0}, \La_{\emptyset,\nu}) \ot V^n_{\emptyset,\nu}.
\end{aligned}
\end{equation*}
where $V(\gl_{\ast}, \La)$ ($\ast =\  >0, \le 0$) denotes the irreducible highest weight $U_q(\gl_\ast)$-module with highest weight $\La$.

Note that $w_{H^n_+(\mu)}$ is a highest weight vector which generates $V(\gl_{>0}, \La_{\mu,\emptyset}) \ot V^n_{\mu,\emptyset}$.
Similarly, $w_{L^n_{-}(\nu)}$ generates $V(\gl_{\le 0}, \La_{\emptyset,\nu}) \ot V^n_{\emptyset,\nu}$.
On the other hand, since $w_{E^n(\mu,\nu)} = w_{H^n_+(\mu)} \ot w_{L^n_-(\nu)}$, a tensor product of highest weight vector and a lowest weight vector, it generates $V^n_{\mu,\emptyset} \ot V^n_{\emptyset,\nu}$ as a $\UpglPlus[n]$-module.
Therefore, $E^n$ generates $\FFn_+ \ot \FFn_-$.
This implies that $\phi_n(\FFn\ot 1)\subset S$ for all $n$.

Next, we claim that $\phi_n(\FFn\ot \MM) \subset S$ for all $n$, which implies that $\FFinftyM = S$.
For $\gamma\in \wt(\MM)$, let $\MM_{\ge \gamma} = \bigoplus_{\beta \ge \gamma} \MM_\beta$. We use induction on the height of $\gamma \in \wt(\MM)$ to show that $\phi_n(\FFn\ot \MM_{\ge \gamma}) \subset S$.
When $\gamma = 0$, it follows from the previous argument.
Let $\gamma <0$ and $m\in \MM_\gamma$ be given.
We have $m = \sum_{i\in \Z} f_im_i$ for some $m_i\in \MM_{\gamma + \alpha_i}$, which are non-zero for only finitely many $i$'s.
For $w\in \FFn$, consider $w\ot m\in \FFn\ot \MM_{\ge \gamma}$.
Since $f_i(w\ot m_i) = f_iw \ot m_i + c_iw\ot f_im_i$ for some $c_i \in \Q(q)$, $S$ contains $f_i(w\ot m_i)$, and $f_iw\ot m_i$ by induction hypothesis, we conclude that $w\ot f_im_i\in S$ for $c_i\ne 0$. 
This completes the induction.

\newcommand{\Filt}[1]{(\FFinftyM)(#1)}
Now we prove the decomposition of $\FFinftyM$.
For $d\in \Z_{>0}$, let $\Filt{d}$ be the submodule of $\FFinftyM$ generated by $w_{E(\mu,\nu)}$ for $(\mu,\nu)\in \cP^2$ with $|\mu|,|\nu|\le d$.
We use induction on $d$ to show that
\begin{equation}\label{eq:main statement of the proof of the isotypic decomposition}
    \Filt{d} = \bigoplus_{|\mu|,|\nu|\le d} V_0(\Lambda_{\mu,\nu}) \otimes \left( V_{\mu,\emptyset} \ot V_{\emptyset,\nu}\right).
\end{equation}

If $d = 0$, then \eqref{eq:main statement of the proof of the isotypic decomposition} clearly holds.
Suppose that \eqref{eq:main statement of the proof of the isotypic decomposition} is true for $d-1$.
Let $(\zeta,\eta)\in \cP^2$ be given with $|\zeta| = d$ or $|\eta|=d$.
We claim that $\ov{w_0}$, the image of $w_0 = w_{E(\zeta,\eta)}$ in ${\FFinftyM}/{\Filt{d-1}}$ generates a submodule isomorphic to
$V_0(\La_{\zeta,\eta})\ot \left(  V_{\zeta,\emptyset}\ot V_{\emptyset,\eta} \right)$.
Since $V_0(\La_{\zeta,\eta})$ is irreducible, it suffices to show that $w_0$ generates $V_0(\La_{\zeta,\eta})$ as a $\UqslZero$-module and $V_{\zeta,\emptyset} \ot V_{\emptyset,\eta}$ as a $\UpglPlus$-module, respectively.
Then \eqref{eq:main statement of the proof of the isotypic decomposition} follows from the semisimplicity of $\FFinftyM \in \mc{O}^{\rm int}_{\UqslZero}$ and the induction hypothesis.

Let us prove the claim. Suppose first that $\ov{w_0}$ is non-zero in ${\FFinftyM}/{\Filt{d-1}}$.
In this case, $w_0$ is a singular vector of weight $\La_{\zeta,\eta}$ in ${\FFinftyM}/{\Filt{d-1}}$.
    Indeed, $\ov{w_0}$ is singular with respect to the $\UqglInfMinusZero$-action.  Moreover, $e_0w_0$ is a linear combination of $w_M \ot 1$'s for $M \in \BFn$, where $M = M_+ \ot M_- \in \BCrys(\FFn_+)\ot \BCrys(\FFn_-)$ such that $M_+$ contains 1's as many as $|\zeta| - 1$ and $M_-$ contains 0's as many as $|\eta| - 1$.
    Therefore, $e_0w_0$ is contained in $\Filt{d-1}$.

So $\ov{w_0}$ generates $V_0(\La_{\zeta,\eta})$ as a $\UqslZero$-module.
On the other hand, recall that $w_0 = \phi_n(w_{E^n(\zeta,\eta)} \ot 1)$ for $n\ge \ell(\zeta) + \ell(\eta)$.
Since $w_{E^n(\mu,\nu)} $ generates $V^n_{\mu,\emptyset} \ot V^n_{\emptyset,\nu}$ and $\phi_{n,n+1}(w_{E^n(\mu,\nu)} \ot 1) = w_{E^{n+1}(\mu,\nu)} \ot 1$, we see that $\ov{w_0}$ generates $V_{\zeta,\emptyset} \ot V_{\emptyset,\eta} $.
Next, suppose that $w_0 \in \Filt{d-1}$. Then by the previous argument, $\FFinftyM$ does not have an isotypic component for $V_0(\La_{\zeta,\eta})$,
which is a contradiction, since $w_{M^{\ell(\zeta) + \ell(\eta)}(\zeta,\eta)} \ot 1$ is a singular vector of weight $\La_{\zeta,\eta}$.
Hence, the claim is proved.

Finally, we obtain a decomposition of $\FFinftyM$ by taking a limit of $\Filt{d}$ as $d\to \infty$.
\qed

\begin{rem}
    Recall that $V_{\mu,\emptyset} \otimes V_{\emptyset, \nu} $ in \Cref{thm:isotypic decomposition of F_inftyM} is $\ov{\ot}_+$ by \eqref{eqref:our choice},
    while we may replace it with $\ov{\ot}_-$ by \Cref{rem:$R$ matrix}.
    The isomorphism in \Cref{thm:isotypic decomposition of F_inftyM} is compatible with the bar-involutions on both sides. Here, the bar-involution on the left hand side is the one in Section 6.4,
    and that on the right hand side is the tensor product of the bar-involutions on each factor. Note that $ V_{\mu,\emptyset} \ot V_{\emptyset, \nu}$ is generated by $ u_{{\mu,\emptyset}} \ot u_{{\emptyset,\nu}}$ as a $\UpglPlus$-module, and has a unique bar-involution fixing $ u_{{\mu,\emptyset}} \ot u_{{\emptyset,\nu}}$.
\end{rem}

\subsection{Filtrations of $\FFinftyM$}\label{subsec: filtration of F}

In this subsection, we introduce a filtration of $\FFnM$ in terms of its isotypic component, which will be used to construct a filtration of $\FFinftyM$.

Let $(\mu,\nu)\in \cP^2$ be given. For $\gamma\in \wt(\MM)$ and $n\ge \ell(\mu)+\ell(\nu)$, let
\begin{equation}
\label{eq:def of Hn gamma}
H^n_\gamma(\mu,\nu) = \{ \, b\in \BCrysZero(\MM)_\gamma \,|\, \Mnlm\ot b\in \HCrysZero{\FFnM} \, \},
\end{equation}
and
$H^n(\mu,\nu) = \bigcup_{\gamma \in \wt(\MM)} H^n_\gamma(\mu,\nu)$,
where $\BFnM$ is understood as a $\UqslZero\ot \UpglPlus[n]$-crystal. 
By \Cref{rem:characterizing highest weight elements}, $b\in H^n_\gamma(\mu,\nu)$ if and only if $\varepsilon_i(b) \le \langle h_i, \wt \Mnlm \rangle$ for all $i\in \Z$. This implies that $H^n_\gamma(\mu,\nu) \subset H^{n+1}_\gamma(\mu,\nu)$, and furthermore it stabilizes as $n\to \infty$, which we denote by $H_\gamma(\mu,\nu)$.

We take a set $\{x_b \,|\, b\in \BCrysZero(\MM) \} \subset \LCrysZero(\MM)$ such that $x_b \equiv b \pmod{q\LCrysZero(\MM)}$.
For $n\ge 0$ and $\gamma \in \wt(\MM)$, we define
\begin{equation} \label{eq:def of FFnM gamma}
       (\FFnM)_{\ge \gamma} =  \UqslZero[-] \ot \QUE{\gl_n}[-][p] \text{-span of } \left\{\, w_{\Mnlm} \ot x_b \, | \,  b\in H^n_{\delta}(\mu,\nu), \delta\ge \gamma \,\right\},
\end{equation}

Note that $\FFnM$ is a semisimple $\UqslZero\ot \UpglPlus[n]$-module.
By \Cref{thm:tensor product rule} and \eqref{thm:level n Fock space decomposition}, each irreducible component of $\FFnM$ is $V_0(\La_{\zeta,\eta})\ot V(\doteps^n_{\mu,\nu})$ for some $(\mu,\nu), (\zeta,\eta) \in \cP^2$ with $\ell(\mu) + \ell(\nu) \le n$.
We denote the corresponding isotypic component by $(\FFnM)^{(\zeta,\eta)}_{(\mu,\nu)}$,
and the canonical projection by $\pi^{(\zeta,\eta)}_{(\mu,\nu)}: \FFnM \to (\FFnM)^{(\zeta,\eta)}_{(\mu,\nu)}$.
Recall the partial order $\le$ on $\cP^2$ in \eqref{eq:weaker partial order on cP^2}.

\begin{lem}\label{lem:imposing relation on partitions in a filtration}
Under the above hypothesis, we have the following.
\begin{enumerate}
\item $(\FFnM)^{(\zeta,\eta)}_{(\mu,\nu)}$ is non-zero only when $(\mu,\nu)\le (\zeta,\eta)$.
\item As a $\UqslZero\ot \UpglPlus[n]$-module,
\[
    (\FFnM)_{\ge \gamma} = \bigoplus_{\substack{(\mu,\nu)\le (\zeta,\eta) \\ \La_{\zeta,\eta} - \La_{\mu,\nu} \ge \gamma}} (\FFnM)^{(\zeta,\eta)}_{(\mu,\nu)}.
\]
\end{enumerate}
\end{lem}
\pf
(1) Let $b\in \HCrys{\FFnM}$ be given whose connected component is isomorphic to $\BCrysZero(\La_{\zeta,\eta}) \ot \BCrys(\doteps^n_{\mu,\nu})$.
Since $\BCrysZero(\FFnM) \cong \BCrys(\FFn) \ot \BCrysZero(\MM) \cong \bigsqcup_{(\mu,\nu)} (\BCrys(\La^n_{\mu,\nu}) \ot \BCrysZero(\MM) ) \ot \BCrys(\doteps^n_{\mu,\nu})$,
we have by \Cref{thm:tensor product rule} and \Cref{rem:characterizing highest weight elements} that $b = u_{\La^n_{\mu,\nu}} \ot b'\ot u_{\doteps^n_{\mu,\nu}}$ for some $b'\in \BCrysZero(\MM)$.
This implies that $\La_{\zeta,\eta} - \La_{\mu,\nu} \in \wt(\MM)$, or equivalently, $(\mu,\nu) \le (\zeta,\eta)$ by \Cref{lem:alternative description of the partial order on P2}.

(2) It follows directly from (1) and an analogue of \Cref{lem:highest weight crystal representatives generate filtration} for $\UqslZero\ot \UpglPlus[n]$-modules. Note that the decomposition is independent of the choice of $x_b$ ($b\in \BCrysZero(\MM)$).
\qed

\medskip
For $(\mu,\nu),(\zeta,\eta)\in \cP^2$ with $(\mu,\nu) \le (\zeta,\eta)$ and $\ell(\mu) + \ell(\nu) \le n$, we let
\begin{equation}\label{eq:def of FFnM filtrations}
\begin{aligned}
    (\FFnM)_{(\mu,\nu)}^{\le(\zeta,\eta)} &= \bigoplus_{(\sigma,\tau)\le (\zeta,\eta)} (\FFnM)_{(\mu,\nu)}^{(\sigma,\tau)}, \\
    (\FFnM)_{\ge (\mu,\nu)}^{(\zeta,\eta)} &= \bigoplus_{(\sigma,\tau) \ge (\mu,\nu)} (\FFnM)_{(\sigma,\tau)}^{(\zeta,\eta)}, \\
    (\FFnM)^{(\zeta,\eta)} &= \bigoplus_{ (\sigma,\tau)} (\FFnM)^{(\zeta,\eta)}_{(\sigma,\tau)}, \\
    (\FFnM)_{(\mu,\nu)} &= \bigoplus_{(\sigma,\tau)} (\FFnM)^{(\sigma,\tau)}_{(\mu,\nu)}.
\end{aligned}
\end{equation}
We define $(\FFnM)_{> (\mu,\nu)}^{(\zeta,\eta)}$ similarly.
For $d\in \Z_{\ge 0}$, let
\[
(\FFnM)_{\ge -d} = \sum_{\substack{\gamma \in \wt(\MM) \\ ( \La_0, \gamma )  \ge -d}} (\FFnM)_{\ge \gamma}.
\]
Then $\left\{(\FFnM)_{\ge -d} \right\}_{d\in \Z_{\ge 0}}$ gives a filtration of $\FFnM$ since we have $(\FFnM)_{\ge -d+1} \subset (\FFnM)_{\ge -d}$.
We define $(\FFnM)_{> -d}$ similarly.

\begin{prop}\label{lem:application of branching rule to filtration}\label{lem:partition filtration permanence}
    Under the above hypothesis, we have
    \[
        \phi_{n,n+1}\left( (\FFnM)^{(\zeta,\eta)}_{(\mu,\nu)} \right) \subset (\FFnplusM)^{(\zeta,\eta)}_{\ge (\mu,\nu)}.
    \]
    Therefore, we have for $\gamma\in \wt(\MM)$,
    \[
        \phi_{n,n+1} \left( (\FFnM)_{\ge \gamma} \right) \subset (\FFnplusM)_{\ge \gamma}.
    \]
\end{prop}
\pf
    Let $(\mu^i,\nu^i),(\zeta^i,\eta^i)\in \cP^2$ ($i = 1,2$) be given such that $(\zeta^i,\eta^i) \ge (\mu^i,\nu^i)$.
    By composing $\phi_{n,n+1}$ with the associated inclusion and projection, we have a map
    \[
        (\FFnM)^{(\zeta^1,\eta^1)}_{(\mu^1,\nu^1)} \longrightarrow (\FFnplusM)^{(\zeta^2,\eta^2)}_{(\mu^2,\nu^2)}.
    \]
    Suppose that it is non-zero. Since it is $\UqslZero$-linear, we have $\zeta := \zeta^1=\zeta^2$ and $\eta:=\eta^1=\eta^2$.
    Also, since it is $\UpglPlus[n]$-linear, we have from the branching rule for the pair $(\lie{gl}_{n+1},\lie{gl}_n)$ (\cite[Theorem 8.1.1]{GW}) that
    \begin{gather*}
        (\mu^2)_1 \ge (\mu^1)_1 \ge (\mu^2)_2 \ge (\mu^1)_2 \ge \cdots (\mu^2)_{\ell(\mu^2)} \ge (\mu^1)_{\ell(\mu^2)}, \\
        (\nu^2)_1 \ge (\nu^1)_1 \ge (\nu^2)_2 \ge (\nu^1)_2 \ge \cdots (\nu^2)_{\ell(\nu^2)} \ge (\nu^1)_{\ell(\nu^2)}.
    \end{gather*}
    In particular, this implies $(\mu^1,\nu^1)\le (\mu^2,\nu^2)$, which proves the first statement.

    The second statement follows from the first one
    by Lemmas \ref{lem:alternative description of the partial order on P2} and \ref{lem:imposing relation on partitions in a filtration}.
\qed
\medskip

Now we define for $\gamma \in \wt(\MM)$
\[
    (\FFinftyM)_{\ge \gamma} = \varinjlim_{n} \, (\FFnM)_{\ge \gamma}
\]
to be the submodule of $\FFinftyM$ induced from $(\FFnM)_{\ge \gamma}$ with respect to $\phi_{n,n+1}$ for $n\ge 0$, which is well-defined by \Cref{lem:partition filtration permanence}.
For $d \in \Z_{\ge 0}$, we define
\begin{equation}\label{def of FFinftyM d}
    (\FFinftyM)_{\ge -d} = \sum_{\substack{\gamma\in \wt(\MM) \\ ( \La_0,\gamma ) \ge -d}} (\FFinftyM)_{\ge \gamma}.
\end{equation}
Then $\left\{ (\FFinftyM)_{\ge -d} \right\}_{d\in \Z_{\ge 0}}$ gives a filtration of $\FFinftyM$.
We also define $(\FFinftyM)_{> \gamma}$ and $(\FFinftyM)_{> -d}$ similarly.

\section{Crystal valuations}
\label{sec:crystal valuation}
\newcommand{\AlgA}{U}

\subsection{Crystal valuation}
\label{subsec:crystal valuation}

We introduce a weaker version of crystal lattice $L$, where $L$ is not necessarily a free $A_0$-module. We describe such an $A_0$-submodule in terms of valuations. It is used to study a structure of a non-semsimple module including its socle filtration in \Cref{sec:Crystal valuation on FFinftyM}.

Let $\val: \Q(q)\to \Z\sqcup \{\infty\}$ be a canonical valuation on $\Q(q)$ given by $\val(a) = \max\{n \,|\, q^{-n}a\in A_0\}$. A valuation on a $\Q(q)$-space $V$ is a function $\val:V \to \Z \sqcup \{\infty\}$ satisfying
\begin{enumerate}
\item $\val(v) = \infty$ if and only if $v=0$,
\item $\val(cv)=\val(c)+\val(v)$ for all $c\in \Q(q)$ and $v\in V$,
\item $\val(v+w)\ge \min\{\val(v),\val(w)\}$ for all $v,w\in V$.
\end{enumerate}
Then $\val$ is naturally in 1-1 correspondence with an $A_0$-submodule $L\subset V$ such that $L$ spans $V$ as a $\Q(q)$-space, and no nonzero element of $L$ is infinitely divisible by $q$, where the correspondence is given by $L_\val = \{v\in V \,|\, \val(v) \ge 0\}$ and $\val_L(m) = \max\{ n \,|\, q^{-n}m\in L\}$.
By abuse of notation, we often identify $\val$ with $L_\val$ when there is no confusion.
Note that $L$ need not be free, while it is true when $L$ is finitely generated as an $A_0$-module (cf. \cite[Theorem 2.4.1]{KT18}).

In this section, $\texttt{U}$ denotes either $U_q(\g)$ or $U_q(\g,\mf{p})$ in \Cref{sec:parabolic q-boson algebra} for a symmetrizable Kac-Moody algebra $\g$ and its parabolic subalgebra $\lie{p}$. 

\begin{df}
Let $V$ be an integrable $\texttt{U}$-module.
A \emph{crystal valuation} on $V$ is a valuation $\val:V \to \Z \sqcup \{\infty\}$ satisfying
\begin{enumerate}
\item  $\val(v) = \min\{ \val(v_\mu) \,|\, \mu\in P\}$ for $v = \sum v_\mu \in V$ with $v_\mu \in V_\mu$,
\item $\val(\te_i v) \ge \val(v)$ and $\val(\tf_i v) \ge \val(v)$ for all $v\in V$ and $i\in I$.
\end{enumerate}
Equivalently, $L_\val = \bigoplus_{\mu} (L_\val)_{\mu}$, and $\te_i(L_\val)\subset L_\val$, $\tf_i(L_\val)\subset L_\val$ for $i\in I$, respectively.
\end{df}

In particular, $L_{\val}$ is a crystal lattice if $L_{\val}$ is a free $A_0$-module.

\begin{rem}
For a crystal valuation $\val$, $L_\val$ is not necessarily $A_0$-free. Let $f\in \Q\llbracket q\rrbracket \setminus A_0$ (for example, $f = \sqrt{1+q}$) and $L = \{\, (r,s)\in \Q(q)^{\oplus 2} \, | \, rf + s \in \Q\llbracket q\rrbracket \,\}$. Then $L$ is an $A_0$-submodule of $\Q(q)^{\oplus 2}$, which is not free and has no nonzero element infinitely divisible by $q$
    (see \cite[Lemma 5.2.7]{KT18}).
    If we regard $\Q(q)L = \Q(q)^{\oplus 2}$ as a two copy of trivial $U_q(\lie{sl}_2)$-representation, then $\val_L$ is a crystal valuation which does not give a crystal lattice.
\end{rem}

\subsection{Decomposition into isotypic components}

Let $\texttt{V}(\la)$ denote an integrable highest weight $\texttt{U}$-module $V(\la)$ or $V_J(\la)$ for $\la$ in $P^+$ or $P_J^+$ with a crystal base $(\texttt{L}(\la),\texttt{B}(\la))$.
Let $\mc{O}^{\rm int}_{\texttt{U}}$ be a category of integrable $\texttt{U}$-modules such that each object is isomorphic to a direct sum of $\texttt{V}(\la)$  whose set of weights is finitely dominated.

Suppose that $V\in\mc{O}^{\rm int}_{\texttt{U}}$ is given and $\val$ is a crystal valuation on $V$. 
\begin{lem}[{cf. \cite[Lemma 5.2.4 (1)]{HK}}]\label{lem:crystal valuation analogue of ei vanishing}
    Suppose that {\em $V = \texttt{V}(\la)^{\oplus m}$} for some $m\ge 1$. Let $v\in V$ be a weight vector that is not of highest weight. If $\val(\te_i v) \ge 0$ for all $i\in I$, then $\val(v) \ge 0$.
\end{lem}
\pf
We have an isomorphism of $\Q(q)$-spaces $V \cong V_\la \ot_{\Q(q)} \texttt{V}(\la)$, which sends $u\in V_\la$ to $u \ot u_\la$.
Under this isomorphism, we may write $v = \sum_k v_k\ot w_k$ for some $v_k\in \mathbb{Q}(q)$ and $w_k\in \texttt{L}(\la)$ such that $w_k \pmod{q\texttt{L}(\la)}$ are pairwise distinct elements of $\texttt{B}(\la)$. Fix $k$. By applying $\te_i$'s to $w_k$ to obtain $u_\la \pmod{q\texttt{L}(\la)}$, we see that $v_k \in L_\val$. This implies $\val(v) \ge 0$.
\qed

\begin{lem}[{cf. \cite[Lemma 2.6.3]{Kas91}}] \label{lem:uniqueness of crystal valuation}
    Suppose that $V=  V_1\oplus V_2$, where $V_1$ and $ V_2$ are {\em $\texttt{U}$}-submodules such that {\em $V_1 = \texttt{V}(\la)^{\oplus m}$} for a maximal weight $\la$ of $V$ and $m\ge 1$.
    Then we have $\val(v_1+v_2) = \min\{\val(v_1),\val(v_2)\}$ for $v_i\in V_i$ $(i=1,2)$.
\end{lem}
\pf It can be proved in the same way as \cite[Lemma 5.2.1]{HK}, where \cite[Lemma 5.2.4 (1)]{HK} is replaced by \Cref{lem:crystal valuation analogue of ei vanishing}. Note that \cite[Lemma 5.2.1]{HK} assumes the existence of crystal base, but this is not necessary in the proof of the decomposition of crystal lattices.
\qed
\medskip

The crystal valuation is compatible with the isotypic decomposition of $V$ as follows.

\begin{prop}\label{cor:uniqueness of crystal valuation}
    Suppose that {\em $V = \bigoplus_{s\in S} \texttt{V}(\la_s)^{\oplus m_s}$} for some $S$ and $m_s \in \Z_{\ge 0}$ $(s\in S)$ such that $\la_s\ne \la_{s'}$ for $s\ne s'$.
    \begin{enumerate}
    \item We have $\val(v) = \min\{\,\val(v_s) \,|\, s\in S\,\}$ for {\em $v = \sum v_s$ with $v_s\in \texttt{V}(\la_s)^{\oplus m_s}$}, equivalently, {\em $L_\val = \bigoplus_{s\in S} L_\val \cap (\texttt{V}(\la_s)^{\oplus m_s})$}.
    \item If $L_\val$ is $A_0$-free, then there exists a unique {\em $\texttt{U}$}-linear automorphism $\phi$ of $V$ such that {\em $\phi(L_\val \cap (\texttt{V}(\la_s)^{\oplus m_s}))=\texttt{L}(\la_s)^{\oplus m_s}$}.
    \end{enumerate}
\end{prop}
\pf
(1) It follows from a repeated application of \Cref{lem:uniqueness of crystal valuation}.
(2) The proof is similar to that of \cite[Theorem 3]{Kas91}.
We have $L_\val = \bigoplus_{s\in S} L_\val \cap (\texttt{V}(\la_s)^{\oplus m_s})$ by (1). Since $L_\val$ is $A_0$-free and $\texttt{B}(\la_s)$ is connected, each $L_\val \cap (\texttt{V}(\la_s)^{\oplus m_s})$ is free and coincides with the crystal lattice generated by $L_s = L_\val \cap (\texttt{V}(\la_s)_{\la_s}^{\oplus m_s})$ under $\tf_i$'s.
For $s\in S$, let $\phi_s$ be the automorphism of $\texttt{V}(\la_s)^{\oplus m_s}$ sending the $A_0$-basis of $L_s$ to the one consisting of $u_{\la_s}$'s in each component. Then $\phi_s(L_\val \cap (\texttt{V}(\la_s)^{\oplus m_s})) = \texttt{L}(\la_s)^{\oplus m_s}$, and $\phi = \bigoplus_{s\in S} \phi_s$ is the required automorphism.
\qed
\medskip

\begin{rem}\label{rem:same holds for UqUp}
    One may also apply \Cref{cor:uniqueness of crystal valuation} when $\texttt{U}$ is an  algebra other than $U_q(\g)$ or $U_q(\g,\mf{p})$. For example, one may take $\texttt{U}=U_q(\mf{sl}_{\infty,0})\ot U_p(\gl_n)$ and $\mc{O}^{\rm int}_{\texttt{U}}$ to be the category of $\texttt{U}$-modules $V$,
    which is isomorphic to a direct sum of $V_0(\La_{\zeta,\eta})\ot V(\doteps^m_{\mu,\nu})$ for $(\mu,\nu), (\zeta,\eta) \in \cP^2$ with $\ell(\mu) + \ell(\nu)  \le n$,
    and whose set of weights is finitely dominated. This will be used in \Cref{sec:Crystal valuation on FFinftyM}.
    \end{rem}

Proposition \ref{cor:uniqueness of crystal valuation} also holds for a direct sum of extremal weight modules over $U_q(\gl_{>0})$, whose set of weights is not necessarily finitely dominated.
\begin{prop}\label{prop:uniqueness of crystal for direct sum of extremal weight modules}
    Let $V = \bigoplus_{(\mu,\nu)\in S} V_{\mu,\nu}^{\oplus m_{\mu,\nu}}$ for some finite set $S$ and $m_{\mu,\nu}\in \Z_{\ge0}$,
    and let $\val$ be a crystal valuation on $V$. Then
    \begin{enumerate}
    \item We have $\val(v) = \min\{\val(v_{\mu,\nu}) \,|\, (\mu,\nu)\in S\}$ for {\em $v = \sum v_{\mu,\nu}$ with $v_{\mu,\nu}\in V_{\mu,\nu}^{\oplus m_{\mu,\nu}}$}.
    \item If $L_\val$ is $A_0$-free, then there exists a unique $\UqglPlus$-linear automorphism $\psi$ of $V$ such that $\psi(L_\val \cap V_{\mu,\nu}^{\oplus m_{\mu,\nu}}) = \LCrys_{\mu,\nu}^{\oplus m_{\mu,\nu}}$.
    \end{enumerate}
\end{prop}
\pf
We may write $V = \bigoplus_{k=1}^l V_{\mu^k,\nu^k}^{\oplus m_k}$ such that $(\mu^i,\nu^i) \gneqq (\mu^j,\nu^j)$ implies $i < j$.
Note that $V$ is a direct limit of $V^n:=\bigoplus_{(\mu,\nu)\in S} (V_{\mu,\nu}^n)^{\oplus m_{\mu,\nu}}$ for $n\ge \max \{ \,\ell(\mu^k) + \ell(\nu^k) \,|\, 1\le k\le l\,\}$ (see the proof of \Cref{prop:irreducibility of extremal weight modules}).

Put $V_1 = V_{\mu^1,\nu^1}^{\oplus n_1}$ and $V_2 = \oplus_{k\ne 1} V_{\mu^k,\nu^k}^{\oplus m_k}$. 
Then $(V_2)_{\eps^n_{\mu^1,\nu^1}} = 0$. Otherwise, we have $(V^{n'}_{\mu^k,\nu^k})_{\eps^n_{\mu^1,\nu^1}} \ne 0$ for some $k>1$ and $n' > n$. 
Applying the action of Weyl group of $\gl_{n'}$ to $\eps^n_{\mu^1,\nu^1}$, we have $(V^{n'}_{\mu^k,\nu^k})_{\eps^{n'}_{\mu^1,\nu^1}} \ne 0$.
This implies $\eps^{n'}_{\mu^k,\nu^k} > \eps^{n'}_{\mu^1,\nu^1}$, which is a contradiction.
Now we apply the same arguments as in the proof of \Cref{lem:uniqueness of crystal valuation} and use induction on $l$ to prove (1).

By (1), it suffices to consider the case $V \cong V_{\mu,\nu}^{\oplus m}$ in (2).
By \Cref{cor:uniqueness of crystal valuation}, there exists a unique $\UqglPlus[n]$-linear automorphism $\psi^n$ of $V^n$ such that $\psi^n(L_\val\cap V^n) = L_{\mu,\nu} \cap (V^n_{\mu,\nu})^{\oplus m}$ for $n\ge \ell(\mu)+\ell(\nu)$. For $n'> n$, the restriction of $\psi^{n'}$ on $V^n$ is $\psi^n$ by the uniqueness, so combining $\psi^n$'s give a desired $\UqglPlus$-linear automorphism $\psi$ of $V$. The uniqueness of $\psi$ follows from that of $\psi^n$.
\qed

\medskip

A non-semisimple module may not have a unique crystal valuation, as illustrated in the following example.

\begin{ex}\label{ex:non-uniqueness of crystal lattices in 10011}
   Consider a $U_q(\gl_{>0})$-module $V_{(1),\emptyset} \ot_-V_{\emptyset, (1)}$, which is not semisimple (cf.\eqref{eq:example of non-split seq}).
   We may regard $V_{(1),\emptyset}=\bigoplus_{n\in\N} \Q(q)v_n$, where an $\UqglPlus$-action is given by $q^{\eps_i}v_j = q^{\delta_{ij}}v_j$, $f_iv_j = \delta_{ij}v_{i+1}$, $e_iv_j = \delta_{i+1,j}v_i$,
    and then $\LCrys_{(1),\emptyset}=\bigoplus_{n\in \N} A_0 v_n$.
    Similarly, we let $V_{\emptyset,(1)} = \bigoplus_{n\in\N} v_{n^\vee}$ with $q^{\eps_i}v_{j^\vee} = q^{-\delta_{ij}} v_{j^\vee}$, $f_iv_{j^\vee} = \delta_{i+1,j}v_{i^\vee}$, $e_iv_{j^\vee} = \delta_{ij}v_{(i+1)^\vee}$, and $\LCrys_{\emptyset,(1)} = \bigoplus_{n\in \N} A_0 v_{n^\vee}$.
   By the tensor product rule, $\LCrys := \LCrys_{(1),\emptyset} \ot \LCrys_{\emptyset, (1)}$ is a crystal lattice of $V_{(1),\emptyset}\ot V_{\emptyset, (1)}$, 
   
    The socle of $V_{(1),\emptyset} \ot V_{\emptyset, (1)}$ is isomorphic to $V_{(1),(1)}$, which is spanned by $B=\{\,v_m \ot v_{n^\vee}, (m\ne n),\ qv_{n}\ot v_{n^\vee} + v_{n+1}\ot v_{(n+1)^\vee}\, (n\ge 1)\}$. 
    The $A_0$-span $\LCrys_{(1),(1)}$ of $B$ is a crystal lattice of $V_{(1),(1)}$, and $\LCrys \cap V_{(1),(1)}=\LCrys_{(1),(1)}$.

    On the other hand, consider for $S \subset \N$   
\begin{equation*}
 \LCrys_S = \LCrys + \sum_{n\in S} A_0 D_n, \quad  \text{where $D_n = \sum_{k = 1}^n (-q)^{-n-1+k}v_k\ot v_{k^\vee} \ (n \ge 1)$.}
\end{equation*}
Then it is not difficult to see that $\LCrys_S$ is also a crystal valuation of $V_{(1),\emptyset}\ot V_{\emptyset, (1)}$.
    Indeed, we have
    \[ 
        \te_k D_n = \begin{cases}
            -\frac{1}{1+q^2}v_{n}\ot v_{(n+1)^\vee} & \text{if } k = n, \\
            0 & k\ne n,
        \end{cases}\quad
        \tf_k D_n = \begin{cases}
            -\frac{1}{1+q^2}v_{n+1}\ot v_{n^\vee} & \text{if } k = n, \\
            0 & k \ne n.
        \end{cases} 
    \]
    We also have $\LCrys_S\cap V_{(1),(1)}=\LCrys_{(1),(1)}$ for all $S$. In fact, $\LCrys_S$ is a crystal lattice, since it is easy to see that its finite rank submodule is finitely generated (cf. \cite[Theorem 2.4.3]{KT18}).
\end{ex}

\subsection{Saturated crystal valuation}

Let $V$ be a $\Q(q)$-vector space. Given valuations $\val_1,\val_2$ on $V$, we define $\val_1\ge \val_2$ when $\val_1(v)\ge \val_2(v)$ for all $v\in V$ (in particular, $L_{\val_1} \supset L_{\val_2}$). This gives a partial order on the set of valuations on $V$.
Let $W$ be a subspace of $V$. For a valuation $\mathbbm{w}$ on $W$, we say that $\val$ is an extension of $\mathbbm{w}$ if $\val|_W = \mathbbm{w}$.

\begin{lem}\label{lem:characterizing saturated crystal valuations}
    The following are equivalent.
    \begin{enumerate}
        \item $\val$ is a maximal extension of $\mathbbm{w}$ with respect to the partial order $\ge$,
        \item the $\Q$-linear map $L_\mathbbm{w} / qL_\mathbbm{w} \to L_\val / qL_\val$, which is induced from the inclusion $L_\mathbbm{w} \subset L_\val$, is an isomorphism,
        \item $L_\val / L_\mathbbm{w}$ is a divisible $A_0$-module, that is, any element of $L_\val / L_\mathbbm{w}$ is divisible by $q$.
    \end{enumerate}
\end{lem}
\pf
Put $L = L_\mathbbm{w}$ and $L' = L_\val$.

(1)$\Rightarrow$(2): 
Let $l'\in L'\setminus qL'$ be given. By the maximality of $L'$, we have $L \subsetneqq (L' + q^{-1}A_0l') \cap W$.
So there exists $l''\in L'$ and $c\in A_0$ such that 
$l'' + q^{-1}cl' \in W\setminus L$. Also, $c$ is a unit of $A_0$. If not, $q^{-1}c\in A_0$ and hence $l'' + q^{-1}cl'\in L$. 
If we put $l := qc^{-1}l'' + l' \in L' \cap W = L$, $l + qL$ is a preimage of $l' + qL'$. Since $L/qL\to L'/qL'$ is clearly injective, it is bijective.

(2)$\Rightarrow$(3): For any $l'\in L'$, there exists $l\in L$ and $l''\in L'$ such that $l' = l + ql''$. Hence (3) holds.

(3)$\Rightarrow$(1): Suppose that there exists a strictly larger extension $L''\supsetneqq L'$. Then there exists $l'\in L'$ such that $q^{-1}l' \in L''\setminus L'$. By assumption, we have $l' = ql'' + l$ for some $l\in L$ and $l''\in L'$. Then $q^{-1}l = q^{-1}l' - l'' \in L'' \cap W = L$, so $q^{-1}l' \in L'$, which is a contradiction.
\qed
\medskip

Suppose that $V,W$ are {integrable $\texttt{U}$-modules}, and $\val$ is a crystal valuation.

\begin{df}\label{df:saturated crystal lattices}\mbox{}
    \begin{enumerate}
        \item $\val$ is called \emph{{a maximal crystal valuation}} with respect to $W$ if it is maximal (with respect to $\ge$) among all crystal valuations on $V$ that extends $\val|_W$.
        \item $\val$ is called \emph{{a saturated crystal valuation}} with respect to $W$ if it is maximal among all valuations on $V$ that extends $\val|_W$.
    \end{enumerate}
\end{df}

\begin{prop}
    Suppose that $V$ is of finite length, and $W = \soc V$ with a crystal valuation $\mathbbm{w}$ on $W$. If there exists at least one crystal valuation on $V$ extending $\mathbbm{w}$, then there exists a maximal crystal valuation on $V$ extending $\mathbbm{w}$.
    In particular, there exists a maximal crystal valuation on $V$ extending $\val|_{\soc V}$.
\end{prop}
\pf
We apply Zorn's lemma. Suppose that $\{\val_s\}_{i\in S}$ are crystal valuations on $V$ extending $\mathbbm{w}$, parametrized by a totally ordered set $S$, satisfying $\val_s \ge \val_t$ when $s\ge t$. We show that $\td{\val}: V\to\Z\sqcup\{\infty\}$ defined by $\td{\val}(v) = \max\{\,\val_s(v) \,|\, s\in S\,\}$ is a crystal valuation. The conditions (2) and (3) are straightforward.

Let $v\in V$ be given. Note that $\texttt{U} v$ has a nontrivial intersection with $W$. Take $w\in (\texttt{U}v \cap W) \setminus\{0\}$. We may write $w = \sum_{\bm{i}} c_{\bm{i}} \td{x}_{i_1}\cdots \td{x}_{i_k} v$ ($x = e,f$), where the summation is over sequences $\bm{i} = (i_1,\cdots, i_k)$ of elements of $I$ with $c_{\bm{i}}\in \mathbb{Q}(q)$.
Then for any crystal valuation $\val'$ on $V$ extending $\mathbbm{w}$, we have 
$\mathbbm{w}(w) = \val'(w) \ge \min\{\,\val(c_{\bm{i}}) \,|\,\bm{i}\,\}  + \val'(v)$, and hence $\val'(v)$ is bounded above.
This verifies the first condition.
\qed
\medskip 

By definition, a saturated crystal valuation $\val$ with respect to $W$ is a maximal crystal valuation with respect to $W$.
By \Cref{lem:characterizing saturated crystal valuations}, we have the following characterization of $\val$.

\begin{prop}\label{cor:saturated iff mod-q bijective}
    The crystal valuation $\val$ is saturated with respect to $W$ if and only if
    $(L_\val\cap W)/q(L_\val\cap W) \to L_\val / qL_\val$ is an isomorphism of $\Q$-spaces.
In particular, when $(\LCrys,\BCrys)$ is a crystal base of $V$, $\val_{\LCrys}$ is saturated with respect to $W$ if and only if $(\LCrys \cap W, \BCrys)$ is a crystal base of $W$.
\end{prop}

\begin{ex}\label{ex:saturated crystal lattice in 10011}
Let $\texttt{U} = \UqglPlus$. Consider $V_{\emptyset,\nu} \ot V_{\mu,\emptyset}$ and its crystal lattice $\LCrys_{\emptyset,\mu}\ot \LCrys_{\lambda,\emptyset}$ for $\mu,\nu\in \cP$.
By \Cref{cor:identical filtration on lowest tensor highest}, there exists a submodule of $V_{\emptyset,\nu}\ot V_{\mu,\emptyset}$ isomorphic to $V_{\mu,\nu}$ which is generated by $u_{\mu,\nu}=u_{\emptyset,\nu}\ot u_{\mu,\emptyset}$.
Note that $\ms{L}_{\mu,\nu}=V_{\mu,\nu}\cap (\LCrys_{\emptyset,\mu}\ot \LCrys_{\lambda,\emptyset})$ can be identified with the $A_0$-span of $\td{x}_{i_1}\dots\td{x}_{i_l}u_{\mu,\nu}$ for $i_i,\dots,i_l$ ($l\ge 0$) with $x=e$ or $f$ for each $i_k$.
    Since $\BCrys_{\lambda,\mu} \simeq \BCrys_{\emptyset,\mu}\ot \BCrys_{\lambda,\emptyset}$ by \Cref{thm: ext wt crystal of >0}, we conclude from \Cref{cor:saturated iff mod-q bijective} that $\LCrys_{\emptyset,\mu}\ot \LCrys_{\lambda,\emptyset}$ is saturated with respect to $V_{\mu,\nu}$. In fact, $V_{\mu,\nu}$ is the socle of $V_{\emptyset,\nu}\ot V_{\mu,\emptyset}$ (see \Cref{cor:main application - crystal basis}). It can be viewed as a generalization of \eqref{eq:example of non-split seq}.
\end{ex}

\begin{ex}\label{ex:saturated crystal lattice in 10011}
    As another example when $\texttt{U} = \UqglPlus$, consider the crystal lattice  
    $\LCrys_\N$ of $V_{(1),\emptyset}\ot V_{\emptyset,(1)}$ in \Cref{ex:non-uniqueness of crystal lattices in 10011}, which is an extension of $\LCrys_{(1),(1)}$.
    We see that $\LCrys_{(1),(1)} / q\LCrys_{(1),(1)} \to \LCrys_\N / q\LCrys_\N$ is bijective since $D_n = -qD_{n+1} - v_{n+1}\ot v_{(n+1)^\vee}$ and hence $D_n + q\LCrys_\N$ has a preimage $-qv_n \ot v_{n^\vee} -v_{n+1}\ot v_{(n+1)^\vee} + q\LCrys_{(1),(1)}$. Thus $\LCrys_\N$ is saturated with respect to its socle by \Cref{cor:saturated iff mod-q bijective}.
\end{ex}

\begin{ex}
    When $\texttt{U} = \UqslZero\ot \UpglPlus$, $\LCrysZero(\FFinftyM)$ is saturated with respect to $(\FFinftyM)_0 \subset \FFinftyM$ by \Cref{thm:identifying crystal basis of FFinftyM with that of its socle}.
\end{ex}

The saturated crystal valuations on $V$ in \Cref{ex:saturated crystal lattice in 10011} are given with respect to ${\rm soc}(V)$.
Conversely, we will show in \Cref{sec:Crystal valuation on FFinftyM} that when $\val$ is saturated with respect to $W$, we have $W = \soc V$ under certain additional assumptions (cf. \Cref{thm:crystal isomorphism implies socle}).

\begin{ex} 
Let $\texttt{U} = \QUE{\g}$ be a quantum affine algebra. Following the notations in \cite{Kas02'},
let $\varpi_i = \La_i - a^\vee_i\La_{0^\vee}$ ($i\in I$) be the $i$-th fundamental weight of level $0$, and let $\la = \sum_{i\in I\setminus \{0\}} m_i\varpi_i$ with $m_i\ge 0$. Then $u_\la=\bigotimes_{i\in I\setminus \{0\}} u_{m_i\varpi_i}$ is an extremal weight vector in $\bigotimes_{i\in I\setminus \{0\}} V(m_i\varpi_i)$, and the induced map $V(\la) \to \bigotimes_{i\in I\setminus \{0\}} V(m_i\varpi_i)$ is injective such that the crystal lattice $\bigotimes_{i\in I\setminus \{0\}} \LCrys(m_i\varpi_i)$ is saturated with respect to $V(\la)$. It was conjectured in \cite{Kas02'} and proved in \cite{BN}.
\end{ex}

From now on, we simply say that $\val$ is saturated if it is saturated with respect to $\soc V$.

\section{Saturated crystal valuations and socle filtration}
\label{sec:Crystal valuation on FFinftyM}

\subsection{Valuation on isotypic components of $\FFn$}
\label{subsec:lemmas on isotypic components}
We prove some important lemmas on the projection of standard monomials of a form $w_{M(\la)\ot A}$ in $\FF^{n+1}$ to its irreducible components in \eqref{thm:level n Fock space decomposition}.

Let $\la\in \Z^n_+$ be given with $\mu,\nu\in \cP$ satisfying \eqref{eq:integer partitions to two partitions}.
Let \[\pi^n_\la = \pi^n_{\mu,\nu} : \FFn \longrightarrow V(\La_\la)\ot V(\doteps_\la)\] be the canonical projection.
Denote $\la\cup \{0\}$ an element of $\Z^{n+1}_+$ obtained from $\la$ by adding $0$ in its entries. Let $A\in \BCrys(\FF)$ be an element of charge $0$.
\begin{lem}\label{lem:projection onto isotypic components of Fn}
For $\zeta\in \Z^{n+1}_+$, we have $\pi^{n+1}_\zeta(w_{M(\la)\ot A}) = 0$ unless $\La_\zeta \le \La_{\la\cup \{0\}}$.
\end{lem}

\pf
Let $M\in \BCrys(\FFnplus)$ be given such that $\dot{\wt}(M) = \dot{\wt}(M(\la)\ot A)$.
By \eqref{eq:level n Fock space crystal decomposition}, $M$ belongs to $\BCrys(\La_\zeta)\ot \BCrys(\doteps_\zeta)$ for some $\zeta\in \Z^{n+1}_+$.
We claim that $\La_\zeta \le \La_{\lambda \cup \{0\}}$.
Note that $M(\zeta) = \dot{\te}_{j_1}\cdots \dot{\te}_{j_l}\te_{i_1}\cdots \te_{i_k} M$ for some $i_s, j_t$'s (cf. \eqref{eq:highest weight element in FFn}). Since $\te_i$'s may only move a single $1$ to one position up in each column of $M$, the charge of each column of $M$ does not change.
This implies that $\wt(M(\zeta)) = \wt(\te_{i_1}\cdots \te_{i_k}M) \le \wt(M(\lambda) \ot |0\rangle) = \wt (M(\la\cup \{0\}))$, which proves the claim.

Therefore, for each $M \in \BCrys(\FFnplus)$ with $\wt(M) = \wt(M(\lambda)\otimes A)$ and $\dot{\wt}(M) = \dot{\wt}(M(\lambda)\otimes A)$, we may choose its representative in $\LCrys(\FFnplus)$ which is contained in $V(\Lambda_\zeta)\otimes V(\doteps_{\zeta})$ with $\Lambda_\zeta \le \Lambda_{\lambda\cup \{0\}}$. Since the correponding weight space is finite dimensional, $w_{M(\lambda)\otimes A}$ can be expressed as a linear combination of such representatives. 
It follows that the weight space of $\FFnplus$ with weight $\left(\wt(M(\lambda)\otimes A),\dot{\wt}(M(\lambda)\otimes A)\right)$ is contained in $U_q^-(\gl_\infty) \ot U^-_p(\gl_n)$-span of $w_{M(\mu)}$ for $\La_\mu\le \La_{\lambda \cup \{0\}}$.
This proves the lemma.
\qed
\medskip

Let $\val_{\LCrys(\FFn)}$ be the crystal valuation associated to $\LCrys(\FFn)$ in \eqref{eq:crystal base of Fock space} (cf. \Cref{subsec:crystal valuation}).
Under the isomorphism of \eqref{thm:level n Fock space decomposition}, $\val_{\LCrys(\FFn)}$ gives a valuation on $\bigoplus V(\La_\la)\ot V(\doteps_\la)$, which we denote by the same symbol.
Then we have the following lemma, which will be crucially used in later sections.

\begin{lem}\label{lem:q-order of projection onto isotypic components of Fn}
    If $A\ne |0\rangle$, then we have
    \[
        \lim_{n\to \infty} \val_{\LCrys(\FFnplus)}(\pi^{n+1}_{\mu,\nu}(w_{M^n(\mu,\nu) \ot A})) = \infty.
    \]
\end{lem}
\pf 
Put $w = w_{M^n(\mu,\nu)\ot A}$ and $w_{\mu,\nu} = \pi^{n+1}_{\mu,\nu}(w_{M^n(\mu,\nu) \ot A})$.
Note that the irreducible components of $\FFnplus$ are mutually orthogonal with respect to the bilinear form $\langle \cdot, \cdot \rangle$ on $\FFnplus$.
Since $\LCrys(\FFnplus) = \{\, v\in \FFnplus \,|\, \langle v, \LCrys(\FFnplus) \rangle \in A_0\,\}$ (cf. \Cref{subsec:crystal base}),
we see that $\val_{\LCrys(\FFnplus)}(w_{\mu,\nu})$ is the minimum of $\val(\langle w, x\rangle)$ for all $x\in \LCrys(\La^{n+1}_{\mu,\nu}) \ot \LCrys(\doteps^{n+1}_{\mu,\nu})$ with $\wt(x) = \wt(w) = \La^n_{\mu,\nu} + \wt(A)$ and $\dot{\wt}(x) = \dot{\wt}(w) = \doteps^n_{\mu,\nu}$, where $\val$ is the valuation on $\Q(q)$ defined in \Cref{subsec:crystal valuation}.

Recall that $\LCrys(\Lambda_{\mu,\nu}^{n+1})\otimes \LCrys(\doteps^{n+1}_{\mu,\nu})$ is the $A_0$-span of $G(b)G(b')w_{M^{n+1}(\mu,\nu)}$ for $b\in \BCrys(\lie{sl}_\infty, \infty)$ and $b'\in \BCrys(\lie{gl}_{n+1},\infty)$, where $\BCrys(\lie{g},\infty)$ denotes the crystal of $\UqgMinus$.
Since $\doteps^n_{\mu,\nu}$ is an extremal weight of $V(\doteps_{\mu,\nu}^{n+1})$,
we have $(\wt(M), \dot{\wt}(M)) = (\La^{n+1}_{\mu,\nu}, \doteps^{n}_{\mu,\nu})$ if and only if $M = \Mnlm \otimes |0\rangle$ for $M\in \BCrys(\La^{n+1}_{\mu,\nu}) \ot \BCrys(\doteps^{n+1}_{\mu,\nu})$. Hence,
we have
\[
    A_0\text{-span of}\left\{\, G(b')w_{M^{n+1}(\mu,\nu)} \, | \, b' \in \BCrys(\lie{gl}_{n+1},\infty), \dot{\wt}(b') = \doteps^{n}_{\mu,\nu} - \doteps^{n+1}_{\mu,\nu}\,\right\} = A_0 w_{M^n(\mu,\nu)\ot |0\rangle}.   
\]
Therefore, it suffices to show that
\begin{equation} \label{eq:fock_space_form_order}
    \lim_{n\to \infty} \min \left\{ \, \val \left(\left\langle w, G(b) w_{M^n(\mu,\nu)\otimes |0\rangle} \right\rangle \right) \, \middle| \, b\in \BCrys(\lie{sl}_\infty, \infty)_{\wt(A)-\La_0 }
    \,\right\} =\infty.
\end{equation}
Since $\BCrys(\lie{sl}_\infty, \infty)_{ \wt(A)-\La_0 }$ is finite, the set of $b \in \BCrys(\lie{sl}_\infty, \infty)_{ \wt(A)- \La_0 }$ with
$G(b)w_{M^n(\mu,\nu)\ot |0\rangle} \neq 0$ (or equivalently, $\ov{\pi}_{\La^{n+1}_{\mu,\nu}}(b)\neq 0$)
stabilizes as $n$ tends to $\infty$, and it is independent of $n$.
So, it suffices to show $\val\left(\langle w, G(b) w_{M^n(\mu,\nu) \ot |0\rangle}\rangle \right)$ diverges as $n\to \infty$ for each $b$.

By definition of $\langle \cdot, \cdot \rangle$ on $\FFn$,
$\langle w, G(b)w_{M^n(\mu,\nu)\ot |0\rangle} \rangle$ is equal to the coefficient of $w$ in $G(b)w_{M^n(\mu,\nu)\otimes |0\rangle}$ given as a $\Q(q)$-linear combination of standard monomials.
Write $G(b) = \sum_{\bm{i}} c_{\bm{i}} f_{i_1}\cdots f_{i_k}$ ($c_{\bm{i}} \in \Q(q)$), where the sum is over all sequences $\bm{i} = (i_1,\ldots, i_k)$ such that $\sum_{l=1}^k \alpha_{i_l} = \Lambda_0 - \wt(A)$. Each action of $f_{i_l}$ on a standard monomial $w_M$ yields new standard monomials obtained from $M$ by moving a certain 1 in $M$ to the right. In order for $G(b)w_{M^n(\mu,\nu)\otimes |0\rangle}$ to have $w$ in its expression, every action of $f_{i_l}$ in $f_{i_1} \cdots f_{i_k}$ should move a 1 at the $n+1$-th row of $M$. Hence, by the formula of $\Delta(f_i)$, the coefficient of $w$ in $f_{i_1}\cdots f_{i_k} w_{M^n(\mu,\nu)\otimes |0\rangle}$ is given by
\[
    q^{\langle \sum_{l=1}^k h_{i_l} , \Lambda^n_{\mu,\nu} \rangle}.   
\]
Since $\Lambda^n_{\mu,\nu} = \Lambda_{\mu,\nu} + (n-\ell(\lambda)-\ell(\mu))\Lambda_0$ and $A$ is of charge $0$, $i_1,\ldots, i_k$ must contain at least $0$. Therefore, $\langle \sum_{l=1}^k h_{i_l} , \Lambda^n_{\mu,\nu} \rangle$ becomes arbitrarily large as $n$ tends to $\infty$. This implies \eqref{eq:fock_space_form_order}.
\qed

\subsection{Semisimple subquotients of $\FFinftyM$}
\label{subsec: semisimple subquotient}

Suppose that $(\mu,\nu)\in \cP^2$ and $\gamma\in \wt(\MM)$ are given such that $\La_{\zeta,\eta} = \La_{\mu,\nu} + \gamma$ for some $(\zeta,\eta)\in \cP^2$.
Let us begin with the following description of images of the generators of $(\FFnM)_{\ge\gamma}$ elements under $\phi_{n,n+1}$.
\begin{lem}\label{thm:filtration_permanence}
    Let $b\in H^n_\gamma(\mu,\nu)$ be given. Then we have
    \[
      \phi_{n,n+1}(w_{\Mnlm}\otimes x_b) \equiv \sum_{b'} c_n(\mu,\nu;b,b') w_{\Mnlm \ot |0\rangle} \ot x_{b'} \pmod{(\FFnplusM)_{>\gamma}}
    \]
    for some  $c_n(\mu,\nu;b,b')\in \Q(q)$, where the sum is over $b' \in H^{n+1}_{\gamma}(\mu,\nu)$ (cf. \eqref{eq:def of Hn gamma}).
\end{lem}

\pf
Note that $x_b = u\cdot 1$ for some $u\in \UqslZero[-]_\gamma$, and $\Delta(u) \equiv t_\gamma \otimes u \pmod{\QUE{\lie{sl}_\infty}[-] \otimes \QUE{\lie{sl}_\infty}[-]_{>\gamma}}$, where $t_\gamma = \prod t_i^{m_i}$ for $\gamma = \sum_i m_i\alpha_i$.
So, we may write
\[
    \phi_{0,1}(x_b) = \sum_{i=1}^N c_i w_{A_i}\otimes x_{b_i} + q^{-(\Lambda_0, \gamma )}|0\rangle \otimes x_b,
\]
for some $A_i\in \BCrys(\FF) \setminus \{ |0\rangle \}$ and $b_i\in \BCrysZero(\MM)$ ($i=1,\cdots, N$). Then 
\begin{equation}
    \label{eq:expressing phi_n}
    \phi_{n,n+1}(w_{\Mnlm}\otimes x_b) = \sum_{i=1}^N c_iw_{\Mnlm \otimes A_i} \otimes x_{b_i} + q^{-( \Lambda_0, \gamma )}w_{\Mnlm\otimes |0\rangle} \otimes x_b.
\end{equation}

We write $A_0 = |0\rangle, b_0=b$, and $c_0 = q^{-( \Lambda_0, \gamma)}$ for convenience.
Fix $i  = 0,\ldots, N$. First, by \Cref{lem:projection onto isotypic components of Fn}, we have
\[
    w_{\Mnlm\otimes A_i} \in \bigoplus_{(\sigma,\tau)} V(\Lambda^{n+1}_{\sigma,\tau}) \otimes V(\doteps^{n+1}_{\sigma,\tau}) \subset \FFnplus
\]
where the sum is over $(\sigma,\tau)\in\cP^2$ such that $\La^{n+1}_{\sigma,\tau} \le \La^{n+1}_{\mu,\nu}$.
Regarding $\BCrysZero(\FFnplusM) = \BCrys(\FFnplus) \ot \BCrysZero(\MM)$ as a $\UqslZero\ot \UpglPlus[n+1]$-crystal, we have by \Cref{lem:highest weight crystal representatives generate filtration},
\begin{equation}\label{eq:monomial description of filtration in filtration permanence}
    \begin{split}
    w_{\Mnlm\otimes A_i} \otimes x_{b_i} \in  & \, \QUE{\lie{sl}_{\infty,0}}[-] \ot \QUE{\gl_{n+1}}[-][p]\text{-span of } \\
    & \left\{ \, w_{M^{n+1}(\sigma,\tau)} \otimes x_{b'} \, | \, (\sigma,\tau)\in \cP^2, \ b'\in H^{n+1}(\sigma,\tau) \text{ with } \La^{n+1}_{\sigma,\tau} \le \La^{n+1}_{\mu,\nu} \,\right\}.
    \end{split}
\end{equation}
We see that $\wt(b') \ge \wt(b)$ since 
\begin{equation}\label{eq:weight inequality in filtration permanence}
\begin{split}
    \wt(M^{n+1}(\sigma,\tau)) + \wt(b')  & \ge \wt(\Mnlm) + \wt(A_i) + \wt(b_i) \\
    &  = \wt(\Mnlm) + \wt(|0\rangle) + \wt(b) \\
    & = \wt(M^{n+1}(\mu,\nu)) + \wt(b). 
\end{split}
\end{equation}

Consider $w_{M^{n+1}(\sigma,\tau)} \otimes x_{b'}$ in \eqref{eq:monomial description of filtration in filtration permanence}.
If $\wt(b') > \wt(b)$, then $w_{M^{n+1}(\sigma,\tau)}\ot x_{b'} \equiv 0 \pmod{ (\FFnplusM)_{> \gamma }}$.
If $\wt(b') = \wt(b)$, then 
\eqref{eq:weight inequality in filtration permanence} and the condition of \eqref{eq:monomial description of filtration in filtration permanence} imply $\La^{n+1}_{\sigma,\tau} = \La^{n+1}_{\mu,\nu}$, so
we have $M^{n+1}(\sigma,\tau) = \Mnpluslm$, and $\wt(\Mnlm) + \wt(A_i) + \wt(b_i) = \wt(\Mnpluslm) + \wt(b')$. Hence
\[
    w_{\Mnlm\otimes A_i} \otimes x_{b_i} \in (\FFnplusM)_{> \gamma} + \QUE{\gl_{n+1}}[-][p]\text{-span of } \{\, w_{M^{n+1}(\mu,\nu)} \otimes x_{b'} \,|\, b'\in H^{n+1}_\gamma(\mu,\nu)\, \}.
\]

On the other hand, $\dot{\wt}(\Mnpluslm)$ is a dominant integral weight and $\xi := \dot{\wt}(\Mnlm \ot A_i) = \dot{\wt}(\Mnlm \ot |0\rangle)$ is in the Weyl group orbit of $\dot{\wt}(\Mnpluslm)$.
This implies that 
\begin{equation}\label{eq:extremal deescription in filtration permanence}
    \left(\QUE{\gl_{n+1}}[-][p] w_{\Mnpluslm}\right) _\xi = \Q(q) \dot{S}_{i_1}\cdots \dot{S}_{i_k} w_{\Mnpluslm} = \Q(q) w_{\Mnlm \ot |0\rangle},
\end{equation}
for some $i_1,\ldots, i_k$, where $\dot{S}_i$ denotes $S_i$ in \eqref{eq:extremal operators} with respect to $\UpglPlus[n+1]$.
Hence we have
\begin{equation*}
    w_{\Mnlm \ot A_i} \ot x_{b_i} \in (\FFnplusM)_{> \gamma} +  \Q(q)\text{-span of } \{\, w_{\Mnlm \ot |0\rangle} \ot x_{b'} \,|\, b'\in H^{n+1}_\gamma(\mu,\nu)\,\}.
\end{equation*}
This completes the proof.
\qed

\begin{rem}
    \Cref{thm:filtration_permanence} provides an alternative proof of \Cref{lem:partition filtration permanence}.
\end{rem}

Define
\begin{equation} \label{eq:def of gamma mu nu}
    (\FFnM)_{\gamma, (\mu,\nu)} = (\FFnM)_{>\gamma} \oplus (\FFnM)^{(\zeta,\eta)}_{(\mu,\nu)}.
\end{equation}
Then
\begin{equation}\label{eq:identifying multiplicity of gamma mu nu}
    \frac{(\FFnM)_{\gamma, (\mu,\nu)}}{(\FFnM)_{>\gamma}} \cong \left( V_0(\Lambda_{\zeta,\eta}) \otimes V^n_{\mu,\nu} \right)^{\oplus |H^n_\gamma(\mu,\nu)|}.
\end{equation}

On the other hand, by applying \Cref{lem:highest weight crystal representatives generate filtration} to $V(\Lambda^n_{\mu,\nu})\ot V(\doteps^n_{\mu,\nu}) \ot \MM$, we have
\begin{equation}\label{eq: standard monomial description of a filtration}
    (\FFnM)_{(\mu,\nu)}^{\le (\zeta,\eta)} = \UqslZero[-] \ot \QUE{\gl_n}[-][p]\text{-span of} \left\{\, w_{\Mnlm} \ot x_b \, \middle| \,
        b\in H^n_{\delta}(\mu,\nu), \ \delta\ge \gamma \, \right\}.
\end{equation}
Therefore, we have
\begin{equation}
    \label{eq:description of subquotient in terms of standard monomials}
    \begin{split}
    & (\FFnM)_{\gamma, (\mu,\nu)} \\ & =   (\FFnM)_{>\gamma} \, + \, \UqslZero[-]\ot\QUE{\lie{gl}_n}[-][p]\text{-span of } \left\{ \, w_{\Mnlm}\otimes x_b \,\middle|\, b\in H^n_\gamma(\mu,\nu) \, \right\}.
    \end{split}
\end{equation}
Hence,  $\left\{w_{\Mnlm}\otimes x_b \pmod{(\FFnM)_{>\gamma} } \, \middle| \, b\in H^n_\gamma(\mu,\nu) \right\}$ forms a $\Q(q)$-basis of the highest weight space of the $(\FFnM)_{\gamma,(\mu,\nu)}/(\FFnM)_{>\gamma}$.

By \Cref{lem:application of branching rule to filtration}, $\phi_{n,n+1}$ induces a map
\newcommand{\phibar}{\phi^{\gamma}}
\begin{equation}
    \phibar_{n,n+1}: \frac{(\FFnM)_{\gamma, (\mu,\nu)}}{(\FFnM)_{>\gamma}} \longrightarrow \frac{(\FFnplusM)_{\gamma, (\mu,\nu)}}{(\FFnplusM)_{>\gamma}}.
\end{equation}
Let 
\begin{equation}\label{eq:def of valuation on FFnM}
\val^n = \val_{\LCrysZero(\FFnM)}
\end{equation}
be the crystal valuation on $\FFnM$ associated to $\LCrysZero(\FFnM) = \LCrys(\FFn) \ot \LCrysZero(\MM)$.

\begin{prop}\label{lem:subquotient_injectivity}
    For all sufficiently large $n$, the following statements hold;
    \begin{enumerate}
    \item The map $\phibar_{n,n+1}$ is injective.
    \item For all $x\in (\FFnM)_{\gamma, (\mu,\nu)}/(\FFnM)_{>\gamma}$, we have
    \[
        \val^{n+1} (\phibar_{n,n+1}(x)) = \val^n(x) - (\Lambda_0, \gamma),
    \]
    \end{enumerate}
    where $x$ and $\phibar_{n,n+1}(x)$ are identified with elements of $(\FFnM)^{(\zeta,\eta)}_{(\mu,\nu)}$ and $(\FFnplusM)^{(\zeta,\eta)}_{(\mu,\nu)}$ in \eqref{eq:def of gamma mu nu}, respectively.
\end{prop}

\pf By \Cref{thm:filtration_permanence}, we have for $b\in H^n_\gamma(\mu,\nu)$,
\begin{equation}\label{eq:phi_n on wMnlm}
    \phi_{n,n+1}(w_{\Mnlm}\otimes x_b) \equiv \sum_{b'\in H^{n+1}_\gamma(\mu,\nu)} c_n(\mu,\nu; b,b') w_{\Mnlm \otimes |0\rangle} \otimes x_{b'} \pmod{(\FFnplusM)_{>\gamma}}.
\end{equation}

Note that $\dot{S}_n\dot{S}_{n-1}\cdots \dot{S}_{{\ell(\mu) + 1}}$ induces an isomorphism on the $\UqslZero$-highest weight space of \eqref{eq:identifying multiplicity of gamma mu nu}, since it is spanned by extremal weight vectors for $\UpglPlus[n+1]$ by \eqref{eq:description of subquotient in terms of standard monomials}.
Also, $\{ \, w_{M^{n+1}(\mu,\nu)} \otimes x_{b'} \,|\, b' \in H^{n+1}_\gamma(\mu,\nu)\,\}$ is a linearly independent set consisting of extremal weight vectors in $(\FFnplusM)_{\gamma,(\mu,\nu)} / (\FFnplusM)_{>\gamma}$.
Since \eqref{eq:phi_n on wMnlm} is obtained by applying $\dot{S}_n\dot{S}_{n-1}\cdots \dot{S}_{\ell(\mu) + 1}$ to
$\sum_{b'} c_n(\mu,\nu; b,b')\allowbreak w_{M^{n+1}(\mu,\nu)} \otimes x_{b'}$ (cf.~\eqref{eq:extremal deescription in filtration permanence}),
it suffices to show that the matrix  \[
    C = \left( c_n(\mu,\nu; b,b') \right)_{b,b'\in H^n_\gamma(\mu,\nu)}
\]
has non-zero determinant for all sufficiently large $n$ in order to prove (1).

To prove $\det C \ne 0$, we show
\begin{equation} \label{eq:limit of the order of c_n}
    \lim_{n\to \infty} \val\left( c_n(\mu,\nu;b, b') \right)= \begin{cases} 
        -(\Lambda_0,\gamma) & \text{ if } b = b', \\
        \infty & \text{ if } b \ne b'.
    \end{cases}
\end{equation}
We use the same notations as in the proof of \Cref{thm:filtration_permanence}. 
Fix $i = 1, \ldots, N$. 
The projection $\pi^{(\zeta,\eta)}_{(\mu,\nu)}$ factors through
\[\begin{tikzcd}[column sep=large]
	\FFnplusM & {V(\La^{n+1}_{\mu,\nu}) \ot V(\doteps^{n+1}_{\mu,\nu}) \ot \MM} & { (\FFnplusM)^{(\zeta,\eta)}_{(\mu,\nu)}}.
	\arrow["{\pi^{n+1}_{\mu,\nu} \ot 1}", from=1-1, to=1-2]
	\arrow["{}", from=1-2, to=1-3]
\end{tikzcd}\]
Since $\LCrys(\FFnplusM) = \LCrys(\FFnplus)\otimes\LCrys(\MM)$, we have $\val^{n+1} \left( \pi^{n+1}_{\mu,\nu}(w_{\Mnlm\otimes A_i})\otimes x_{b_i} \right) = \val_{\LCrys(\FFnplus)}\left( \pi^{n+1}_{\mu,\nu}(w_{\Mnlm\otimes A_i}) \right)$, 
and thus
\[
    \lim_{n\to\infty} \val^{n+1} \left( \pi^{(\zeta,\eta)}_{(\mu,\nu)}( w_{\Mnlm\ot A_i} \ot x_{b_i}) \right)= \infty,
\]
by \Cref{lem:q-order of projection onto isotypic components of Fn}.
For $i=0$, the term $q^{-( \La_0, \gamma )} w_{\Mnlm \ot |0\rangle} \ot x_b$ in \eqref{eq:expressing phi_n} contributes
on the diagonal of $C$ by $q^{-( \La_0, \gamma )}$.
Applying these observations to \eqref{eq:expressing phi_n} and \eqref{eq:extremal deescription in filtration permanence}, we obtain \eqref{eq:limit of the order of c_n}.
This also implies (2).
\qed

\begin{cor} \label{cor:subquotient_injectivity of isotypic components}
    For all sufficiently large $n$, the following statements hold;
    \begin{enumerate}
        \item The map $ \pi^{(\zeta,\eta)}_{(\mu,\nu)} \circ \phi_{n,n+1}$ induces an injection
        \[
            (\FFnM)^{(\zeta,\eta)}_{(\mu,\nu)} \longrightarrow (\FFnplusM)^{(\zeta,\eta)}_{(\mu,\nu)}.
        \]
        \item For all $x\in (\FFnM)^{(\zeta,\eta)}_{(\mu,\nu)}$, we have \[ \val^{n+1}\left( \pi^{(\zeta,\eta)}_{(\mu,\nu)}(\phi_{n,n+1}(x)) \right) = \val^n(x) - ( \Lambda_0, \gamma). \]
    \end{enumerate}
\end{cor}

\medskip

Let $(\FFinftyM)^{(\zeta,\eta)}_{\ge (\mu,\nu)} = \varinjlim_n (\FFnM)^{(\zeta,\eta)}_{\ge (\mu,\nu)}$, which is well-defined by \Cref{lem:application of branching rule to filtration}, and $(\FFinftyM)^{(\zeta,\eta)}_{> (\mu,\nu)}$ defined similarly.
Then 
\[
(\FFinftyM)_{\ge \gamma} = \sum_{\substack{(\mu,\nu)\le (\zeta,\eta), \\ \La_{\zeta,\eta} - \La_{\mu,\nu} \ge \gamma}} (\FFinftyM)^{(\zeta,\eta)}_{\ge (\mu,\nu)}.
\]
We have the following.

\begin{thm}\label{thm:isotypic decomposition of subquotients}
    As a $\UqslZero\ot \UpglPlus$-module, we have
    \[
        \frac{(\FFinftyM)_{\ge -d}}{(\FFinftyM)_{> -d }} \cong
        \bigoplus_{
            \substack{
                (\mu,\nu)\le (\zeta,\eta), \\
                |\zeta| - |\mu| = |\eta| - |\nu| = d
            }
        } \left(V_0(\Lambda_{\zeta,\eta}) \ot V_{\mu,\nu} \right)^{\oplus |H_\gamma(\mu,\nu)|},
    \]
    for $d\in \Z_{\ge 0}$
    where $\gamma = \Lambda_{\zeta,\eta} - \Lambda_{\mu,\nu} \in \wt(\MM)$ for each $(\mu,\nu) \le (\zeta,\eta)$.
\end{thm}

\pf
Let $d\in \Z_{\ge 0}$ and $(\zeta,\eta), (\mu,\nu)\in \cP^2$ be given such that $(\zeta,\eta) \ge (\mu,\nu)$ and $|\zeta| - |\mu| = |\eta| - |\nu| = d$.
For each $n\ge \ell(\mu)+\ell(\nu)$, we have an exact sequence
\begin{equation}\label{eq:ffnm exact sequence}
    0 \longrightarrow (\FFnM)^{(\zeta,\eta)}_{> (\mu,\nu)} \longrightarrow (\FFnM)^{(\zeta,\eta)}_{\ge (\mu,\nu)} \longrightarrow (\FFnM)^{(\zeta,\eta)}_{(\mu,\nu)} \longrightarrow 0.
\end{equation}
Here we regard $\left\{ (\FFnM)^{(\zeta,\eta)}_{\ge (\mu,\nu)} \right\}_{n \ge \ell(\mu)+\ell(\nu)}$ as a directed system 
whose associated maps can be identified with $\UqslZero\ot\UpglPlus[n]$-linear maps
\begin{equation} \label{eq:direct limit map of zeta eta mu nu}
    V_0(\Lambda_{\zeta,\eta})\otimes V(\doteps^n_{\mu,\nu}) \otimes \Q(q)^{\oplus H^n_\gamma(\mu,\nu)} \to
    V_0(\Lambda_{\zeta,\eta})\otimes V(\doteps^{n+1}_{\mu,\nu}) \otimes \Q(q)^{\oplus H^{n+1}_\gamma(\mu,\nu)},
\end{equation}
by \eqref{eq:identifying multiplicity of gamma mu nu}.
Taking a directed limit, we obtain
\[
    \frac{(\FFinftyM)^{(\zeta,\eta)}_{\ge (\mu,\nu)}}{(\FFinftyM)^{(\zeta,\eta)}_{> (\mu,\nu)}}
    \cong \varinjlim_{n} \, (\FFnM)^{(\zeta,\eta)}_{(\mu,\nu)}.
\]
By the branching rule for the pair $(\lie{gl}_{n+1},\lie{gl}_n)$ (\cite[Theorem 8.1.1]{GW}), any $\UpglPlus[n]$-linear inclusion $V(\doteps^n_{\mu,\nu}) \to V(\doteps^{n+1}_{\mu,\nu})$ is a scalar multiple of the canonical inclusion. Therefore, there exists a $\Q(q)$-linear map
$\iota_n: \Q(q)^{\oplus H^n_\gamma(\mu,\nu)} \to \Q(q)^{\oplus H^{n+1}_\gamma(\mu,\nu)}$ such that
the map \eqref{eq:direct limit map of zeta eta mu nu} is given by $u \ot v \ot w \mapsto u \ot v \ot \iota_n(w)$.

By \Cref{lem:subquotient_injectivity}, $\iota_n$ is injective for all sufficiently large $n$, and hence
\[
  \varinjlim_n \Q(q)^{\oplus |H^n_\gamma(\mu,\nu)|} \cong \Q(q)^{\oplus |H_\gamma(\mu,\nu)| },
\]
where the directed limit is taken with respect to $\iota_n$.
Thus, we have
\begin{equation} \label{eq:subquotients in terms of partitions}
    \frac{(\FFinftyM)^{(\zeta,\eta)}_{\ge (\mu,\nu)}}{(\FFinftyM)^{(\zeta,\eta)}_{> (\mu,\nu)}} \cong  V_0(\Lambda_{\zeta,\eta}) \ot V_{\mu,\nu} \ot \Q(q)^{\oplus |H_\gamma(\mu,\nu)|}.
\end{equation}

Next, observe that if $( \Lambda_0,\gamma) = -d$, then
\[
    (\FFinftyM)^{(\zeta,\eta)}_{\ge (\mu,\nu)}\cap (\FFinftyM)_{\ge -d+1} = (\FFinftyM)^{(\zeta,\eta)}_{> (\mu,\nu)}.
\]
Indeed, if $(\FFnM)^{(\zeta,\eta)}_{(\sigma,\tau)} \subset (\FFinftyM)_{\ge -d+1}$ for some $(\zeta,\eta) \ge (\sigma,\tau) \ge (\mu,\nu)$, then we have $(\sigma,\tau) \gneq (\mu,\nu)$.
Conversely, if $(\sigma,\tau) \gneq (\mu,\nu)$, then $\La_{\sigma,\tau} - \La_{\mu,\nu} \in \wt(\MM) \setminus \{0\}$, so $( \La_0, \La_{\sigma,\tau} - \La_{\mu,\nu} ) < 0$.

Thus,
\begin{eqnarray*}
\frac{(\FFinftyM)_{\ge -d}}{(\FFinftyM)_{> -d}} & = & \frac{(\FFinftyM)_{> -d } + \sum_{\substack{
    (\mu,\nu)\le (\zeta,\eta), \\
    |\zeta| - |\mu| = |\eta| - |\nu| = d
}} (\FFinftyM)^{(\zeta,\eta)}_{\ge (\mu,\nu)}}{(\FFinftyM)_{>-d}} \\ & \cong &  \bigoplus_{\substack{
    (\mu,\nu)\le (\zeta,\eta), \\
    |\zeta| - |\mu| = |\eta| - |\nu| = d
}} \frac{(\FFinftyM)^{\zeta,\eta}_{\ge (\mu,\nu)}}{(\FFinftyM)^{\zeta,\eta}_{> (\mu,\nu)}}.
\end{eqnarray*}
Combining with \eqref{eq:subquotients in terms of partitions}, the assertion follows.
\qed

\begin{rem}
    When $d = 0$, we have
    \[
        (\FFinftyM)_{\ge 0} = \bigoplus_{(\mu,\nu)\in \cP^2} V_0(\La_{\mu,\nu}) \ot V_{\mu,\nu},
    \]
    which is equal to $(\FFinftyM)_0$ generated by $w^n_{\mu,\nu}$ \eqref{eq:extremal weight vectors of zeroth filtration} (cf. \Cref{prop:decomposition of zeroth filration}).
\end{rem}

\subsection{A crystal valuation on socle quotients of $\FFinftyM$}
\label{subsec: a crystal valuation on socle quotients}
In this subsection, we introduce a crystal valuation on $(\FFinftyM)/(\FFinftyM)_{> -d}$.
We remark that the image of $\LCrysZero(\FFinftyM)$ under the projection $\FFinftyM\longrightarrow (\FFinftyM)/(\FFinftyM)_{> -d}$ may have an element which is infinitely divisible by $q$ (see \Cref{lem:subquotient_injectivity}), and hence it does not give a well-defined crystal valuation. 

Recall that $\val^n$ is the crystal valuation on $\FFnM$ \eqref{eq:def of valuation on FFnM}.
For $d\in \Z_{\ge 0}$, let $(\FFnM)_{-d}$ be the direct summand of $\FFnM$ isomorphic to $(\FFnM)_{\ge -d} / (\FFnM)_{>-d}$, that is,
\[
(\FFnM)_{-d} = \bigoplus_{\substack{(\mu,\nu)\le (\zeta,\eta) \\ |\zeta| - |\mu| = |\eta| - |\nu| = d}} (\FFnM)^{(\zeta,\eta)}_{(\mu,\nu)}.
\]
Let $\pi^n_{-d}:\FFnM \to (\FFnM)_{-d}$ be the canonical projection.
Since $(\FFnM)_{-d}$ exhausts isotypic components, that is, $(\FFnM)_{-d}$ and $\FFnM/(\FFnM)_{-d}$ do not have a common isotypic component, $(\FFnM)_{-d}$ has a canonical crystal lattice given as the restriction of $\LCrysZero(\FFnM) = \LCrys(\FFn)\ot \LCrysZero(\MM)$ by \Cref{cor:uniqueness of crystal valuation}.
Let $\val^n_{-d}$ be the associated crystal valuation, that is,
\[
    \val^n_{-d} (x) := \val^n \left( \pi^n_{-d}(x) \right),
\]
for $x\in \FFnM$.

The goal of this subsection is to prove the following.

\begin{thm}\label{thm:valuation limit exists}
    Let $x\in \FFinftyM$ be given such that $\phi_N(x_N) = x$ for some $x_N\in \FF^N\ot \MM$.
    For $d\in \Z_{\ge 0}$, the limit
    \[
        \val^\infty_{-d} (x) := \lim_{n\to \infty} \big( \val^n_{-d}(x_n) - dn \big)
    \]
    exists and lies in $\Z \sqcup \{\infty\}$, where $x_n = \phi_{N,n}(x_N)$ for $n > N$.
    Moreover, $\val^\infty_{-d}(x)$ is finite if and only if $x\not\in (\FFinftyM)_{> -d}$.
\end{thm}

In order to prove the theorem,
we first need a series of lemmas leading to an injectivity of $\phi_{N,n}$ followed by the projection to $(\FF^{\bullet} \ot M)_{-d+1}$ with `one degree up' (\Cref{lem:derivation_lowest_degree_coefficient}).

For $d\in \Z_{\ge 0}$, let $\MM_{-d} = \bigoplus_{\gamma\in \wt(\MM), (\gamma, \Lambda_0 ) = -d} \MM_\gamma$, and let $\MM_{\ge -d}, \MM_{> -d}$ and so on be defined in a similar way.

\begin{lem}\label{lem:the approximation lemma}
    Let $d\in \Z_{\ge 0}$, $\rho\in P_0$, and $M\in \Z_{\ge 0}$ be given. 
    Then there exists $N=N(d,\rho,M )$ such that the following statement holds for all $n \ge N$;
    For $x\in \FFn\ot (\MM_{\ge -d}) \cap \LCrysZero(\FFnM)_{n\La_0+ \rho}$,
    there exists $x^\circ\in (\FFnM)_{\ge -d} \cap \LCrysZero(\FFnM)_{n\La_0 + \rho}$ such that
    \begin{enumerate}
    \item $x\equiv x^\circ \pmod{q\LCrysZero(\FFnM)}$,
    \item $\val^n(x^\circ_{< -d}) >M$, where $x^\circ_{< -d}$ denotes the projection of $x^\circ$ in $\FFn\ot (\MM_{< -d})$.
    \end{enumerate}
\end{lem}

\pf
We may assume $x \in \BCrysZero(\FFnM) \pmod{q\LCrysZero(\FFnM)}$, say $M\ot b$. Suppose that $M\ot b$ is connected to $\Mnlm \ot b'$ for some $(\mu,\nu)\in \cP^2$ and $b'\in H^n(\mu,\nu)$.
Let $(\zeta,\eta)\in \cP^2$ be such that $\La^n_{\zeta,\eta} = \La^n_{\mu,\nu} + \wt(b')$.

Since
\begin{equation}\label{eq:proof of approximation lemma first weight bound}
    \La^n_{\mu,\nu} \ge \La^n_{\mu,\nu} + \wt(b') \ge \wt(M\ot b) = n\La_0 + \rho,
\end{equation}
(note that $0 \ge w$ for all $w\in \wt(\MM)$) and hence there are only finitely many $(\mu,\nu)\in \cP^2$ satisfying \eqref{eq:proof of approximation lemma first weight bound},
it suffices to prove the case when $(\mu,\nu)$ is fixed.

Let $v = v_{\Mnlm \ot b'} \in \LCrysZero(\FFnM)$ be the projection of $w_{\Mnlm} \ot x_{b'}$ onto the
isotypic component whose highest weight is given by the weight of $\Mnlm \ot b'$.
Then it is a highest weight vector.
Note that $M \ot b = \dot{\tf}_{j_1} \cdots \dot{\tf}_{j_l} \tf_{i_1} \cdots \tf_{i_k} ( \Mnlm\ot b')$
for some $i_1,\ldots, i_k$ and $j_1, \ldots, j_l$.
Let  $b''=  \tf_{i_1} \cdots \tf_{i_k} 1 \in \BCrys(\infty)$, where $\BCrys(\infty)$ is the crystal of $\QUE{\gl_\infty}[-]$.
Note that $\wt(M\ot b) = n\La_0 + \rho$ is fixed, and
\begin{equation}\label{eq:proof of approximation lemma second weight bound}
    \wt(b'') = \wt(M\ot b) - \wt(\Mnlm \ot b') \ge n\La_0 + \rho- \La^n_{\mu,\nu},
\end{equation}
where the rightmost term is independent of $n$.
So, there are only finitely many $b''$'s satisfying \eqref{eq:proof of approximation lemma second weight bound},
and it suffices to prove the case when $b'' \in \BCrys(\infty)$ is fixed.

Put
\[
    x^\circ = \dot{\tf}_{j_1} \cdots \dot{\tf}_{j_l} G(b'') v \in \LCrysZero(\FFnM).
\]
We claim that there exists $N' \gg 0$ such that
$x^\circ$ satisfies the required properties for all $n \ge N'$ and for all $M\ot b$ subject to the assumptions on $(\mu,\nu)$, $b'$, and $b''$.

By construction, we have $x^\circ\equiv M\ot b \equiv x \pmod{q\LCrysZero(\FFnM)}$. Also we have $v\in (\FFnM)_{\ge -d}$ since $\wt b' \ge \wt b$ and $( \wt b,\Lambda_0 ) \ge -d$. This implies $x^\circ\in (\FFnM)_{\ge -d}$.

So it remains to prove that $\val^n({x^\circ_{< -d}}) > M$, where $x^\circ_{< -d}$ denotes the projection of $x^\circ$ in $\FFn\ot (\MM_{< -d})$.
Moreover, we may assume that $x^\circ = G(b'')v$, since $\FFn \ot (\MM_{< -d})$ is stable under $\dot{\tf}_j$ and $\val^n(\dot{\tf}_j x) \ge \val^n(x)$.

\medskip
\textit{Step 1.}
First, we may write 
\begin{equation}\label{eq:proof of approximation lemma first expression}
    \Delta(G(b'')) = \sum_i x_i t_{\rho_i} \ot y_i,
\end{equation}
for some $x_i \in \QUE{\lie{gl}_{\infty}}[-]$, $y_i \in  \QUE{\lie{gl}_{\infty}}[-]_{-\rho_i}$
and $t_{\rho_i} = \prod_k t_k^{(\rho_i, \La_k)}$.
Similarly, we write
\begin{equation}\label{eq:proof of approximation lemma second expression}
    v = \sum_{j} c_j w_{M_j} \ot x_{b_j},
\end{equation}
for some $c_j\in \Q(q)$, $M_j\in \BFn$ and $b_j\in \BCrysZero(\MM)$.
Here, we have $c_j\in A_0$ since $v\in \LCrysZero(\FFnM)$.
By \eqref{eq:proof of approximation lemma first expression} and \eqref{eq:proof of approximation lemma second expression}, we have
\begin{equation}\label{eq:proof of approximation lemma third expression}
    x^\circ = G(b'')v = \sum_{i,j} c_j (x_it_{\rho_i} w_{M_j}) \ot (y_ix_{b_j}).
\end{equation}

\medskip
\textit{Step 2.}
We have  $\wt(b_j) \ge \wt(b') \ge n\La_0 + \rho - \La^n_{\mu,\nu}$ for all $j$, since $v \in (\FFnM)^{(\zeta,\eta)}_{(\mu,\nu)} \subset (\FFnM)^{\le(\zeta,\eta)}_{(\mu,\nu)}$ (cf. \eqref{eq: standard monomial description of a filtration} and \eqref{eq:proof of approximation lemma second expression}).
Therefore, as $n$ varies, there are only finitely many elements in $\BCrysZero(\MM)$ which may appear in \eqref{eq:proof of approximation lemma second expression} as $b_j$, depending only on $\rho$ and $(\mu,\nu)$ as in the case of $b''$.
    Write $\wt(M_j) = n\La_0 + \pi_j$. Since $\wt(M_j \ot b_j) = \wt(v) = \La^n_{\zeta,\eta}$ and there are only finitely many choices for $\wt(b_j)$, there are only finitely many elements $\pi \in P$ which may appear as $\pi_j$,
depending only on $\rho$ and $(\mu,\nu)$ as well.

\medskip
\textit{Step 3.}
Consider $x_i w_{M_j}$ and $y_i x_{b_j}$'s in \eqref{eq:proof of approximation lemma third expression}.

First, we take $N_1$ such that $q^{N_1} y_ix_{b_j} \in \LCrysZero(\MM)$ for all $i,j$. By \textit{Step 2}, we may assume that $N_1$ depends only on $\rho, (\mu,\nu)$, and $b''$.

Next, consider $x_i w_{M_j}$.
For a fixed $\pi \in P$, the set
\[
    \left\{ \, \varepsilon_k(M) \,|\, k\in \Z, M\in \BFn\text{ with } \wt(M) = n\La_0 + \pi \, \right\}
\]
is bounded above, whose upper bound is independent of $n$.
Indeed, for $M\in \BFn$ with $\wt(M) = n\La_0 + \pi$, we have $(\pi, \eps_k) \ge 0$ for $k> 0$ and $(\pi, \eps_k) \le 0$ for $k\le 0$. Considering the actions of $\te_k, \tf_k$ ($k\in \Z)$ on $\BFn = \BCrys(\FF)^{\ot n}$, it is not difficult to see that
\[
    \varepsilon_k(M) \le \begin{cases}
        (\pi, \eps_{k+1}) & \text{ for } k \ge 0,\\
        -(\pi, \eps_k) & \text{ for } k < 0. \\
    \end{cases}
\]
By applying \eqref{eq:approximating action of divided powers}, there exists $N''$ such that $q^{N''} f_k w_{M_j} \in \LCrys(\FFn)$ for all $k, j$, and $n$. Repeating similar arguments, we conclude that there exists $N_2$ depending only on $\rho, (\mu,\nu)$, and $b''$ such that $q^{N_2} c_jx_i w_{M_j} \in \LCrys(\FFn)$ for all $i,j$, and $n$.

\medskip
\textit{Step 4.} Finally, consider $t_{\rho_i} w_{M_j}$ in \eqref{eq:proof of approximation lemma third expression}.
Fix $i$ and let $l_{k,j}$ be such that $t_k^{(\rho_i, \La_k)} w_{M_j} = q^{l_{k,j}} w_{M_j}$.
By \textit{Step 2}, $\{\, l_{k,j} \,|\, k,j \ (k\ne 0)\, \}$ has a lower bound independent of $n$.
On the other hand, the projection of $y_i x_{b_j}$ onto $\MM_{< -d}$ is nonzero if and only if
$(\wt(y_ix_{b_j}), \La_0) < -d$, which is equivalent to
\[
    -(\rho_i, \La_0) = (\wt(y_i), \La_0 ) < -d - (\wt(b_j), \La_0) \le -d - (\wt(b'), \La_0) \le -d - (\wt(b), \La_0) \le 0.
\]
So we have $(\rho_i,\La_0) > 0$ in this case.
Furthermore, this implies that $l_{0,j}$ can be arbitrarily large as $n\to \infty$, since 
the number of columns in $M_j$ not equal to $|0\rangle$ is bounded above by $-(\pi_j, \La_0)$, and the number of columns equal to $|0\rangle$ is arbitrarily large.

Therefore, combining with \textit{Step 3}, we conclude that there exists $N=N(d,\rho,M)$ such that $\val^n(x^\circ_{< -d}) > M$ for $n\ge N$.
This completes the proof.
\qed

\medskip

The following lemma says that an element of $\FFn\ot (\MM_{\ge -d})$ can be approximated by an element of $(\FFnM)_{\ge -d}$ arbitrarily closely with respect to the $q$-adic topology in $\FFnM$.

\begin{lem} \label{cor:arbitrarily close approximation of elements of FFinftyM}
    Let $d,\rho$ and $M$ be given as above. 
    The following statement holds for all sufficiently large $n$; For $x\in \FFn\ot (\MM_{\ge -d})$ of weight $n\Lambda_0 + \rho$,
    there exists $x^\circ \in (\FFnM)_{\ge -d}$ such that $\val^n (x) = \val^n(x^\circ)$ and $\val^n(x) + M < \val^n (x - x^\circ)$.
\end{lem}
\pf
Suppose that $n \ge N(d,\rho,M)$, where $N(d,\rho,M)$ is as in \Cref{lem:the approximation lemma}. We may assume that $\val^n(x) = 0$.
By using \Cref{lem:the approximation lemma}, we define $x^{(0)}, x^{(1)}, \cdots, x^{(M)}\in (\FFnM)_{\ge -d} \cap \LCrysZero(\FFnM)_{n\La_0+\rho}$ inductively by
\[q^{-k-1}\left(x - \sum_{i=0}^{k} q^{i}x^{(i)}_{\ge -d}\right) \equiv x^{(k+1)} \pmod{q\LCrysZero(\FFnM)},\] and $\val^n(x^{(k+1)}_{< -d}) > M - k-1$.
Put $x^\circ = \sum_{i=0}^{M} q^{i}x^{(i)}$.
Then $\val^n(x) = \val^n(x^\circ)$ and
\[\val^n\left(x - \sum_{i=0}^{M} q^{i}x^{(i)}\right) \ge \min\left\{\, \val^n \left( x - \sum_{i=0}^{M} q^{i}x^{(i)}_{\ge -d} \right), \val^n \left(\sum_{i=0}^{M} q^{i}x^{(i)}_{< -d} \right) \,\right\} > M.\]
\qed
\medskip

\newcommand{\rzero}{{}_0r^-}
\newcommand{\rzerobar}{{}_0\overline{r}^-}

Next, we prove a lemma on a valuation of a projection on certain $\UqslZero$-modules, which is analogous to \Cref{lem:q-order of projection onto isotypic components of Fn}.

\begin{lem}\label{lem:UqslZero projection valuation diverges}
    Let $(\mu,\nu)  , (\zeta,\eta)\in \cP^2$ be given with $(\mu,\nu) \le (\zeta,\eta) $. Let $\uppi^n_{\zeta,\eta}$ be the projection of $V(\La^n_{\mu,\nu})\ot \MM$ onto its isotypic component
    having highest weight $\La_{\zeta,\eta}$.
    Let $\val^n_{\mu,\nu}$ be the crystal valuation associated to $\LCrys(\La^n_{\mu,\nu})\ot \LCrysZero(\MM)$.
    Let $(\sigma,\tau)\in \cP^2$ be given such that $(\zeta,\eta) \lneq  (\sigma,\tau)$ 
    with $\gamma = \La^n_{\sigma,\tau} - \La^n_{\mu,\nu}$. 
    Then for $x\in \MM_\gamma$,
    \[
        \lim_{n\to \infty} \val^n_{\mu,\nu}(\uppi^n_{\zeta,\eta}(u_{\La^n_{\mu,\nu}} \ot x)) = \infty.
    \]
\end{lem}

\medskip
\pf
Recall that there exist $q$-Shapovalov forms on $V(\La^n_{\mu,\nu})$ in \eqref{eq:q-Shapovalov form on uqg} and $\MM$ in \Cref{prop:q-Shapovalov form on parabolic Verma modules}.
Define a bilinear form $\langle\cdot, \cdot \rangle$ on $V(\La^n_{\mu,\nu})\ot \MM$ by
$\langle v_1\ot x_1 , v_2\ot x_2 \rangle = \langle v_1,v_2\rangle\langle x_1,x_2\rangle$ for $v_i\ot x_i\in V(\La^n_{\mu,\nu}) \ot \MM$.
By \Cref{prop:distributivity of tau on parabolic boson algebra}, we have $\langle uv, w\rangle = \langle v, \tau(u)w\rangle$ for $v,w\in V(\La^n_{\mu,\nu})\ot \MM$ and $u\in \UqslZero$. 
Since $\LCrys(V(\La^n_{\mu,\nu}))\ot \LCrysZero(\MM)$ is a crystal lattice, 
we also have an analogue of 
\Cref{thm:characterization of crystal lattice of parabolic Verma modules in terms of q-Shapovalov form} for $\LCrys(V(\La^n_{\mu,\nu}))\ot \LCrysZero(\MM)$.

Let $\gamma' = \La^n_{\zeta,\eta} - \La^n_{\mu,\nu}$.
We may assume that $H^n_{\gamma'}(\mu,\nu) = H_{\gamma'}(\mu,\nu)$ by letting $n$ sufficiently large.
For $b\in H^n_{\gamma'}(\mu,\nu)$, 
let $v_b$ the projection of $u_{\La^n_{\mu,\nu}} \ot x_b$ onto the isotypic component of $V_0(\La_{\zeta,\eta})$.
One can easily see that $v_b$ is singular, $v_b \in \LCrys(V(\La^n_{\mu,\nu}))\ot \LCrysZero(\MM)$,
and 
\begin{equation}\label{eq:UqslZero projection expression 1}
    v_b - u_{\La^n_{\mu,\nu}} \ot x_b \in V(\La^n_{\mu,\nu}) \ot \MM_{> \gamma}.
\end{equation}
By similar arguments as in the proof of \Cref{lem:q-order of projection onto isotypic components of Fn}, we have
\[
    \val^n_{\mu,\nu}(\uppi^n_{\zeta,\eta}(u_{\La^n_{\mu,\nu}} \ot x)) = \min \left\{\, \langle G(b')v_b, u_{\La^n_{\mu,\nu}} \ot x \rangle \,\middle|\, b\in H_{\gamma'}(\mu,\nu), b'\in \BCrys(\infty)_{\gamma - \gamma'} \,\right\}.
\]
We see from \eqref{eq:UqslZero projection expression 1} that the weight of the first tensor component in $G(b')(v_b - u_{\La^n_{\mu,\nu}} \ot x_b)$ is strictly less than $\La^n_{\mu,\nu}$.
Thus, this term cannot contribute to the pairing with $u_{\La^n_{\mu,\nu}} \ot x$.
Similarly, the only term in $G(b')(u_{\La^n_{\mu,\nu}} \ot x)$ that may contribute to the pairing is $t_{-\gamma + \gamma'}u_{\La^n_{\mu,\nu}} \ot G(b')x_b$.
Thus, we have $\langle G(b')v_b, u_{\La^n_{\mu,\nu}} \ot x \rangle = q^{(-\gamma + \gamma', \La^n_{\mu,\nu})}\langle G(b')x_b, x\rangle$, and hence
\[
    \val^n_{\mu,\nu}(\uppi^n_{\zeta,\eta}(u_{\La^n_{\mu,\nu}} \ot x)) = \min \left\{\,
    q^{(-\gamma + \gamma', \La^n_{\mu,\nu})}\langle G(b')x_b, x\rangle \,\middle|\, b\in H_{\gamma'}(\mu,\nu), b'\in \BCrys(\infty)_{\gamma - \gamma'} \,\right\}.
\]
Since $(\zeta,\eta) \lneq  (\sigma,\tau)$, we have $(-\gamma+\gamma', \La_0) \ge 1$, and hence the above value diverges as $n\to \infty$.
\qed

\medskip

The following lemma, analogous to \Cref{lem:subquotient_injectivity}, is a key ingredient to the proof of our main theorem.
It essentially depends on an analogous property of ${}_0r^{-}$ proved in \Cref{lem:ri-annihilator is 1}.

\begin{lem}\label{lem:derivation_lowest_degree_coefficient}
    Let $(\mu,\nu), (\zeta,\eta)\in \cP^2$ be given with $(\mu,\nu) \lneq (\zeta,\eta)$ and $d = |\zeta| - |\mu| = |\eta| - |\nu|$, and let $\gamma = \La^n_{\zeta,\eta} - \La^n_{\mu,\nu}$.
    Then the following statements hold for all sufficiently large $n$ and $N$;
    \begin{enumerate}
        \item The map
        \[
            \pi^{n+N}_{-d+1} \circ \phi_{n,n+N}: (\FFnM)^{(\zeta,\eta)}_{(\mu,\nu)} \longrightarrow (\FF^{ n + N} \ot \MM)_{-d + 1}
        \]
        is injective.
        \item For all $x\in (\FFnM)^{(\zeta,\eta)}_{(\mu,\nu)}$, we have
        \[
            \val^{n+N}_{-d+1}\left( \phi_{n,n+N}(x) \right) \in \val^n(x) + N(d-1) + [r_1, r_2],
        \]
        where $[r_1, r_2]$ is an interval independent of $x$, $n$ and $N$, depending only on $(\mu,\nu)$ and $(\zeta,\eta)$.
    \end{enumerate}
\end{lem}

\medskip

\pf
Let $N$ be given.

\medskip
\textit{Step 1.}
Let $V_N$ be the irreducible highest weight $\QUE{\gl_\infty}$-module with highest weight $N\La_0$.
By \Cref{prop:comodule}, one can check the following commutative diagram of $\QUE{\gl_\infty}[-]$-linear maps:
\begin{equation}\label{eq:expressing phiN in terms of rzerominus}
    \begin{tikzcd}[column sep=large]
	\MM & {\MM \ot \MM} \\
	\MM & {V_N\ot \MM} \\
	\MM & {\FF^{ N} \ot \MM}
	\arrow["{\Delta}", from=1-1, to=1-2]
	\arrow[equal, from=1-1, to=2-1]
	\arrow["{(\pi_{\Z\setminus 0, N\Lambda_0} \ot 1) \circ \Xi_N}", from=1-2, to=2-2]
	\arrow[from=2-1, to=2-2]
	\arrow[equal, from=2-1, to=3-1]
	\arrow[hook', from=2-2, to=3-2]
	\arrow["{\phi_{0,N}}", from=3-1, to=3-2]
\end{tikzcd}
\end{equation}
where $\pi_{\Z\setminus 0, N\Lambda_0} :\MM = V_{\Z\setminus 0}(0) \to V(N\Lambda_0)$ is the canonical projection \eqref{eq:pi lower J lambda},
$\Xi_N(x\ot y) =  q^{-N ( \wt y, \Lambda_0 )} x\ot y$,
and $V_N \hookrightarrow \FF^{N}$ is the unique embedding sending $u_{N\Lambda_0}\mapsto |0\rangle^{\ot N}$. 

On the other hand, let $W$ be the ireducible highest weight $\UqglInfMinusZero$-module with highest weight $-\alpha_0$, which is isomorphic to $\MM_{-1}\subset \MM$ with highest weight vector $f_0$.
Let $\pi_W :\MM \to W$ be the canonical projection which is $\UqglInfMinusZero$-linear. Recall that
\[\begin{tikzcd}[column sep=large]
	{r:= {}_0r^-: \MM} & {\MM\ot \MM} & {W \ot \MM}
	\arrow["{\Delta}", from=1-1, to=1-2]
	\arrow["{\pi_W\ot \mathrm{id}}", from=1-2, to=1-3]
\end{tikzcd}\]
(cf. \eqref{eq:def of derivations}). 
Since the map $W\cong \MM_{-1} \to V_N \hookrightarrow \FF^{ N}$ is injective, where the first map is ${\pi_{\Z\setminus 0, N\La_0}}$,
we may regard $r$ as a map
\begin{equation}\label{eq:identifying rzerominus with part of phiN}
    r: \MM \longrightarrow \FF^{ N} \ot \MM.
\end{equation}

\medskip
\textit{Step 2.}
Note that $(\gamma,\La_0) = -d$ by assumption.
We may assume that $H^n_\gamma(\mu,\nu) = H_\gamma(\mu,\nu)$. Let $\phi = \pi^{n+N}_{-d+1} \circ \phi_{n,n+N}$.
Let $w_{\Mnlm} \ot x_b \in (\FFnM)^{\le (\zeta,\eta)}_{(\mu,\nu)}$ be given for $b\in H_\gamma(\mu,\nu)$.

Consider $\phi(w_{\Mnlm} \ot x_b)$.
Since $\MM = \bigoplus_{d\ge 0} \MM_{-d}$, we have by \eqref{eq:expressing phiN in terms of rzerominus} and \eqref{eq:identifying rzerominus with part of phiN}
\[
    \phi_{0,N}(x_b) = q^{Nd} |0\rangle^{\ot N} \ot x_b + q^{N(d-1)} r(x_b) + y,
\]
for some $y\in \FF^N \ot (\MM_{\ge -d + 2})$, while $r(x_b) \in \FF^N\ot (\MM_{-d+1})$, and hence
\[
\begin{split}
    \phi_{n,n+N}&(w_{\Mnlm} \ot x_b) \\ 
    & = q^{Nd} w_{\Mnlm} \ot  |0\rangle^{\ot N}\ot x_b  + q^{N(d-1)} w_{\Mnlm} \ot r(x_b) + w_{\Mnlm} \ot y.
\end{split}
\]

Let $N_b$ be such that $q^{-N_b}r(x_b) \in \LCrysZero(\FFinftyM)$ and $q^{-N_b}r(b) \equiv \overline{r}(x_b) \pmod{q\LCrysZero(\FFinftyM)}$ (cf. \eqref{eq:rzerominus induces crystal map up to scaling}, \eqref{eq:crystal morphism induced by r-}).
We choose $N$ to be sufficiently large so that $Nd > N(d-1) + N_b > 0$. Note that this choice of $N$ is independent of $n$.

By \Cref{cor:arbitrarily close approximation of elements of FFinftyM}, 
for all sufficiently large $n$, there exists $y^\circ \in (\FF^{ n+N} \ot \MM)_{\ge -d+2}$ such that
\[
    \val^{n+N}\left(w_{\Mnlm} \ot y - y^\circ\right) > \val^{n+N}\left(w_{\Mnlm} \ot y\right) + Nd \ge Nd.
\]
Put $X = \phi_{n,n+N}(w_{\Mnlm} \ot x_b) - y^\circ$.
Then $X$ and $\phi_{n,n+N}(w_{\Mnlm} \ot x_b)$ have the same image under $\pi^{n+N}_{-d+1}$, and
\begin{equation} \label{eq:non-vanishing thm lowest deg coeff}
    q^{-N(d-1) - N_b} X \equiv \Mnlm \ot \overline{r} (b) \pmod{q\LCrysZero(\FF^{n+N}\ot \MM)}.
\end{equation}

Since $\overline{r}$ commutes with $\te_i,\tf_i$ for $i\ne 0$, we have $\te_i (\Mnlm \ot \overline{r}(b)) = 0$ for $i\ne 0$.
Moreover, if $n$ is sufficiently large, $\te_0(\Mnlm \ot \overline{r}(b)) = 0$ by tensor product rule \eqref{eq:tensor product rule}.
So $\Mnlm \ot \overline{r}(b)$ is connected to 
$M^{n+N}(\mu_b,\nu_b) \ot b' \in \HCrysZero{(\FF^{n+N} \ot \MM)^{(\zeta,\eta)}_{(\mu_b,\nu_b)}}$ as a $\UpglPlus[n+N]$-crystal for some $(\mu_b,\nu_b)$ and $b' \in H_{\gamma'}(\mu_b,\nu_b)$ with $\gamma' = \La^{n+N}_{\zeta,\eta} - \La^{n+N}_{\mu_b,\nu_b}$.

By construction of $\overline{r}$ and the fact that $\UpglPlus[n+N]$-crystal operators do not affect the tensor component in $\MM$, $b'$ depends only on $b$ and $(\wt(b'),\La_0) = -d + 1$. This implies that $X$ has non-zero image under the projection onto $(\FF^{n+N}\ot \MM)^{(\zeta,\eta)}_{(\mu_b,\nu_b)} \subset (\FF^{n+N}\ot \MM)_{-d+1}$, and hence so does $\phi_{n,n+N}(w_{\Mnlm}\ot x_b)$, that is, $\phi(w_{\Mnlm}\ot x_b)\ne 0$.

\medskip
\textit{Step 3.}
For each $z\in \BCrysZero:=\BCrysZero\left( (\FFnM)^{(\zeta,\eta)}_{(\mu,\nu)} \right)$, choose sequences of crystal operators
$\tilde{F}_z = \tf_{i_{z,1}}\cdots \tf_{i_{z,x_{z}}}$ and
$\dot{\tilde{F}}_z = \dot{\tf}_{j_{z,1}}\cdots\dot{\tf}_{j_{z,y_{z}}}$
such that
\[
z = \tilde{F}_z \dot{\tilde{F}}_z (\Mnlm \ot b_z),
\]
for some $b_z \in H_\gamma(\mu,\nu)$.
For $b\in H_\gamma(\mu,\nu)$, let $v_b = v_{\Mnlm\ot b}$ be the projection of $w_{\Mnlm}\ot x_b$ onto 
$(\FFnM)^{(\zeta,\eta)}_{(\mu,\nu)}$ (cf. the proof of \Cref{lem:the approximation lemma}).

Suppose a non-zero weight vector $x\in (\FFnM)^{(\zeta,\eta)}_{(\mu,\nu)}$ is given.
We have
\begin{equation}\label{eq:non-vanishing thm choice of x}
    x = \sum_{z \in \BCrysZero} c_z F_z \dot{F}_z v_{b_z},
\end{equation}
for some $c_z\in \Q(q)$ and $1\le j_{z,k} \le n$, 
since $\{\tilde{F}_z\dot{\tilde{F}}_zv_{b_z}\}_{z\in \BCrysZero}$ form a basis of $(\FFnM)^{(\zeta,\eta)}_{(\mu,\nu)}$ at $q = 0$, and the weight spaces of $\FFnM$ are finite-dimensional.

Let us consider
\[
    y = \sum_{z\in \BCrysZero} c_z \tilde{F}_z \dot{\tilde{F}}_z \left( w_{\Mnlm} \ot x_{b_z} \right),
\]
which will be used as an approximation of $x$.
Then $x - y \in (\FFnM)^{ < (\zeta,\eta)}_{(\mu,\nu)}$, and $\phi_{n,n+N}(x - y) \in (\FF^{n+N}\ot \MM)^{<(\zeta,\eta)}_{\ge (\mu,\nu)}$ (cf. \Cref{lem:application of branching rule to filtration}).
Thus, we have
\[ \pi_{-d+1}^{n+N}(\phi_{n,n+N}(x-y)) = \pi_{-d+1}^{n+N}(\phi_{n,n+N}(\pi_{-d+1}^{n}(x-y))), \]
and
\begin{equation}\label{eq:non-vanishing thm bounding the garbage term}
    \begin{split}
\val^{n+N}_{-d+1}(\phi_{n,n+N}(x-y)) & = \val^{n+N}(\pi^{n+N}_{-d+1}( \phi_{n,n+N}(x-y)) ) \\
 & = \val^{n+N}(\pi^{n+N}_{-d+1}( \phi_{n,n+N}(\pi^n_{-d+1}(x-y))) ) \\
 & = \val^n(\pi_{-d+1}^n(x-y)) + N(d-1) \quad\text{by \Cref{lem:subquotient_injectivity}}
    \end{split}
\end{equation}
for all sufficiently large $n$, depending only on $(\zeta,\eta)$.

\medskip
\textit{Step 4.}
Next, consider $\phi(y)$.
Let $z_1,\cdots, z_l \in \BCrysZero$ be such that
\[
    N(d-1) + N_{b_{z_i}} + \val(c_{z_i}) = \min \left\{\, N(d-1)+N_{b_{z}} + \val(c_{z}) \,\middle|\, z \in\BCrysZero \,\right\} =: m \quad (i=1,\ldots, l),
\]
where $N_b$ is defined in \textit{Step 2}.
Since $H_\gamma(\mu,\nu)$ is finite, we may choose $N$ sufficiently large so that assumptions in \textit{Step 2} are satisfied for all $b\in H_\gamma(\mu,\nu)$.
Then, by the same argument as in \textit{Step 2} (see \eqref{eq:non-vanishing thm lowest deg coeff}), we have
\begin{equation}\label{eq:non-vanishing thm crystal expansion}
    q^{-m} \phi(y) \equiv \sum_{i=1}^l q^{-\val(c_{z_i})} c_{z_i}
    F_{z_i} \dot{F}_{z_i}
    \left( \Mnlm \ot \overline{r}(b_{z_i}) \right) \pmod{q\LCrysZero(\FF^{n+N} \ot \MM)}.
\end{equation}
On the other hand, as in \textit{Step 2}, if $n$ is sufficiently large, then $\Mnlm\ot \overline{r}(b)$ is of highest weight for all $b\in H_\gamma(\mu,\nu)$ as an element of $\UqslZero$-crystal, and there exists a morphism of $\UqslZero\ot \UpglPlus[n]$-crystals

\[\begin{tikzcd}[row sep=tiny]
	{\psi:\BCrysZero\left( (\FFnM)^{(\zeta,\eta)}_{(\mu,\nu)}\right) } & {\BCrysZero\left(\FF^{n+N}\ot \MM\right)} \\
	{\Mnlm \ot b} & {\Mnlm \ot \overline{r}(b)}
	\arrow[from=1-1, to=1-2]
	\arrow[maps to, from=2-1, to=2-2],
\end{tikzcd}\]
which is injective since $\overline{r}$ is injective.
In particular, since the summands of the right hand side of \eqref{eq:non-vanishing thm crystal expansion} are linearly independent, we have $\phi(y) \ne 0$ and $\val^{n+N}(\phi(y)) = m$.

\medskip
\textit{Step 5.}
We may assume that $\val^n(x) = 0$ so that $\min \{\,\val(c_z) \,|\, z\in \BCrysZero \, \} = 0$. 
Let $r_1 = \min\{ N_b \,|\, b\in H_\gamma(\mu,\nu)\}$ and $r_2 = \max \{N_b \,|\, b\in H_\gamma(\mu,\nu)\}$.
Then we have
\begin{equation}\label{eq:non-vanishing thm bounding m}
    N(d-1) + r_1 \le m  \le N(d-1) + r_2,
\end{equation}
since $\val(c_z)\ge 0$ for $z\in \BCrysZero$ and $m\le N(d-1) + N_{b_{z'}} + \val(c_{z'}) = N(d-1) + N_{b_{z'}} \le N(d-1) + r_2$ where $z'\in \BCrysZero$ is such that $\val(c_{z'}) = 0$.

By \eqref{eq:UqslZero projection expression 1} and then applying \Cref{lem:UqslZero projection valuation diverges}, we see that if we take a sufficiently large $n$, then $\val^n(\pi^n_{-d+1}(w_{\Mnlm} \ot x_b - v_b)) > r_2$ for all $b\in H_\gamma(\mu,\nu)$. Combining with \eqref{eq:non-vanishing thm bounding the garbage term} and \eqref{eq:non-vanishing thm bounding m}, it follows that 
\[ \val^{n+N}_{-d+1}(\phi_{n,n+N}(x-y)) > m = \val^{n+N}(\phi(y)). \]
Since $\phi(x) = \phi(y) + \pi_{-d+1}^{n+N} (\phi_{n,n+N}(x - y))$, we conclude that $\val^{n+N}(\phi(x)) = m$.

In particular, we have $\phi(x)\ne 0$, and $\val^{n+N}(\phi(x)) \in N(d-1) + [r_1, r_2]$ by \eqref{eq:non-vanishing thm bounding m}. This proves (1) and (2).
\qed

\begin{rem}
    We emphasize that there exist $n'$ and $N'$, depending only on $(\mu,\nu)$ and $(\zeta,\eta)$,  and not depending on each other such that the statements (1) and (2) of \Cref{lem:derivation_lowest_degree_coefficient} hold for all $n\ge n'$ and $N\ge N'$.
\end{rem}

\begin{rem}\label{rem:consequences of non-vanishing thm in terms of lowest degree coefficients}
    For a weight vector $x\in \FFnM$ with $\val^n(x) = l$,
    we have $q^{-l}x \equiv x_0 \pmod{q\LCrysZero(\FFnM)}$ for some $x_0 \in \LCrysZero(\FFnM)/q\LCrysZero(\FFnM)= \Q$-span of $\BCrysZero(\FFnM)$.
    Let us call $x_0$ the \textit{lowest degree coefficient of $x$}.
    We have seen in the proof of \Cref{lem:derivation_lowest_degree_coefficient}, more precisely, by \eqref{eq:non-vanishing thm crystal expansion}, that the the lowest degree coefficient of $\pi^{n+N}_{-d+1}\circ \phi_{n,n+N}(w_{M^n(\mu,\nu)} \ot x_b)$ is in the $\Q$-span of the connected component of $M^n(\mu,\nu) \ot \overline{r}(b)$ ($b\in H_\gamma(\mu,\nu)$) with respect the crystal operators for $\UqslZero\ot U_p(\gl_n)$.
\end{rem}

For $n\ge 1$ and $(\mu,\nu)\in \cP^2$, let $\pi_{(\mu,\nu)}:\FFnM \to (\FFnM)_{(\mu,\nu)}$ denote the canonical projection.

\begin{lem}\label{lem:definition of Wnx}
    For $x\in (\FFnM)^{(\zeta,\eta)}$, let
    \begin{equation*}
        W^n(x) = \left\{ \, (\sigma,\tau) \in \cP^2 \,\middle|\, (\sigma,\tau) \le (\zeta,\eta), \ \pi_{(\mu,\nu)}(x) \ne 0 \text{ for some } (\mu,\nu) \le (\sigma,\tau) \,\right\}.
    \end{equation*}
    Then we have $W^{n+1}\left( \phi_{n,n+1}(x) \right) \subset W^n(x)$.
\end{lem}

\pf
If $(\mu^l,\nu^l)$ ($l=1, \cdots, k$) are the minimal elements in $W^n(x)$, then $x \in \sum_{l=1}^k (\FFnM)^{(\zeta,\eta)}_{\ge (\mu^l,\nu^l)}$.
By \Cref{lem:application of branching rule to filtration}, we have $\phi_{n,n+1}(x) \in \sum_{l=1}^k (\FFnplusM)^{(\zeta,\eta)}_{\ge (\mu^l,\nu^l)}$, which proves the claim.
\qed

\medskip

\noindent{\bfseries Proof of \Cref{thm:valuation limit exists}.}
Suppose that a nonzero $x\in \FFinftyM$ is given. Let $x_n\in \FFnM$ be given for all sufficiently large $n$ such that $\phi_n(x_n) = x$ and $\phi_{n,n+1}(x_n) = x_{n+1}$.
We may also assume that $x_n \in (\FFnM)^{(\zeta,\eta)}$ for some $(\zeta,\eta)\in \cP^2$.
For each $n$, we put $x_{n, -d} = \pi^n_{-d}(x_n)$ for simplicity. Note that $\val^n_{-d}(x_n) = \val^n(x_{n,-d})$.

Let $W(x) = \bigcap_{n\ge1} W^n(x_n)$, which is a well-defined non-empty set by \Cref{lem:definition of Wnx}.
For $d\in \Z_{\ge 0}$, let $W_{-d}(x) = \left\{ \,(\sigma,\tau) \in W(x) \,|\, |\zeta| - |\sigma| = d\, \right\}$, and similarly define $W_{\ge -d}(x), W_{> -d}(x)$, and so on.

\medskip
\textit{Case 1.}
Suppose that $x \in (\FFinftyM)_{> -d}$, that is, $W(x) = W_{>-d}(x)$.
Then $x_{n,-d}=0$ for all sufficiently large $n$.
Hence $\val^n_{-d}(x_n) =\infty$ for all sufficiently large $n$, which implies that $\val^\infty_{-d}(x) = \infty$.

\medskip
\textit{Case 2.}
Suppose that $x\in (\FFinftyM) \setminus (\FFinftyM)_{> -d}$, that is, $W_{\le -d}(x) \ne \emptyset$.
We prove $\val^\infty_{-d}(x) \in \Z$ by using descending induction on $d$.
More precisely, we assume the following:
\begin{equation}\label{eq:main theorem induction hypothesis}
    \parbox{0.9\textwidth}{
        For all $e > d$ with $W_{\le -e}(x) \ne \emptyset$,
        the limit $\lim_{n\to \infty} \big( \val^n_{-e}(x_n) - dn \big)$ exists and takes value in $\Z$.
    }
\end{equation}

\medskip
\textit{Step 1.}
If $W_{< -d}(x) = \emptyset$, then $x_n \in (\FFnM)_{\ge -d}$ for all sufficiently large $n$, and hence there exists, $n'$, independent of $x$, such that $\val^{n+1}_{-d}(x_{n+1}) = \val^n_{-d}(x_n) + d$ for $n\ge n'$ by \Cref{lem:subquotient_injectivity}, which implies $\val^\infty_{-d}(x) \in \Z$.

Note that in this case, we have $\val^\infty_{-d}(x) = \val^n_{-d}(x) - nd$ for sufficiently large $n$. We remark that $n'$ depends only on $(\zeta,\eta)$ and not on $x$, while it depends on $x$ in \textit{Step 2}. 

\textit{Step 2.}
Suppose that $W_{< -d}(x) \ne \emptyset$.
Then  $\pi_{(\mu,\nu)}(x_n) \ne 0$ for some $n\in \Z_{\ge 0}$ and $(\mu,\nu)\in \cP^2$ with $(\mu,\nu) \le (\zeta,\eta)$ and $|\zeta| - |\mu| > d$.
Since there exists $(\sigma,\tau)\in \cP^2$ satisfying $(\mu,\nu)\le (\sigma,\tau) \le (\zeta,\eta)$ and $|\zeta| - |\sigma| = d+1$,
we have $(\sigma,\tau)\in W_{-d-1}(x)$ by definition, and hence 
$W_{\le -d-1}(x) \supseteq W_{-d-1}(x) \ne \emptyset$.
By \eqref{eq:main theorem induction hypothesis},
$\val^\infty_{-d-1}(x)$ is finite.
This, in particular, implies that $x_{n,-d-1} \ne 0$ for all sufficiently large $n$.

Fix a sufficiently large $N$ so that \Cref{lem:derivation_lowest_degree_coefficient} holds for all $(\mu,\nu) \in W_{-d-1}(x)$.
For simplicity, we write $\psi = \pi^{n+N}_{-d} \circ \phi_{n,n+N}$.
Since
\[
    x_n = x_{n, >-d} + x_{n, -d} + x_{n, -d-1} + x_{n, < -d-1},
\]
where $x_{n, > -d} = \sum_{ -d' > -d} x_{n, -d'}$ and $x_{n, < -d-1}$ is defined similarly,
we have
\begin{equation}\label{eq:expressing n+N as components from n}
    x_{n+N, -d} = \psi(x_{n, -d}) + \psi(x_{n, -d-1}) + \psi(x_{n, < -d-1}).
\end{equation}
Note that $\psi(x_{n, > -d}) = 0$ by \Cref{lem:partition filtration permanence}.

\textit{Step 3.}
We proceed with finding a bound of valuation of each component in \eqref{eq:expressing n+N as components from n}.
First, we have
\begin{equation}\label{eq:lower bound of d+2 term}
    \begin{aligned}
        \val^{n+N}_{-d}(\psi(x_{n, < -d-1})) & \ge \val^n( x_{n, < -d-1}) \ge \min \left\{ \, \val^n_{-d'}(x_{n, < -d-1}) \,\middle|\, d' \ge d+2 \,\right\} \\
        &= \min \left\{\, \val^\infty_{-d'}(x) + d'n \,\middle|\, d' \ge d+2\, \right\}
    \end{aligned}
\end{equation}
for all sufficiently large $n$.
Here, we have the first inequality since $\psi$ preserves crystal lattices, and the last equality by applying
\eqref{eq:main theorem induction hypothesis} to $e = d'$ if $W_{\le -d'}(x) \ne \emptyset$,
and \textit{Case 1} otherwise.

On the other hand, by applying \Cref{lem:derivation_lowest_degree_coefficient} to all $(\mu,\nu) \in W_{-d-1}(x)$, we see that there exists $r_1, r_2$ such that
\[
    \val^n_{-d-1}(x_{n, -d-1}) + dN + r_1 \le \val^{n+N}_{-d}(\psi(x_{n, -d-1})) \le \val^n_{-d-1}(x_{n, -d-1}) + dN + r_2,
\]
for all sufficiently large $n$.
By \eqref{eq:main theorem induction hypothesis}, we have
$\val^n_{-d-1}(x_{n, -d-1}) = \val^\infty_{-d-1}(x) + (d+1)n$ and hence
\begin{equation}
    \val^\infty_{-d-1}(x) + (d+1)n + dN + r_1 \le \val^{n+N}_{-d}(\psi(x_{n,-d-1})) \le \val^\infty_{-d-1}(x) + (d+1)n + dN + r_2.\label{eq:upper bound of d+1 term}
\end{equation}

\medskip
\textit{Step 4.}
We claim that
\begin{equation}\label{eq:estimating order of N-step projection}
    \val^{n+N}_{-d} (x_{n+N, -d})= \min \left\{\, \val^{n+N}_{-d} \left(\psi(x_{n, -d}) \right), \val^{n+N}_{-d} \left(\psi(x_{n,-d-1}) \right) \,\right\},
\end{equation}
for all sufficiently large $n$ (cf. \eqref{eq:expressing n+N as components from n}).

First, the lower bound of \eqref{eq:lower bound of d+2 term} is $(d+2)n + c$ ($c\in \Z$), since there are only finitely many $d'\ge d+2$ such that $\val^\infty_{-d'}(x) < \infty$,
while the upper bound in \eqref{eq:upper bound of d+1 term} is of the form $(d+1)n + c'$ ($c'\in \Z$).
This implies that 
\begin{equation}\label{eq:d+2 term is smaller than d+1 term}
\val^{n+N}_{-d}(\psi(x_{n,-d-1})) < \val^{n+N}_{-d}(\psi(x_{n, < -d-1})),
\end{equation}
for all sufficiently large $n$, and hence
\begin{equation}\label{eq:estimating order preliminary inequality}
    \val^{n+N}_{-d} (x_{n+N, -d}) \ge \min \left\{ \,\val^{n+N}_{-d} \left(\psi(x_{n, -d}) \right), \val^{n+N}_{-d} \left(\psi(x_{n,-d-1}) \right)\, \right\}.
\end{equation}

Now, we claim that the lowest degree coefficients of $\psi(x_{n,-d})$ and $\psi(x_{n, -d-1})$ (see \Cref{rem:consequences of non-vanishing thm in terms of lowest degree coefficients}) do not cancel each other in $\psi(x_{n, -d}) + \psi(x_{n, -d-1})$.
This implies that the inequality in \eqref{eq:estimating order preliminary inequality} becomes equality.

Indeed, the lowest degree coefficient of $\psi(x_{n,-d})$ is in the 
$\Q$-linear span of elements of $\BCrysZero\left( \FF^{n+N}\ot \MM \right)$ connected to
\begin{equation}\label{eq:ldc for xnd}
    \left\{ M^n(\mu,\nu) \ot |0\rangle ^{\ot N} \ot b \,\middle|\, (\mu,\nu) \le (\zeta,\eta),\  |\zeta|-|\mu| = d,\  b \in H_\gamma(\mu,\nu) \right\},
\end{equation}
under the crystal operators of $\UqslZero\ot \UpglPlus[n]$
(see the proof of \Cref{lem:subquotient_injectivity}, where $\gamma = \La_{\zeta,\eta} - \La_{\mu,\nu}$).
On the other hand, by \Cref{rem:consequences of non-vanishing thm in terms of lowest degree coefficients}, the lowest degree coefficient of $\psi(x_{n, -d-1})$ is in the $\Q$-span of elements of $\BCrysZero\left( \FF^{n+N}\ot \MM \right)$ connected to
\begin{equation}\label{eq:ldc for xnd+1}
    \left\{ \,M^n(\mu',\nu') \ot \overline{r}(b') \,\middle|\, (\mu',\nu') \le (\zeta,\eta),\  |\zeta| - |\mu'| = d + 1, \  b'\in H_\gamma(\mu',\nu')\, \right\},
\end{equation} 
under the crystal operators of $\UqslZero\ot \UpglPlus[n]$.
Note that $\overline{r}:\BCrysZero(\MM) \to \BCrys(\FF^N)\ot \BCrysZero(\MM)$ and $\overline{r}(b') = M\ot b''$ for some $M\ne |0\rangle^{\ot N}$ and $b''\in \BCrysZero(\MM)$.
Also recall that $M^n(\mu',\nu')\ot \overline{r}(b')$ is of highest weight as an element of a $\UqslZero$-crystal for sufficiently large $n$ (see \textit{Step 2} in the proof of \Cref{lem:derivation_lowest_degree_coefficient}).
Similarly, $M^n(\mu,\nu)\ot |0\rangle^{\ot N} \ot b$ is of highest weight as an element of a $\UqslZero$-crystal.
Hence, if $M^n(\mu',\nu') \ot \overline{r}(b')$ and $M^n(\mu,\nu)\ot |0\rangle^{\ot N} \ot b$ are connected with respect to the $\UqslZero\ot \UpglPlus[n]$-crystal operators, they must be connected solely through the crystal operators of $\UpglPlus[n]$, which is impossible, since the components $|0\rangle^{\ot N} \ot b$ and $M\ot b'$ are clearly different.
Note that they may still be connected with respect to $\UpglPlus[n+N]$-crystal operators, and the restriction to $\UpglPlus[n]$ is crucial for this argument to work.
This proves the claim and hence \eqref{eq:estimating order of N-step projection}. 

\medskip
\textit{Step 5.}
Now, we claim that
\begin{equation}\label{eq:estimating order of N-step projection 2}
    \val^{n+N}_{-d}(x_{n+N, -d}) = \val^{n+N}_{-d}(\psi(x_{n,-d})),
\end{equation}
for all sufficiently large $n$.
We observe from \eqref{eq:upper bound of d+1 term} and \eqref{eq:estimating order of N-step projection} that $x_{n+N, -d}\ne 0$ for all sufficiently large $n$.

By \Cref{lem:subquotient_injectivity},
\begin{equation}\label{eq:proof of main thm dN}
\val^{n+N}_{-d}(\psi(x_{n,-d})) = \val^n_{-d}(x_{n,-d}) + dN,
\end{equation} for all sufficiently large $n$, which implies
$\val^{n+N}_{-d}(x_{n+N, -d}) \le \val^n_{-d}(x_{n,-d}) + dN$ by \eqref{eq:estimating order of N-step projection}.
Suppose that this holds for $n\ge n_0$. Then the sequence $\{ \val^n_{-d}(x_{n,-d}) - dn \}$ takes maximum value for $n \le n_0 + N$, 
so there exists $C$ such that $\val^m_{-d}(x_{m,-d}) \le dm + C$ for all $m$.
Combining this with $m$ replaced by $n+N$ and \eqref{eq:upper bound of d+1 term}, we have
$\val^{n+N}_{-d}(x_{n+N,-d}) < \val^{n+N}_{-d}\left(\psi(x_{n,-d-1})\right)$ for all sufficiently large $n$, since the left hand side is bounded above by a linear function on $n$ with slope $d$, while the right hand side is bounded below by a linear function on $n$ with slope $d+1$.
Hence, it follows from \eqref{eq:estimating order of N-step projection} that 
\begin{equation}\label{eq:d+1 term is smaller than d term}
    \val^{n+N}_{-d}(x_{n+N, -d}) = \val^{n+N}_{-d}(\psi(x_{n,-d})) < \val^{n+N}_{-d}\left(\psi(x_{n,-d-1})\right)
\end{equation}
for all sufficiently large $n$, which proves \eqref{eq:estimating order of N-step projection 2}.

\medskip
\textit{Step 6.}
By \eqref{eq:estimating order of N-step projection 2} and \eqref{eq:proof of main thm dN}, we have \[\val^{n+N}_{-d}(x_{n+N, -d}) = \val^n_{-d}(x_{n,-d}) + dN,\] for all sufficiently large $n$.
We may repeat the preceeding argument with $N+1$ instead of $N$ to obtain $\val^{n+N+1}_{-d}(x_{n+N+1, -d}) = \val^n_{-d}(x_{n,-d}) + d(N+1)$ for all sufficiently large $n$. Combining these two equations, we conclude that the sequence $\{\val^n_{-d} (x_{n,-d})\}_{n\in \Z_{\ge 0}}$ satisfies \[\val^{n+1}_{-d}(x_{n+1,-d}) = \val^n_{-d}(x_{n,-d}) + d\] for all sufficiently large $n$.
Thus, $\val^\infty_{-d}(x)$ is well-defined and is finite.
This completes the induction step, and hence the proof of \Cref{thm:valuation limit exists}.
\qed

\begin{rem}\label{rem:consequences of the proof of valuation limit exists}
    Summarizing from \eqref{eq:d+2 term is smaller than d+1 term}, \eqref{eq:estimating order of N-step projection 2}-\eqref{eq:d+1 term is smaller than d term}, we see that given $x$,
    \begin{align*}
        \val^n_{-d}(x_{n,-d}) + dN &= \val^n_{-d}(\psi(x_{n,-d})) = \val^{n+N}_{-d}(x_{n+N, -d})
        < \val^{n+N}_{-d}(\psi(x_{n, -d-1})) \\
        & = \val^{n+N}_{-d}(\psi(x_{n, < -d})) < \val^{n+N}_{-d}(\psi(x_{n, < -d-1}))
    \end{align*}
    (cf. \eqref{eq:expressing n+N as components from n}) for all sufficiently large $n$ and $N$, where $N$ depends on $x$, and $n$ depends on $x$ and $N$.
\end{rem}

\subsection{Socle filtration of $\FFinftyM$}
\label{subsec:Socle filtration}
\newcommand{\FFinftyMSubqSoc}{\mathcal{V}^\circ_d}
\newcommand{\FFinftyMSubq}{\mathcal{V}_d}
For $d\in \Z_{\ge 0}$, let
\begin{equation}\label{eq:definition of FFinftyMSubq}
    \FFinftyMSubqSoc = \frac{(\FFinftyM)_{\ge -d}}{(\FFinftyM)_{>-d}}, \quad \FFinftyMSubq = \frac{\FFinftyM}{(\FFinftyM)_{>-d}}.
\end{equation}
Then we define
\begin{align*}
    \LCrys(\FFinftyMSubq) &= \left\{\, \overline{x}\in \FFinftyMSubq \,\middle|\, x\in \FFinftyM, \val^\infty_{-d} (x) \ge 0 \,\right\}, \\
    \LCrys(\FFinftyMSubqSoc) &= \FFinftyMSubqSoc \cap \LCrys(\FFinftyMSubq).
\end{align*}
\begin{thm} \label{thm:crystal lattice of a socle become isomorphic mod q}
    For $d\in \Z_{\ge 0}$, we have the following. 
    \begin{enumerate}
        \item $\LCrys(\FFinftyMSubq)$ is a crystal valuation.
        \item $\LCrys(\FFinftyMSubqSoc)$ is a crystal lattice.
        \item The map
        \[
            \LCrys(\FFinftyMSubqSoc)/q\LCrys(\FFinftyMSubqSoc) \longrightarrow \LCrys(\FFinftyMSubq) / q\LCrys(\FFinftyMSubq)
        \]
        induced from the inclusion $\FFinftyMSubqSoc \subset \FFinftyMSubq$ is an isomorphism of $\Q$-spaces,
        that is, $\LCrys(\FFinftyMSubq)$ is saturated with respect to $\FFinftyMSubqSoc$.
    \end{enumerate}
\end{thm}

\pf
(1) It is clear that $\val^\infty_{-d}$ induces a well-defined crystal valuation on $\FFinftyMSubq$ by definition of $\val^\infty_{-d}$ and \Cref{thm:valuation limit exists}.

(2) Since $\val^{\infty}_{-d}$ is also a crystal valuation on $\FFinftyMSubqSoc$, it is enough to show that $\LCrys(\FFinftyMSubqSoc)$ is a free $A_0$-module.

Let $(\zeta,\eta)\in \cP^2$ be given.
Let $(\FFnM)^{(\zeta,\eta)}_{\ge -d} = (\FFnM)^{(\zeta,\eta)} \cap (\FFnM)_{\ge -d}$.
Then $\frac{(\FFinftyM)^{(\zeta,\eta)}_{\ge -d}}{(\FFinftyM)^{(\zeta,\eta)}_{>-d}}$ is a direct limit of $\frac{(\FFnM)^{(\zeta,\eta)}_{\ge -d}}{(\FFnM)^{(\zeta,\eta)}_{>-d}}$.
Let \[ \LCrys^{\infty, (\zeta,\eta)} = \LCrys(\FFinftyMSubqSoc) \cap  \frac{(\FFinftyM)^{(\zeta,\eta)}_{\ge -d}}{(\FFinftyM)^{(\zeta,\eta)}_{>-d}}, \quad \LCrys^{n,(\zeta,\eta)} = \LCrys(\FFinftyM) \cap \frac{(\FFnM)^{(\zeta,\eta)}_{\ge -d}}{(\FFnM)^{(\zeta,\eta)}_{>-d}} .\]
Then by \Cref{lem:subquotient_injectivity} (1) and \Cref{thm:valuation limit exists},
\[
    \LCrys^{\infty,(\zeta,\eta)}= \varinjlim_n q^{-dn} \LCrys^{n,(\zeta,\eta)}.
\]
Note that $\LCrys^{n,(\zeta,\eta)}$ is a free $A_0$-module, as it is given by a restriction of a crystal lattice on $\FFinftyM$.
Any submodule of $\LCrys^{\infty,(\zeta,\eta)}$ with finite rank is free, since it is contained in $\LCrys^{n,(\zeta,\eta)}$ for some $n$.
Therefore, we conclude that $\LCrys^{\infty,(\zeta,\eta)}$ is free by \cite[Theorem 2.4.3]{KT18}, since it is countably generated as an $A_0$-module.

(3) Since the map is clearly injective, it is enough to show that it is surjective.
Suppose that $\overline{x}\in \FFinftyMSubq$ is given. We show that $\overline{x}\equiv \overline{x^\circ} \pmod{q\LCrys(\FFinftyMSubq)}$ for some $x^\circ\in (\FFinftyM)_{\ge -d}$.
We keep the notations and assumptions in the proof of \Cref{thm:valuation limit exists}.
So we assume that $x\in (\FFinftyM)^{(\zeta,\eta)} = \varinjlim_n (\FFnM)^{(\zeta,\eta)}$. 

For $m\ge n$, denote $\psi_m = \pi^m_{-d} \circ \phi_{n,m}$. We claim that for all sufficiently large $m$,
\[
    \val^m_{-d} (x_{m, -d}) < \val^m_{-d} \left( \psi_m(x_{n, < -d}) \right).
\]
Let us denote the above statement by $P(m,n)$.

\medskip
\textit{Step 1.}
We first prove that $P(m,n)$ implies $P(m+N,n)$ for all sufficiently large $m$ and $n$, and $N$, where lower bounds for $m$ and $n$ depend only on $N$, and does not depend on each other.

We have
\begin{equation}\label{eq:saturatedness expansion}
\begin{aligned}
\psi_{m+N}( x_{n, < -d} ) & = \pi^{m+N}_{-d}\circ \phi_{m,m+N}\left( \psi_m(x_{n, < -d}) + \pi^m_{< -d} (\phi_{n,m}(x_{n, < -d}))\right) \\
& = \pi^{m+N}_{-d} (\phi_{m,m+N} (\psi_m(x_{n, < -d}))) + \pi^{m+N}_{-d} (\phi_{m,m+N}(x_{m, < -d})),
\end{aligned}
\end{equation}
where we have used $\pi^{m+N}_{-d}(\phi_{m,m+N}(\pi^m_{> -d}(\phi_{n,m}(x_{n, < -d})))) = 0$ at the first equality.

For all sufficiently large $m$,
\begin{equation}\label{eq:saturatedness first bound}
    \begin{aligned}
    \val^{m+N}_{-d} \big( \pi^{m+N}_{-d}(\phi_{m,m+N}&(\psi_m(x_{n, < -d})))\big) 
     \\
    & = \val^m_{-d} \left( \psi_m(x_{n, < -d}) \right) + dN \quad \text{by \Cref{lem:subquotient_injectivity}}\\
    & > \val^m_{-d} (x_{m,-d}) + dN \quad \text{by $P(m,n)$}\\
    & = \val^{m+N}_{-d} (x_{m+N, -d}) \quad \text{by \Cref{thm:valuation limit exists}},
    \end{aligned}
\end{equation}
On the other hand, from \Cref{rem:consequences of the proof of valuation limit exists}, we have
\begin{equation}\label{eq:saturatedness second bound}
\val^{m+N}_{-d}( \pi^{m+N}_{-d}(\phi_{m,m+N}(x_{m, < -d}))) > \val^{m+N}_{-d}(x_{m+N, -d}),
\end{equation}
for all sufficiently large $m$ and $N$, where $m$ depends on $N$.
Combining \eqref{eq:saturatedness first bound} and \eqref{eq:saturatedness second bound}, we have $P(m+N,n)$ by \eqref{eq:saturatedness expansion}.

\medskip
\textit{Step 2.}
Recall that $P(n+N, n)$ holds for all sufficiently large $n$ and $N$ by \Cref{rem:consequences of the proof of valuation limit exists}.
Also, by applying \textit{Step 1} to $N,N+1,\cdots, 2N-1$, it follows that $P(m,n)$ implies $P(m+N,n), \cdots P(m+2N-1,n)$ for all sufficiently large $m$ and $N$ where $m$ depends on $N$. Thus, if $n$ is sufficiently large, then $P(m,n)$ holds for all $m \ge n+N$.
This proves the claim.

Now, let $x^\circ = \phi_n(x_{n, \ge -d})$. Then $x^\circ \in (\FFinftyM)_{\ge -d}$.
By $P(m,n)$, we have \[\val^m_{-d}(x_{m, -d}) < \val^m_{-d}(\psi_m(x_{n,< -d})) = \val^m_{-d}(\psi_m(x_n - x_{n, \ge -d})) =  \val^m_{-d}( x_m - x^\circ_m),\] for all sufficiently large $m$, where $x^\circ_m = \phi_{n,m}(x_{n,\ge -d})$. By taking a limit, we have $0 = \val^\infty_{-d}(x) < \val^\infty_{-d}(x - x^\circ)$, and hence $\val^\infty_{-d}(x) = \val^{\infty}_{-d}(x^\circ) = 0$ and $\overline{x} \equiv \overline{x^\circ} \pmod{q\LCrys(\FFinftyMSubq)}$. This completes the proof of the surjectivity.
\qed

\begin{rem}
    The $A_0$-modules $\LCrys(\FFinftyMSubq)$ and $\LCrys(\FFinftyMSubqSoc)$ are given in a non-constructive way. We do not know an effective algorithm to compute $\val^\infty_{-d}$.
\end{rem}

\begin{thm}
    For $d\in \Z_{\ge 0}$, $\FFinftyMSubqSoc$ has a crystal base $(\LCrys(\FFinftyMSubqSoc),\BCrys(\FFinftyMSubqSoc))$ with
    \[
        \BCrys(\FFinftyMSubqSoc) \cong \bigoplus_{
            \substack{
                (\mu,\nu) \le 
                (\zeta,\eta),\\
                |\zeta| - |\mu| = |\eta| - |\nu| = d
            }
        } \left( \BCrys_0(\Lambda_{\zeta,\eta}) \ot \BCrys_{\mu,\nu}\right)^{\oplus|H_\gamma(\mu,\nu)|},
    \]
    where $\gamma= \La_{\zeta,\eta} - \La_{\mu,\nu}$.
\end{thm}
\pf
By \Cref{thm:isotypic decomposition of subquotients}, we have
\[
\FFinftyMSubqSoc \cong \bigoplus_{\substack{(\mu,\nu) \le (\zeta,\eta) \\ |\zeta|-|\mu|=|\eta|-|\nu|=d}} \left( V_0(\Lambda_{\zeta,\eta}) \ot V_{\mu,\nu}\right)^{\oplus |H_\gamma(\mu,\nu)|}.
\]
By uniqueness of a crystal lattice for an integrable $\UqslZero$-module (\Cref{cor:uniqueness of crystal valuation} (2)) and an extremal weight modules $\UpglPlus$-module (\Cref{prop:uniqueness of crystal for direct sum of extremal weight modules}), 
any crystal lattice of $\FFinftyMSubqSoc$ is isomorphic to $\bigoplus \left(\LCrysZero(\La_{\zeta,\eta}) \ot \LCrys_{\mu,\nu}\right)^{\oplus |H_\gamma(\mu,\nu)|}$,
which has a crystal given above. The conclusion follows from \Cref{thm:crystal lattice of a socle become isomorphic mod q} (2).
\qed

\begin{rem}
    One may regard \Cref{thm:crystal lattice of a socle become isomorphic mod q} as a generalization of \Cref{thm:identifying crystal basis of FFinftyM with that of its socle}.
By \Cref{thm:crystal lattice of a socle become isomorphic mod q},
$\BCrys(\FFinftyMSubqSoc)$ is a $\Q$-basis of $\LCrys(\FFinftyMSubq)/q\LCrys(\FFinftyMSubq)$.
However, we do not know how to construct a basis of $\FFinftyMSubq$ lifting the crystal $\BCrys(\FFinftyMSubqSoc)$, while there exists a basis $B=\{x_b \,|\,b\in \BCrys(\FFinftyMSubqSoc)\}$ of $\FFinftyMSubqSoc$ such that $x_b \equiv b \pmod{q\LCrys(\FFinftyMSubqSoc)}$. 
\end{rem}

Now we are in a position to show that the filtration $\{\,(\FFinftyM)_{\ge -d}\,\}_{d\in \Z_{\ge 0}}$ is indeed a socle filtration.
We first consider the following general situation.

Let $\left( \{M_n\}_{n\ge 1}, \{ \phi_n \}_{n\ge 1}\right)$ be a directed system of $\Q(q)$-spaces with $M=\varinjlim_n M_n$ such that
\begin{enumerate}
    \item $M_n$ is a $\UqslZero\ot\UpglPlus[n]$-module, whose set of weights is finitely dominated,
    \item $M_n$ is a direct sum of $V_0(\Lambda_{\zeta,\eta})\ot V^n_{\mu,\nu}$ with finite multiplicity,
    \item $\phi_n: M_n \to M_{n+1}$ is an injective $\UqslZero\ot \UpglPlus[n]$-linear map.
\end{enumerate}
By (3), $M$ becomes a $\UqslZero\ot \UpglPlus$-module.
Let $\left( \{N_n\}_{n\ge 1}, \{ \phi_n|_{N_n} \}_{n\ge 1}\right)$ be a sub-directed system with $N = \varinjlim_n N_n$ satisfying the following:
\begin{enumerate}
    \item $N_n$ is a $\UqslZero\ot\UpglPlus[n]$-submodule of $M_n$,
    \item $N_n$ and $M_n/N_n$ do not have common isotypic component,
    \item $N$ is a semisimple $\UqslZero\ot\UpglPlus$-submodule of $M$,
    \item $M$ has a crystal valuation $\val$, which is saturated with respect to $N$.
\end{enumerate}

\begin{thm}\label{thm:crystal isomorphism implies socle}
    We have $N = \soc M$.
\end{thm}

\pf
Let $\LCrys(M) = L_\val$, and let $\LCrys(M') = M'\cap \LCrys(M)$ for any $\Q(q)$-subspace $M'$ of $M$.
We may regard $M_n\subset M$ since $\phi_n$ is injective for $n\ge 1$. 
Since $\LCrys(M_n)$ is a crystal valuation,
we may apply an analogue of \Cref{prop:uniqueness of crystal for direct sum of extremal weight modules} for $\UqslZero\ot \UpglPlus[n]$-modules $M_n$'s (see \Cref{rem:same holds for UqUp}).

Suppose that there exists a non-zero submodule $N'\subset M$ such that $N\cap N' = \{0\}$.
Put $N'_n = N'\cap M_n$ for $n\ge 1$.
By the condition (2) on $N$,  \Cref{cor:uniqueness of crystal valuation}, and \Cref{rem:same holds for UqUp},
we have $\LCrys(N_n\oplus N'_n) = \LCrys(N_n) \oplus \LCrys(N'_n)$ for $n\ge 1$.
Taking a direct limit, we have $\LCrys(N\oplus N') = \LCrys(N)\oplus \LCrys(N')$.
Since the inclusion $N\oplus N'\hookrightarrow M$ induces an embedding
\[
    \frac{\LCrys(N\oplus N')}{q\LCrys(N\oplus N')} \cong \frac{\LCrys(N)}{q\LCrys(N)} \oplus \frac{\LCrys(N')}{q\LCrys(N')} \longhookrightarrow \frac{\LCrys(M)}{q\LCrys(M)},
\]
and $\LCrys(N)/q\LCrys(N) \hookrightarrow \LCrys(M)/q\LCrys(M)$ is an isomorphism by \Cref{cor:saturated iff mod-q bijective}, we conclude that $\LCrys(N')/q\LCrys(N') = 0$.
This contradicts the fact that no nonzero element in $\LCrys(M)$ is infinitely divisible by $q$.
Therefore, $N$ is a maximal semisimple submodule.
\qed

\newcommand{\IsotypicSocle}[1]{\frac{(#1)^{\lambda,\mu}}{(#1)^{\lambda,\mu}_{> -d}}}
\begin{thm} \label{thm:it is the socle filtration}
    $\{\, (\FFinftyM)_{\ge -d} \,\}_{d\in \Z_{\ge 0}}$ is the socle filtration of $\FFinftyM$.
\end{thm}
\pf
It is equivalent to showing that $\FFinftyMSubqSoc = \soc \FFinftyMSubq$ for $d\in \Z_{\ge 0}$,
where $\FFinftyMSubqSoc$ and $\FFinftyMSubq$ are as in \eqref{eq:definition of FFinftyMSubq}.
For $(\zeta,\eta)\in \cP^2$, let
\[
    \begin{aligned}
        (\FFinftyM)^{(\zeta,\eta)}_{\ge -d} &= (\FFinftyM)^{(\zeta,\eta)} \cap (\FFinftyM)_{\ge -d}, \\
        (\FFinftyM)^{(\zeta,\eta)}_{> -d} &= (\FFinftyM)^{(\zeta,\eta)} \cap (\FFinftyM)_{> -d}.
    \end{aligned}
\]
Since $\FFinftyM = \bigoplus_{(\zeta,\eta)} (\FFinftyM)^{(\zeta,\eta)}$, it is enough to show
\begin{equation}\label{eq:socle of isotypic component}
\soc \left( \frac{(\FFinftyM)^{(\zeta,\eta)}}{(\FFinftyM)^{(\zeta,\eta)}_{> -d}} \right) = \frac{(\FFinftyM)^{(\zeta,\eta)}_{\ge -d}}{(\FFinftyM)^{(\zeta,\eta)}_{> -d}}.
\end{equation}
First, consider
\[
    \phi_{n,n+1}^d: \frac{(\FFnM)^{(\zeta,\eta)}}{(\FFnM)^{(\zeta,\eta)}_{> -d}} \longrightarrow \frac{(\FFnplusM)^{(\zeta,\eta)}}{(\FFnplusM)^{(\zeta,\eta)}_{> -d}}.
\]
If $n$ is sufficiently large, then the map $(\FFnM)^{(\zeta,\eta)}_{(\mu,\nu)} \to (\FFnplusM)^{(\zeta,\eta)}_{(\mu,\nu)}$ induced from $\phi_{n,n+1}^d$ is injective for all $(\mu,\nu)\le (\zeta,\eta)$ with $|\zeta| - |\mu| = |\eta|-|\nu| = d$ by \Cref{cor:subquotient_injectivity of isotypic components}, and we have
\[
    \frac{(\FFinftyM)^{(\zeta,\eta)}}{(\FFinftyM)^{(\zeta,\eta)}_{> -d}} = \varinjlim_{n}\frac{(\FFnM)^{(\zeta,\eta)}}{(\FFnM)^{(\zeta,\eta)}_{> -d}}.
\]

By letting \[M_n = \frac{(\FFnM)^{(\zeta,\eta)}}{(\FFnM)^{(\zeta,\eta)}_{> -d}}, \quad N_n = \frac{(\FFnM)^{(\zeta,\eta)}_{\ge -d}}{(\FFnM)^{(\zeta,\eta)}_{> -d}}, \] we can check that $M_n$ and $N_n$ satisfy all the other conditions, and then apply \Cref{thm:crystal isomorphism implies socle} to obtain \eqref{eq:socle of isotypic component}.
\qed

\begin{cor} \label{cor:saturated crystal lattice of FFinftyM}
For $d\in \Z_{\ge 0}$, $\LCrys(\FFinftyMSubq)$ is a saturated crystal valuation.
\end{cor}
\pf
It follows from \Cref{thm:crystal lattice of a socle become isomorphic mod q}.
\qed

\medskip
As an important application, we obtain the following, which is one of our main results.

\newcommand{\lohimod}[1][\ot]{V_{\emptyset,\nu} #1 V_{\mu,\emptyset}}
\newcommand{\hilomod}[1][\ot]{V_{\mu,\emptyset} #1 V_{\emptyset,\nu}}

\begin{thm}\label{cor:main application - crystal basis} 
    For $\mu,\nu\in \cP$, the $\UpglPlus$-module $V_{\mu,\emptyset} \ot V_{\emptyset,\nu}$ has a socle filtration with the Loewy length $\min\{ |\mu|, |\nu| \} + 1$, and its subquotients are given by
    \begin{equation}\label{eq:main application - socle quotients}
        \frac{\soc^{d+1} \left(\hilomod\right)}{\soc^{d} \left(\hilomod\right)} \cong \bigoplus_{
            \substack{
                (\mu,\nu) \ge (\zeta,\eta) \\
                |\mu| - |\zeta| = |\nu| - |\eta| = d
            }
        } V_{\zeta,\eta}^{\oplus n^{\mu,\nu}_{\zeta,\eta}},
    \end{equation}
    for $0\le d\le \min\{|\mu|,|\nu|\}$, where $n^{\mu,\nu}_{\zeta,\eta} = |H_\gamma(\zeta,\eta)|$ with $\gamma = \La_{\mu,\nu} - \La_{\zeta,\eta}$.
    In particular, $\hilomod$ is indecomposable.
\end{thm}

\pf
By \Cref{thm:isotypic decomposition of F_inftyM} and \Cref{thm:it is the socle filtration},
we have
\begin{equation}\label{eq:identifying multiplicity spaces}
    \soc\left( \frac{\hilomod}{\soc^d(\hilomod)}\right) \cong \Hom_{\UqslZero}\left( V_0(\La_{\mu,\nu}) , \FFinftyMSubqSoc \right).
\end{equation}
By \Cref{thm:isotypic decomposition of subquotients}, the right-hand side of \eqref{eq:identifying multiplicity spaces}
is isomorphic to the right-hand side of \eqref{eq:main application - socle quotients}. Since the socle of $\hilomod$ is simple,
we conclude that $\hilomod$ is indecomposable.
\qed

\begin{rem}
    The tensor product in \Cref{cor:main application - crystal basis} is $\ot = \ov{\ot}_+$ by our convention \eqref{eqref:our choice}.
    The socle filtration of $\lohimod$ has the same subquotients (cf. \Cref{cor:identical filtration on lowest tensor highest} and \Cref{rem:R-matrix on lohimod}).
    Moreover, combining with \Cref{prop:filtration on highest tensor lowest} and \eqref{eq:tensor multiplicity}, we have 
    a combinatorial formula for $n^{\mu,\nu}_{\zeta,\eta}$ as follows:
    \begin{equation}\label{eq:identifying_m_and_n}
        n^{\mu,\nu}_{\zeta,\eta} = m^{\mu,\nu}_{\zeta,\eta} = \sum_\sigma c^\mu_{\sigma\zeta}c^\nu_{\sigma\eta}.
    \end{equation}
    We may also obtain \eqref{eq:identifying_m_and_n} by applying \Cref{thm:tensor product rule} (cf. \eqref{eq:tensor product rule}) to $H_\gamma(\zeta,\eta)$.
\end{rem}

\begin{cor} \label{cor:saturated crystal valuation on lohimod}
    For $\mu,\nu\in \cP$ and $d\in \Z_{\ge 0}$, the $\UqglPlus$-module    
    $\hilomod[\ot_{\pm}] / \soc^d(\hilomod[\ot_{\pm}])$ has a
    saturated crystal valuation $\val$ such that $L_\val/qL_\val$ has a $\Q$-basis, which forms a crystal isomorphic to
    \[
        \bigoplus_{
            \substack{
                (\mu,\nu)\ge (\zeta,\eta) \\
                |\mu| - |\zeta| = |\nu| - |\eta| = d
            }
        } \ms{B}_{\zeta,\eta}^{\oplus n^{\mu,\nu}_{\zeta,\eta}}.
    \]
\end{cor}
\pf
By \Cref{thm:isotypic decomposition of F_inftyM},
the $\UqslZero$-weight space of $\FFinftyMSubq$ with weight $\La_{\mu,\nu}$ is isomorphic to $\hilomod[\ov{\ot}_+]/\soc^d(\hilomod[\ov{\ot}_+])$ as a $\UpglPlus$-module,
and $\LCrys(\FFinftyMSubq)_{\La_{\mu,\nu}}$ is its $A_0$-submodule stable under the crystal operators $\dot{\te}_j$, $\dot{\tf}_j$ ($j\in \Z_{\ge 0}$), which are pullbacks of the upper crystal operators for $\UqglPlus$ under the map \eqref{eq:pull back Uqgln}.
Hence, we have an \emph{upper} crystal valuation (crystal valuation with upper crystal operators) of $\hilomod[\ot_+]/\soc^d(\hilomod[\ot_+])$ at $q = 0$, under \eqref{eq:pull back Uqgln}.
By \cite[Lemma 2.4.1]{Kas91}, this gives a \emph{lower} crystal valuation with the same properties.
\qed

\begin{rem}
    By \Cref{rem:R-matrix on lohimod}, \Cref{cor:saturated crystal valuation on lohimod} also holds for $\lohimod[\ot_{\pm}]$.
\end{rem}

\begin{ex}
    The crystal valuation constructed in \Cref{cor:saturated crystal valuation on lohimod} for $V_{(1),\emptyset}\ot V_{\emptyset,(1)}$ agrees with $\LCrys_\N$ in \Cref{ex:saturated crystal lattice in 10011} up to scaling.

    The inclusion $V_0(\Lambda_{(1),(1)}) \ot (V_{\emptyset, (1)} \ot_- V_{(1),\emptyset}) \subset \FFinftyM$ is given by
    \[
        \begin{tikzcd}
            { v_{\Lambda_{(1),(1)}} \ot u_{\emptyset, (1)} \ot u_{(1),\emptyset}} \arrow[r, maps to] & f_0 |0\rangle \ot 1 + (q - q^{-1})|0\rangle \ot f_0
            \end{tikzcd}
    \]
    up to scaling.
    Using the same notations as in \Cref{ex:saturated crystal lattice in 10011}, the image of elements of $V_{\emptyset, (1)} \ot_- V_{(1),\emptyset}$ are given by
    \begin{align*}
        v_{m^\vee}\ot v_n \, (m\ne n)\, \longmapsto \quad & p^2 |0\rangle \ot \cdots |0\rangle \ot |1\rangle \ot \cdots \ot | {-1}\rangle \cdots \ot |0\rangle \ot 1 \\
         &  (\text{$| {-1}\rangle$ at $m$-th position, $|1\rangle$ at $n$-th position}) \\
        v_{1^\vee}\ot v_1 \, \longmapsto \quad & f_0|0\rangle \ot 1 + (q-q^{-1}) |0\rangle \ot f_0 \\
        v_{n^\vee}\ot v_n + pv_{(n+1)^\vee}\ot v_{n+1} \,\longmapsto \quad & \left( p^2 |0\rangle^{\ot (n-1)} \ot f_0|0\rangle \ot |0\rangle + p^3 |0\rangle^{\ot n} \ot f_0|0\rangle \right) \ot 1
    \end{align*}
    Then for example, one can show that $q^{-N+2}v_{1^\vee}\ot v_1 + \cdots q v_{N^\vee}\ot v_N$ is mapped to an element of $\LCrysZero(\FFinftyM)$.
    Using this, one can compute the restriction of $\LCrysZero(\FFinftyM)$, which is constructed in the proof of \Cref{cor:main application - crystal basis}, onto $V_{\emptyset, (1)} \ot_- V_{(1),\emptyset}$.
\end{ex}

\subsection{Socle filtration of tensor product of extremal weight modules}
For $\sigma,\tau,\mu,\nu\in \cP$, let
\begin{equation*}\label{eq:transition matrix}
p_{\mu,\nu}^{\sigma,\tau}=\sum_{\la}(-1)^{|\la|}c^{\sigma}_{\la\mu}c^{\tau}_{\la'\nu},
\end{equation*}
where $\la'$ denotes the conjugate of $\la$.
For $V\in \mc{C}$, let $[V]$ denote the isomorphism class of $V$ in $K(\mathcal{C})$. 
Given $\la,\mu\in\cP$ with $\mu\subset\la$, write $[V_{\la/\mu,\emptyset}]=\sum_{\nu}c^\la_{\mu\nu}[V_{\nu,\emptyset}]$ for simplicity, and write $[V_{\emptyset,\la/\mu}]$ in a similar way.
Then we have the following character formula. 

\begin{prop}\label{prop:character formula}
For $\sigma,\tau\in \cP$, we have
\begin{equation*}
 [V_{\sigma,\tau}]=\sum_{\mu,\nu\in\cP}p_{\mu,\nu}^{\sigma,\tau}[V_{\mu,\emptyset}\ot V_{\emptyset,\nu}]=\sum_{\la\subset \sigma, \tau'}(-1)^{|\la|}[V_{\sigma/\la,\emptyset}\ot V_{\emptyset,\tau/\la'}].
\end{equation*}
\end{prop}
\pf Note that $I=\{\,[V_{\mu,\emptyset}\ot V_{\emptyset,\nu}]\,|\,\mu,\nu\in\cP\,\}$ and $S=\{\,[V_{\sigma,\tau}]\,|\,\sigma,\tau\in \cP\,\}$ are $\Z$-bases of $K(\mc{C})$. By Proposition \ref{prop:filtration on highest tensor lowest}, we have
\begin{equation}\label{eq:transition-1}
 [V_{\mu,\emptyset}\ot V_{\emptyset,\nu}]=\sum_{\sigma,\tau\in\cP}m^{\mu,\nu}_{\sigma,\tau}[V_{\sigma,\tau}].
\end{equation}
On the other hand, let $\mc{A}$ be the $\Z$-algebra generated by $h^+_r, h_s^-$ ($r,s\in \Z_{>0}$) subject to the relations $h_r^+h^-_s=h_s^-h_r^+ + h_{s-1}^-h_{r-1}^+ \cdots + h^-_{s-m}h^+_{r-m}$ with $m=\min\{r,s\}$ \cite[I.5 Example 29]{Mac}.
Note that $\mc{A}$ can be viewed as a  $\Z$-subalgebra of $\textrm{End}_{\Z}(\Lambda)$, where $\Lambda$ is the ring of symmetric functions, $h_r^-$ is the multiplication by the complete symmetric function of degree $r$ and $h_s^+$ is the adjoint of $h_s^-$ with respect to the Hall inner product on $\Lambda$.
For $\la\in \cP$, we put $s^\pm_\la={\rm det}(h^\pm_{\la_i-i+j})_{1\le i,j\le \ell(\la)}$.
Then $\{\,s^+_\mu s^-_{\nu}\,|\,\mu,\nu\in\cP\,\}$ and $\{\,s^-_\tau s^+_\sigma\,|\,\sigma,\tau\in \cP\,\}$ are $\Z$-bases of $\mc{A}$.
Then we have 
\begin{equation}\label{eq:transition-2}
\begin{split}
 s^+_\mu s^-_{\nu} &= \sum_{\sigma,\tau\in\cP}m^{\mu,\nu}_{\sigma,\tau}\, s^-_\tau s^+_\sigma,\\
 s^-_\tau s^+_\sigma &=\sum_{\mu,\nu\in\cP}n_{\mu,\nu}^{\sigma,\tau}s^+_\mu s^-_{\nu},
\end{split}
\end{equation}
for $\mu,\nu,\sigma,\tau\in\cP$ by \cite[Corollary 7.3]{K09}. 
Comparing \eqref{eq:transition-1} and the first equation in \eqref{eq:transition-2}, we see that the second equation in \eqref{eq:transition-2} also gives the transition matrix from $I$ to $S$.
\qed
\medskip

Now, we can describe the socle filtration of the tensor product of any two extremal weight modules.

\begin{thm}\label{thm:socle filtration of tensor product}
Let $\alpha,\beta,\gamma,\delta\in\cP$ be given with $|\alpha|+|\gamma|=M$ and $|\beta|+|\delta|=N$. For $0\le d\le \min\{M,N\}$, we have
\begin{equation*}
 \frac{{\rm soc}^{d+1}(V_{\alpha,\beta}\ot V_{\gamma,\delta})}{{\rm soc}^{d}(V_{\alpha,\beta}\ot V_{\gamma,\delta})}=
 \bigoplus_{\substack{\phi,\psi\in\cP \\ |\phi|=M-d,\, |\psi|=N-d}}V_{\phi,\psi}^{\ \oplus c^{(\phi,\psi)}_{(\alpha,\beta)(\gamma,\delta)}},
\end{equation*}
where 
\[
c^{(\phi,\psi)}_{(\alpha,\beta),(\gamma,\delta)}
:=
\sum_{\mu,\nu,\theta,\pi}
\left(
    \sum_{\kappa}
    c_{\kappa\,\mu}^{\alpha} \,
    c_{\kappa\,\nu}^{\gamma}
\right)
\left(
    \sum_{\epsilon}
    c_{\epsilon\,\theta}^{\beta} \,
    c_{\epsilon\,\pi}^{\delta}
\right)
c_{\mu\,\pi}^{\phi} \,
c_{\nu\,\theta}^{\psi}.
\]
\end{thm}
\pf We have 
\begin{equation*}
 U\coloneqq V_{\alpha,\beta}\ot V_{\gamma,\delta} \subset 
 V_{\alpha,\emptyset}\ot V_{\emptyset,\beta}\ot V_{\gamma,\emptyset}\ot V_{\emptyset,\delta} \cong
 V_{\alpha,\emptyset}\ot V_{\gamma,\emptyset}\ot V_{\emptyset,\beta}\ot V_{\emptyset,\delta}\eqqcolon V,
\end{equation*}
where the first inclusion follows from \Cref{cor:main application - crystal basis} and the second isomorphism follows from applying the $R$-matrix.

By \Cref{cor:main application - crystal basis}, ${\rm soc}^{d+1}(V)/{\rm soc}^d(V)$ has a direct summand $V_{\phi,\psi}$ only when $|\phi|=M-d$ and $|\psi|=N-d$ for $0\le d\le \min\{M,N\}$. Since ${\rm soc}^{d}(U)=U\cap {\rm soc}^d(V)$, the same holds for ${\rm soc}^{d+1}(U)/{\rm soc}^d(U)$.

Let $\mc{C}^\pm$ be the subcategories of $\mc{C}$ generated by highest weight modules $V_{\mu,\emptyset}$ and lowest weight modules $V_{\emptyset,\nu}$, respectively. They are semisimple and  $K(\mc{C}^\pm)$ are isomorphic to the ring of symmetric functions, say $\Lambda_{\pm}$, where $[V_{\mu,\emptyset}]$ and $[V_{\emptyset,\nu}]$ correspond to Schur functions in $\Lambda_+$ and $\Lambda_-$, respectively.

Under this identification, by \cite[Theorem 2.3]{Ko}, $\{\,[V_{\sigma,\tau}]\,|\,\sigma,\tau\in \cP\,\}$ is equal to the basis $\{\,[\sigma,\tau]_{GL}\,|\,\sigma,\tau\in \cP\,\}$ of $\Lambda_+\ot\Lambda_-$ in \cite[Definition 2.1]{Ko}. Hence the formula for multiplicity of $V_{\phi,\psi}$ in $U$ follows from \cite[Theorem 2.4]{Ko}.
\qed
\medskip

\begin{rem}
The element $[\sigma,\tau]_{GL}\in \Lambda_+\ot\Lambda_-\cong K(\mc{C})$ $(\sigma,\tau\in \cP)$ is known as a universal characters of type $GL$, which was introduced to explain the multiplicity of a rational irreducible representation in a tensor product of two irreducible rational representations of $GL_n$ and its stability for a large rank $n$. Hence our result categorifies the ring of universal character of type $GL$ by $\mc{C}$.
\end{rem}

\begin{rem}
Although our construction of saturated crystal valuations using an embedding into $\FFinftyM$ does not immediately extend to the case of $V_{\alpha,\beta}\ot V_{\gamma,\delta} / \soc^d(V_{\alpha,\beta} \ot V_{\gamma,\delta})$, we expect that it has a saturated crystal valuation.
\end{rem}

\begin{rem}
We follow the notations in the proof of Proposition \ref{prop:character formula}. Let $x=(x_1,x_2,\dots)$ and $y=(y_1,y_2,\dots)$ be two sets of formal variables.
Let $e^\pm_k$ ($k\in \Z_{>0}$) be another set of generators of $\mc{A}$ determined by $s^\pm_\la={\rm det}(e^\pm_{\la'_i-i+j})_{1\le i,j\le \ell(\la)}$, which also satisfy $e_r^+e^-_s=e_s^-e_r^+ + e_{s-1}^-e_{r-1}^+ \cdots + e^-_{s-m}e^+_{r-m}$ with $m=\min\{r,s\}$.
Put $\texttt{E}^\pm(t)=\sum_{k\geq 0}e^\pm_k t^k$.
Then we have the following identity called {\em non-symmetric Cauchy identity} \cite[I.5 Example 29]{Mac}:
\begin{equation}\label{eq:non-symmetric Cauchy}
\texttt{E}^+(y)\texttt{E}^-(x)=\texttt{E}^-(x)\texttt{E}^+(y)\dfrac{1}{\prod_{i,j\geq 1}(1-x_iy_j)},
\end{equation}
where $\texttt{E}^+(y) =\texttt{E}^+(y_1)\texttt{E}+(y_2)\texttt{E}^+(y_3)\cdots$ and $\texttt{E}^-(x)  =\texttt{E}^-(x_1)\texttt{E}^-(x_2)\texttt{E}^-(x_3)\cdots$.
It can be viewed as a non-commutative character identity associated to the decomposition in Proposition \ref{prop:decomposition of zeroth filration} (see \cite[Theorem 7.10]{K09} for a crystal-theoretic proof and also \cite{K11} for a bijective proof). 

Put $\texttt{H}^\pm(t)=\sum_{k\geq 0}h^\pm_k t^k$, $\texttt{H}^+(y)$ and $\texttt{H}^-(x)$ in a similar way. Then we have
\begin{equation}\label{eq:non-symmetric Cauchy-2}
\texttt{E}^+(y)\texttt{H}^-(x)=\texttt{H}^-(x)\texttt{E}^+(y){\prod_{i,j\geq 1}(1+x_iy_j)}.
\end{equation}
We may also have a similar representation theoretic interpretation of the above identity. In this case, we replace $\mc{F}^n_+$ in $\mc{F}^n=\mc{F}^n_+\ot \mc{F}^n_+$ with a bosonic Fock space so that $U_q(\mf{sl}_{\infty,0})$ is replaced by a $q$-boson algebra associated to a general linear Lie superalgebra of infinite rank with respect to its even part, and $V_0(\La_{\mu,\nu})$ with a $q$-deformed Kac-module. 

On the other hand, a categorification of the identities including  \eqref{eq:non-symmetric Cauchy} and \eqref{eq:non-symmetric Cauchy-2} are given in case when $x$ and $y$ are replaced by single variables $z$ and $-w$, respectively \cite{FPS}. It is given in terms of a certain non-semisimple tensor category of $\mf{sl}_\infty$-modules \cite{PSt}, which has properties similar to $\mc{C}$.
\end{rem}
\medskip

\section*{List of Notations}
\begin{longtable}{ll}
$U_q(\g)$ &: the quantized enveloping algebra (Section \ref{subsec:quantum group}) \\
$V(\la)$ &: the extremal weight module over $U_q(\g)$ of weight $\la$ (Section \ref{subsec:crystal base})\\
$(\LCrys(\la),\BCrys(\la))$ &: the crystal base of $V(\la)$ (Section \ref{subsec:crystal base})  \\
$G(\la)$ &: the global crystal basis or canonical basis of $V(\la)$ (Section \ref{subsec:crystal base}) \\
$\Uqgp$ &: the parabolic $q$-boson algebra (Section \ref{subsec:B_q}) \\
$\mc{O}^{\rm int}_{\Uqgp}$ &: the category of integrable $\Uqgp$-modules (Section \ref{subsec:integrable representations of parabolic q-boson algebras}) \\
$V_J(\la)$ &: the parabolic Verma module in $\mc{O}^{\rm int}_{\Uqgp}$ ($\la\in P_J^+$) (Sections \ref{subsec:integrable representations of parabolic q-boson algebras})\\
$(\LCrys_J(\la),\BCrys_J(\la))$ &: the crystal base of $V_J(\la)$ (Section \ref{subsec:crystal base of parabolic Verma modules}) \\
$U_q(\gl_S)$ &: the quantized enveloping algebra for an interval $S$ (Section \ref{sec:gln}) \\
$\cP$ &: the set of partitions \\
$\Z_+^n$ &: the set of generalized partitions of length $n$ \\
$\e_\la$ &: a weight associated to $\la\in \cP$ or $\Z_+^n$ (Section \ref{sec:gln}) \\
$\e^n_{\mu,\nu}$ &: a weight associated to $\mu,\nu\in \cP$ with $\ell(\mu)+\ell(\nu)\le n$ (Section \ref{sec:gln}) \\
$V_{\mu,\nu}$ &: the extremal weight module of $U_q(\gl_{>0})$ ($(\mu,\nu)\in \cP^2$) (Section \ref{subsec:ext wt module}) \\
$(\ms{L}_{\mu,\nu},\ms{B}_{\mu,\nu})$ &: the crystal base of $V_{\mu,\nu}$ (Section \ref{subsec:monoidal cat C}) \\
$\succeq, \ge$ &: partial orders on $\cP^2$ (Section \ref{subsec:monoidal cat C})\\
$\mc{C}$ &: the monoidal category generated by $V_{\mu,\nu}$'s (Section \ref{subsec:monoidal cat C})\\
$U_p(\dot{\g})$ &: the quantized enveloping algebra  over $\Q(p)$ (Section \ref{subsec:q-deformed exterior and symmetric algebras})\\
$\mc{A}(V\ot {W})$ &: a $U_q(\g)\ot U_p(\dot{\g})$-module (Section \ref{subsec:q-deformed exterior and symmetric algebras}) \\
$\tbgwed_{S,T}$ &: the $q$-deformed exterior algebra  (Section \ref{subsec:q-ext alg of type A}) \\
$B_{S,T}$ &: the set of $S\times T$ matrices of binary entries (Section \ref{subsec:q-ext alg of type A})\\
$w_M$ &: a monomial in $\tbgwed_{S,T}$ ($M\in B_{S,T}$) (Section \ref{subsec:q-ext alg of type A}) \\
$\FF_{S,T}$ &: the limit of $\tbgwed_{S',T}$ (Section \ref{sec:semi-infinite limit})\\
$\FFn$ &: $\FF_{\mathbb{R},[1,n]}$, the Fock space of level $n$ (Section \ref{sec:semi-infinite limit}) \\
$(\ms{L}(\FFn),\ms{B}(\FFn))$ &: the crystal base of $\FFn$ (Section \ref{sec:semi-infinite limit})\\
$M(\la)=M^n(\mu,\nu)$ &: an element in $\ms{B}(\FFn)$ associated to $\la\in \Z_+^n$ (Section \ref{sec:semi-infinite limit})\\
$\Lambda_\la=\Lambda^n_{\mu,\nu}$ &: a weight associated to $\la\in \Z_+^n$ (Section \ref{sec:semi-infinite limit}) \\
$\UqslZero$ &: a parabolic $q$-boson algebra (Section \ref{sec:uqslzero and mm})\\
$V_0(\La)$ &: an integrable irreducible $\UqslZero$-module (Section \ref{sec:uqslzero and mm})\\ 
$(\LCrys_0(\La), \BCrys_0(\La))$ &: the crystal base of $V_0(\La)$  (Section \ref{sec:uqslzero and mm})\\
$\mc{M}$ &: $V_0(\La)$ with $\Lambda=0$ (Section \ref{sec:uqslzero and mm})\\
$\phi_{n,n+1}$ &: a $\UqslZero$-linear map from $\FFnM$ to $\FFnplusM$ (Section \ref{subsec:limit of Fock space})\\
$\FFinftyM$ &: the infinite-level Fock space, a limit of $\FFnM$ (Section \ref{subsec:limit of Fock space})\\
$(\FFinftyM)_0$ &: a subspace of $\FFinftyM$ (Section \ref{subsec:limit of Fock space})\\
$H^n_\gamma(\mu,\nu)$ &: a subset of $\BCrysZero(\MM)_\gamma$ (Section \ref{subsec: filtration of F})\\ 
$(\FFnM)_{\ge \gamma}$ &: a subpace of $\FFnM$ ($\gamma \in \wt(\MM)$)(Section \ref{subsec: filtration of F})\\
$(\FFnM)^{(\zeta,\eta)}_{(\mu,\nu)}$ &: an isotypic component of $\FFnM$ (Section \ref{subsec: filtration of F}) (see also \eqref{eq:def of FFnM filtrations})\\
$(\FFnM)_{\ge -d}$ &: a filtration of $\FFnM$ (Section \ref{subsec: filtration of F})\\
$(\FFinftyM)_{\ge \gamma}$ &: a limit of $(\FFnM)_{\ge \gamma}$ (Section \ref{subsec: filtration of F})\\
$(\FFinftyM)_{\ge -d}$ &: a limit of $(\FFnM)_{\ge -d}$ (Section \ref{subsec: filtration of F})\\
$\texttt{U}$  &: $U_q(\g)$ or $U_q(\g,\mf{p})$ (Section \ref{subsec:crystal valuation})\\
$\val$ &: a valuation on an integrable $\texttt{U}$-module $V$  (Section \ref{subsec:crystal valuation})\\
$L_{\val}$ &: an $A_0$-submodule of $V$ associated to $\val$ (Section \ref{subsec:crystal valuation})\\
$\val_{\ms{L}(\mc{F}^n)}$ &: the crystal valuation associated to $\LCrys(\FFn)$ (Section \ref{subsec:lemmas on isotypic components})\\
$(\FFnM)_{\gamma, (\mu,\nu)}$ &: a subspace of $\FFnM$ (Section \ref{subsec: semisimple subquotient})\\
$\phibar_{n,n+1}$ &: a map on $(\FFnM)_{\gamma, (\mu,\nu)}/(\FFnM)_{>\gamma}$ (Section \ref{subsec: semisimple subquotient})\\
$\val^n$ &: the crystal valuation associated to $\LCrysZero(\FFnM)$ (Section \ref{subsec: semisimple subquotient})\\
$(\FFinftyM)^{(\zeta,\eta)}_{\ge (\mu,\nu)}$ &: a limit of $(\FFnM)^{(\zeta,\eta)}_{\ge (\mu,\nu)}$ (Section \ref{subsec: semisimple subquotient})\\
$(\FFnM)_{-d}$ &: a summand of $\FFnM$ (Section \ref{subsec: a crystal valuation on socle quotients})\\
$\pi^n_{-d}$ &: the canonical projection $\FFnM \to (\FFnM)_{-d}$ (Section \ref{subsec: a crystal valuation on socle quotients})\\
$\val^n_{-d}$ &: a map $\val^n_{-d} := \val^n \circ \pi^n_{-d}$ (Section \ref{subsec: a crystal valuation on socle quotients})\\
$\val^\infty_{-d}$ &: a limit given by $\val^\infty_{-d} (x) := \lim_{n\to \infty} \big( \val^n_{-d}(x_n) - dn \big)$ (Section \ref{subsec: a crystal valuation on socle quotients})\\
$\FFinftyMSubq$, $\FFinftyMSubqSoc$ &: a subspace and a subquotient of $\FFinftyM$ (Section \ref{subsec:Socle filtration})\\
$\LCrys(\FFinftyMSubq)$, $\LCrys(\FFinftyMSubq^\circ)$ &: $A_0$-submodules associated to $\val^\infty_{-d}$ (Section \ref{subsec:Socle filtration})\\
\end{longtable}

\subsection*{Acknowledgements}
The authors would like to thank anonymous referees for very careful reading of the manuscript and many helpful comments.

\appendix

{\small
\ifdeclaration
\section*{Declarations}
\subsection*{Funding and/or Conflicts of interests/Competing}
This work is supported by the National Research Foundation of Korea(NRF) grant funded by the Korea government(MSIT) (No.\discretionary{2020R1}{A5A1016126}{2020R1A5A1016126} and RS-2024-00342349). The authors declare no competing interests.
\fi

}

\end{document}